\documentclass[a4paper,10pt]{article}
\usepackage{graphicx}        % standard LaTeX graphics tool
\usepackage{multicol}        % used for the two-column index
\usepackage[bottom]{footmisc}% places footnotes at page bottom
\usepackage{epsfig}
\usepackage{amsfonts,amsthm}
\usepackage{authblk}
\usepackage{amssymb}
\usepackage{graphicx}
\usepackage{float}
\usepackage{diagbox}
\usepackage[applemac]{inputenc}
\usepackage[T1]{fontenc}
\usepackage{blkarray}
\usepackage{color}
\usepackage{amsmath,bm}
\usepackage{mathtools}
\usepackage{mathabx}
\usepackage{nicefrac}
\usepackage[hmargin={15mm,15mm},vmargin={17mm,22mm}]{geometry}
\usepackage{amssymb}
\usepackage{array}
\usepackage{pdfpages}
\usepackage{subfig}
\usepackage{multirow}
\usepackage{enumitem}
\usepackage{bbold}
\usepackage{xparse}								% for \DeclareDocumentCommand
\usepackage{stmaryrd}							% for symbols $\llbracket a+b\rrbracket$
\usepackage{sidecap}								% for SCfigure
\usepackage{upgreek}
\usepackage{accents}
\allowdisplaybreaks
\usepackage{tcolorbox}
\usepackage{algorithmic}
\usepackage[ocgcolorlinks,allcolors={blue}]{hyperref}

\usepackage[bbgreekl]{mathbbol}

\DeclareSymbolFontAlphabet{\mathbb}{AMSb}
\DeclareSymbolFontAlphabet{\mathbbl}{bbold}

\usepackage{bbm}

\def\N{\mathbb N}
\def\R{\mathbb R}

\def\P{\mathbb P}

\def\D{\mathcal{D}}
\def\S{\mathcal{S}}

\def\<{\langle}
\def\>{\rangle}

\def\lsim{\lesssim}

\def\Chi{\raise .3ex \hbox{\large $\chi$}}

\def\[{\Bigl [}
\def\]{\Bigr ]}
\def\({\Bigl (}
\def\){\Bigr )}
\def\l{\iota}

\def\dsp{\displaystyle}
\def\x{{\bf x}}

\def\n{{\bf n}}

\def\U{{\bf U}}

\def\q{{\bf q}}

\def\g{{\mathbf{g}}}
\def\aa{{\mathfrak a}}

\def\x{{\bf x}}
\def\dsp{{\displaystyle x}}
\def\d{{\rm d}}
\def\div{{\rm div}}
\def\a{\alpha}
\def\q{\mathbf{q}}
\def\K{\mathbb{K}}
\def\n{\mathbf{n}}
\def\dsp{\displaystyle}
\def\bu{\mathbf{u}}
\def\bv{\mathbf{v}}

\def\O{\Omega}
\def\G{\Gamma}

\def\bbsig{\bbsigma}
\def\bbeps{\bbespilon}
\newcommand{\jump}[1]{\llbracket #1 \rrbracket}

\newcommand{\dtn}{\delta t^{n+\frac{1}{2}}}
\newcommand{\dtnmun}{\delta t^{n-\frac{1}{2}}}

\newcommand{\dtzero}{\delta t^{\frac{1}{2}}}

\newcommand{\nf}{\nicefrac}

\newcommand{\email}[1]{\href{mailto:#1}{#1}}

\begingroup
    \makeatletter
    \@for\theoremstyle:=definition,remark,plain\do{%
        \expandafter\g@addto@macro\csname th@\theoremstyle\endcsname{%
            \addtolength\thm@preskip\parskip
            }%
        }
\endgroup

\newtheorem{theorem}{Theorem}
\numberwithin{theorem}{section}
\newtheorem{proposition}[theorem]{Proposition}
\newtheorem{lemma}[theorem]{Lemma}

\theoremstyle{remark}
\newtheorem{remark}[theorem]{Remark}
\theoremstyle{definition}

\newtheorem{definition}[theorem]{Definition}

\newenvironment{acknowledgements}
      {\bigskip\bigskip~\newline\textbf{Acknowledgements}}

\def\g{{\rm nw}}
\def\l{{\rm w}}
\def\rt{{\rm rt}}

\def\trace{\gamma}
\def\grad{\nabla}
\def\del{\partial}
\def\div{{\rm div}}

\def\weakto{\rightharpoonup}

\definecolor{labelkey}{rgb}{0.6,0,1}

\usepackage[normalem]{ulem}
 \newcounter{corr}
 \definecolor{violet}{rgb}{0.580,0.,0.827}
 \newcommand{\corr}[3]{\typeout{Warning : a correction remains in page
 \thepage}
 				\stepcounter{corr}        
 				{\color{blue}\ifmmode\text{\,\sout{\ensuremath{#1}}\,}\else\sout{#1}\fi}
         {\color{red}#2}
         {\color{violet} \ifmmode\text{#3}\else #3\fi }}

\usepackage{parskip}
\usepackage{algorithmic}

\usepackage{diagbox}
\usepackage{booktabs}

\makeatletter
\def\thm@space@setup{%
  \thm@preskip=\parskip \thm@postskip=0pt
}
\makeatother

\begin{document}
\title{Gradient discretization of two-phase poro-mechanical models with discontinuous pressures at matrix fracture interfaces %\footnote{This work has been supported by SIR Research Grant no.~RBSI14VTOS funded by MIUR -- Italian Ministry of Education, Universities, and Research, and by ``National Group of Scientific Computing'' (GNCS-INdAM).}
}
\author[1]{{Francesco Bonaldi}\footnote{Corresponding author, \email{francesco.bonaldi@univ-cotedazur.fr}}}
\author[1]{{Konstantin Brenner}\footnote{\email{konstantin.brenner@univ-cotedazur.fr}}}
\author[2]{{J\'er\^ome Droniou}\footnote{\email{jerome.droniou@monash.edu}}}
\author[1]{{Roland Masson}\footnote{\email{roland.masson@univ-cotedazur.fr}}}
\author[3]{{\\ Antoine Pasteau}\footnote{\email{Antoine.Pasteau@andra.fr}}}
\author[3]{{Laurent Trenty}\footnote{\email{Laurent.Trenty@andra.fr}}}
\affil[1]{Universit\'e C\^ote d'Azur, Inria, CNRS, Laboratoire J.A. Dieudonn\'e, team Coffee, France}%
\affil[2]{School of Mathematics, Monash University, Victoria 3800, Australia}%
\affil[3]{Andra, Chatenay-Malabry, France}%
\date{}
\maketitle
\begin{abstract}
\noindent
We consider a two-phase Darcy flow in a fractured and deformable porous medium for which the fractures are described as a network of planar surfaces leading to so-called hybrid-dimensional models. The fractures are assumed open and filled by the fluids and small deformations with a linear elastic constitutive law are considered in the matrix.
As opposed to \cite{bonaldi:hal-02549111}, the phase pressures are not assumed continuous at matrix fracture interfaces, which raises new challenges in the convergence analysis related to the additional interfacial equations and unknowns for the flow.
As shown in \cite{BHMS2018,gem.aghili}, unlike single-phase flow, discontinuous pressure models for two-phase flows provide a better accuracy than continuous pressure models even for highly permeable fractures. This is due to the fact that fractures fully filled by one phase can act as barriers for the other phase, resulting in a pressure discontinuity at the matrix fracture interface. 

The model is discretized using the gradient discretization method \cite{gdm}, which covers a large class of conforming and non conforming schemes. This framework allows for a generic convergence analysis of the coupled model using a combination of discrete functional tools. In this work, the gradient discretization of \cite{bonaldi:hal-02549111} is extended to the discontinuous pressure model and the convergence to a weak solution is proved.  Numerical solutions provided by the continuous and discontinuous pressure models are compared on gas injection and suction test cases using a Two-Point Flux Approximation (TPFA) finite volume scheme for the flows and $\P_2$ finite elements for the mechanics.  
\end{abstract}
%
%%-----------------------------
%%      your text
%%-----------------------------

\section{Introduction}

Coupled flow and geomechanics play an important role in many subsurface processes such as water management, geothermal energy, CO$_2$ sequestration, oil and gas production or nuclear waste storage. This is particularly the case in the presence of fractures which have a strong impact both on the flow and on the rock mechanical behavior.
This work considers the so called hybrid-dimensional or Discrete Fracture Matrix (DFM) models representing the fractures as a network of co-dimension one surfaces coupled with the surrounding matrix domain. The reduced flow model is then obtained by averaging both the unknowns and the equations in the fracture width and by imposing appropriate transmission conditions at both sides of the matrix fracture interfaces. The mechanical model is set on the matrix domain with appropriate boundary conditions on both sides of the fracture interfaces. 
This type of hybrid-dimensional models has been the object of intensive researches over the last twenty years due to the ubiquity of fractures in geology and their considerable impact on the flow and transport of mass and energy in porous media, and on the mechanical behavior of rocks. For the derivation and analysis of such models, let us refer to \cite{MAE02,FNFM03,KDA04,MJE05,ABH09,BGGLM16,BHMS2016,NBFK2019} for single-phase Darcy flows,  \cite{BMTA03,RJBH06,MF07,Jaffre11,BGGM2015,DHM16,BHMS2018,gem.aghili} for two-phase Darcy flows, and  \cite{MK2013,KTJ2013,JJ14,GWGM2015,hanowski2016simulation,GKT16,JZ17,GFSZ17,UKBN2018,Girault2019} for poroelastic models. \\

As in \cite{bonaldi:hal-02549111}, this work focuses on a hybrid-dimensional two-phase Darcy flow model coupled with a linear poroelastic deformation of the matrix. The fractures are assumed to remain open and fully filled by the fluids, {and their propagation over time is neglected}. The Poiseuille law is used for the tangential velocity along the fracture network and extended to two-phase flow based on generalized Darcy laws. As in \cite{coussy}, the concept of equivalent pressure, used to extend  the poro-mechanical coupling to two-phase flow, is based  on the capillary energy density. This is a crucial choice to obtain the stability of the coupled model. 

In \cite{bonaldi:hal-02549111}, the continuity of both phase pressures is used as a transmission condition at matrix fracture interfaces. This is a classical assumption in the case of highly permeable fractures such as open fractures. As shown e.g. in \cite{hennicker:hal-02437030}, this choice is fully justified for the case of single-phase flows. On the other hand, in the case of two-phase flow, this assumption can lead to inaccurate solutions at the matrix fracture interfaces \cite{BHMS2018,gem.aghili}. This is in particular the case when the fractures are fully filled by one phase and act as barriers for the other phase due to its very low relative permeability within the fractures, hence leading to a pressure discontinuity. Let us refer to \cite{gem.aghili} for striking examples including the desaturation by suction at the interface between the atmosphere and a low permeable and fractured storage rock.

This potential inaccuracy of continuous pressure models motivates us to consider the extension of the analysis carried out in \cite{bonaldi:hal-02549111} to  hybrid-dimensional discontinuous pressure flow models \cite{Jaffre11,DHM16,BHMS2018,gem.aghili}. For such flow models, the Darcy fluxes between the matrix fracture interface and the fracture are modelled using a two-point flux approximation combined with an upwind approximation of the mobilities \cite{DHM16,BHMS2018,gem.aghili}. Following \cite{DHM16}, the model also includes a layer of damaged rock at matrix fracture interfaces. This additional accumulation term plays a major role in the numerical analysis of the model and also improves the nonlinear convergence at each time step of the simulation \cite{DHM16,quenjel:hal-02957054}. It must be kept sufficiently small to maintain the accuracy of the solution (see \cite{DHM16}). 
Following \cite{bonaldi:hal-02549111} and \cite{DHM16}, this new hybrid-dimensional poro-mechanical model is discretized using the gradient discretization method \cite{gdm}.
This framework is based on abstract vector spaces of discrete unknowns combined with reconstruction operators.
The gradient scheme is then obtained by substitution of the continuous operators by their discrete counterparts in the weak formulation of the coupled model. The main asset of this framework is to allow a generic convergence analysis based on general properties of the reconstruction operators that hold for a large class of conforming and non conforming discretizations.
{Let us point out that, with respect to~\cite{bonaldi:hal-02549111}, additional trace and jump operators need to be defined in this framework, along with new definitions of coercivity, consistency, limit-conformity, and compactness.}
The two main ingredients to discretize the coupled model are the discretizations of the hybrid-dimensional discontinuous pressure two-phase Darcy flow  and the discretization of the mechanics. Let us briefly mention, in both cases, a few families of discretizations typically satisfying the gradient discretization properties. 

For the discretization of the Darcy flow, the gradient discretization framework covers the case of cell-centered finite volume schemes with Two-Point Flux Approximation on strongly admissible meshes \cite{KDA04,ABH09,gem.aghili}, or some symmetric Multi-Point Flux Approximations \cite{TFGCH12,SBN12,AELHP153D} on tetrahedral or hexahedral meshes. It also accounts for the families of Mixed Hybrid Mimetic and Mixed or Mixed Hybrid Finite Element discretizations such as in \cite{MJE05,BHMS2016,AFSVV16,Girault2019}, and for vertex-based discretizations such as  the Vertex Approximate Gradient scheme \cite{BHMS2016,DHM16,BHMS2018}.
For the discretization of the elastic mechanical model, the gradient discretization framework covers conforming finite element methods such as in \cite{GWGM2015}, the Crouzeix-Raviart discretization ~\cite{hansbo-larson,dipietro-lemaire}, the Hybrid High Order discretization ~\cite{dipietro-ern}, and the Virtual Element Method~\cite{beirao-brezzi-marini}.

The main objective of this work is to introduce the gradient discretization of the hybrid-dimensional poro-mechanical model with discontinuous pressure at matrix fracture interfaces. Then, we prove the convergence of the discrete solution to a weak solution of the model.
Compared with  \cite{bonaldi:hal-02549111}, new difficulties arise from the interfacial additional nonlinear flux and accumulation terms including the damaged rock type. Assuming that the fracture normal transmissivity is fixed in the interfacial two-point fluxes, i.e. that its fracture aperture dependence is frozen, we are able to prove the convergence of the gradient scheme solution to a weak solution of the model. This assumption is rather mild since, in practice, the solution depends only weakly on this fracture normal transmissivity as long as it remains much larger than the matrix transmissivity. 
{Concerning compactness estimates, the same techniques as in~\cite{bonaldi:hal-02549111} are used: time translates, uniform-in-time $L^2$-weak estimates, and a discrete version of the Ascoli-Arzel\`a theorem. In~\cite{bonaldi:hal-02549111}, where fields defined in matrix and fracture are related, matrix and fracture contributions have to be separated by a cut-off argument (since the fracture width vanishes at tips). On the other hand, in this work, such a separation stems from the model itself, but the damaged rock layer has to be embedded in the time translates of the saturations, by using ad-hoc test functions combining the matrix and damaged layer rock types.}

As in \cite{bonaldi:hal-02549111}, the proof additionally assumes that the matrix porosity remains bounded from below  by a strictly positive constant, that the fracture aperture remains larger than a fixed aperture vanishing only at the tips, and that the mobility functions are bounded from below by strictly positive constants. The assumptions on the porosity and fracture aperture cannot be avoided since the continuous model does not ensure these properties, which are needed to ensure its well-posedness. The assumption on the mobilities are classical to carry out the stability and convergence analysis of two-phase Darcy flows with heterogeneous rock types (see \cite{EGHM13,BGGM2015,DHM16}). 

The second objective of this work is to compare the discontinuous pressure poro-mechanical model investigated in this work to the continuous pressure poro-mechanical model presented in 
\cite{bonaldi:hal-02549111}. Two test cases are considered. As in \cite{bonaldi:hal-02549111}, the first test case simulates the gas injection in a cross-shaped fracture network immersed in an initially water saturated porous medium. The second test case models the desaturation of a low permeable medium by suction at the interface with a ventilation tunnel. The data set of this second test case is based on the Callovo-Oxfordian argilite rock properties of the nuclear waste storage prototype facility of Andra. The geometry uses an axisymmetric DFM model based on a simplified version of the fracture network at the interface between the storage rock and the ventilation tunnel. In both cases the discretization is based on the Two-Point Flux Approximation finite volume scheme for the flows and second-order finite elements for the mechanical deformation. \\

The rest of the article is organized as follows. Section~\ref{sec:modeleCont} introduces the continuous hybrid-dimensional coupled model with discontinuous pressures at matrix fracture interfaces. Section~\ref{sec:gradientscheme} describes the gradient discretization method for the coupled model including the definition of the reconstruction operators, the discrete variational formulation and the properties of the gradient discretization needed for the subsequent convergence analysis.
Section~\ref{sec:convergence} proceeds with the convergence analysis. The a priori estimates are established in Subsection~\ref{subsec:apriori}, the compactness properties in Subsection~\ref{subsec:compactness} and the convergence to a weak solution is proved in Subsection~\ref{subsec:convergence}. This convergence falls short, in general, from identifying the limit matrix--fracture nonlinear fluxes; this issue is discussed in Subsection~\ref{identification}, in which an assumption is given on the limit fracture width under which the fluxes can be fully identified.
In Section~\ref{sec:numerical.example}, devoted to numerical experiments, the discontinuous pressure model is compared to the continuous pressure model presented in \cite{bonaldi:hal-02549111}.

\section{Continuous model}\label{sec:modeleCont}

\begin{figure}[h!]
\begin{center}
\includegraphics[scale=.55]{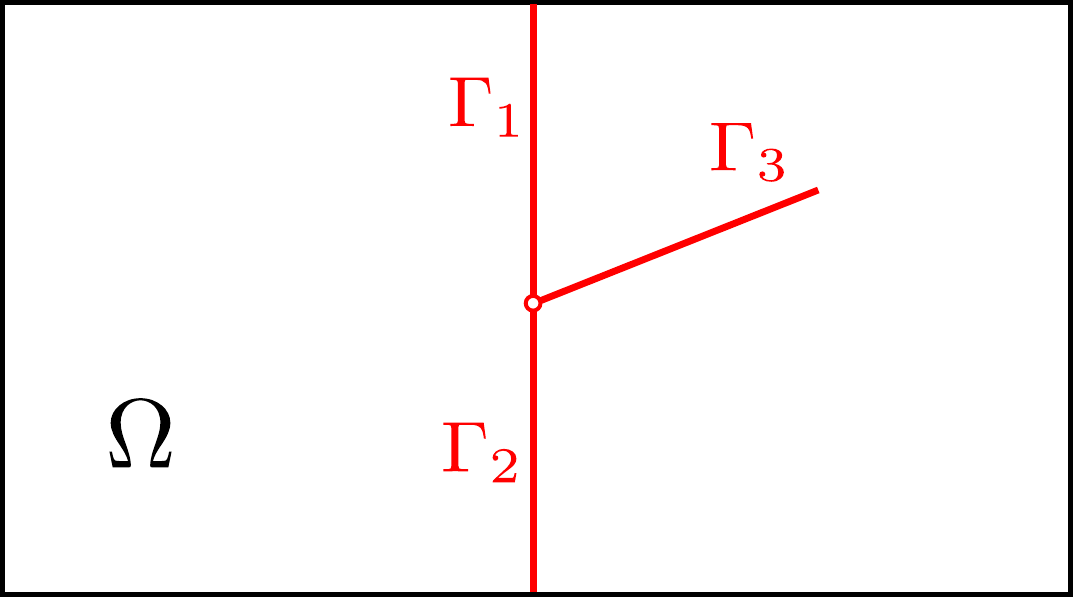}   
\caption{Example of a 2D domain $\Omega$ with three intersecting fractures $\Gamma_i$, $i=1,2,3$.}
\label{fig_network}
\end{center}
\end{figure}

We consider a bounded polytopal domain $\Omega$ of $\R^d$, $d\in\{2,3\}$, partitioned
into a fracture domain $\Gamma$ and a matrix domain $\O\backslash\overline\G$.
The network of fractures is defined by 
$$
\overline \Gamma = \bigcup_{i\in I} \overline \Gamma_i
$$  
where each fracture $\Gamma_i\subset \Omega$, $i\in I$ is a planar polygonal simply connected open domain. Without restriction of generality, we will assume that the fractures may intersect exclusively at their boundaries (see Figure \ref{fig_network}), that is, for any $i,j \in I, i\neq j$ one has $\Gamma_i\cap \Gamma_j = \emptyset$, but not necessarily $\overline{\Gamma}_i\cap \overline{\Gamma}_j = \emptyset$.

The two sides of a given fracture of $\Gamma$ are denoted by $\pm$ in the matrix domain, with unit normal vectors $\n^\pm$ oriented outward of the sides $\pm$. We denote by $\gamma_\aa$ the trace operators on the sides $\aa =\pm$ of $\Gamma$ for functions in $H^1(\Omega\setminus\overline\Gamma)$, by $\gamma_{\del\Omega}$ the trace operator for the same functions on $\del\O$, and by $\jump{\cdot}$ the normal trace jump operator on $\Gamma$ for functions in $H_\div(\O\backslash\overline\Gamma)$, defined by
$$
\jump{\bar\bu} = \bar\bu^+ \cdot \n^+ + \bar\bu^- \cdot \n^-  \, \mbox{ for all } \, \bar\bu \in H_\div(\O\backslash\overline\Gamma). 
$$
We denote by $\nabla_\tau$ the tangential gradient and by $\div_\tau$ the tangential divergence on the fracture network $\Gamma$. The symmetric gradient operator $\bbeps$ is defined such that $\bbeps(\bar\bv) = {1\over 2} (\nabla \bar\bv +^t\!(\nabla \bar\bv))$ for a given vector field $\bar\bv\in H^1(\O\backslash\overline\Gamma)^d$.

Let us fix a continuous function $d_0: \Gamma \to (0,+\infty)$  with zero limits at
$\partial \Gamma \setminus (\partial\Gamma\cap \partial\Omega)$ (i.e.~the tips of $\Gamma$) and strictly positive limits at $\partial\Gamma\cap \partial\Omega$.
The fracture aperture, denoted by $\bar d_f$ and such that $\bar d_f = - \jump{\bar \bu}$ for a displacement field
$\bar\bu \in H^1(\O\backslash\overline\Gamma)^d$, will be assumed to satisfy the following open fracture condition 
$$
\bar d_f(\x) \geq d_0(\x) \mbox{ for a.e.\ } \x \in \Gamma.    
$$
Let us introduce some relevant function spaces. First, we denote by $H_{d_0}^1(\Gamma)$ the space made of functions $v_\Gamma$ in $L^2(\Gamma)$, such that $d_0^{\nf 3 2} \nabla_\tau v_\Gamma$ belongs to  $L^2(\Gamma)^{d-1}$, 
and whose traces are continuous at fracture intersections $\del\G_i\cap\del\G_j$, $(i,j)\in I\times I$ ($i\neq j$) and vanish on the boundary $\partial \Gamma\cap \partial\Omega$.
We then introduce the space
\begin{equation}\label{displacement}
\U^0 =\{ \bar \bv\in (H^1(\O\backslash\overline\Gamma))^d \mid \trace_{\del\O} \bar \bv = 0\}
\end{equation}
for the displacement vector, and
\begin{equation}\label{pressures}
V^0 = V^0_m \times V^0_f,%  \{\bar v\in H^1_0(\O) \mid \gamma \bar v \in H^1_{d_0}(\G)\}
\end{equation}
where
$$
V^0_m = \{ \bar v\in H^1(\Omega\setminus\overline\Gamma) \,|\, \trace_{\del\O} \bar v = 0\},
$$
for each matrix phase pressure, and 
$$
V^0_f = H_{d_0}^1(\Gamma),
$$
for each fracture phase pressure. For $\bar v =(\bar v_m,\bar v_f) \in V^0$, let us denote by
$$
\jump{\bar v}_\aa = \gamma_\aa \bar v_m - \bar v_f, 
$$
the jump operator on the side $\aa=\pm$ of the fractures. 

The matrix, fracture and damaged rock types are denoted by the indices ${\rm rt} = m$, ${\rm rt} = f$, and ${\rm rt} = \pm$, {respectively}, and the non-wetting and wetting phases by the superscripts $\alpha =\g$ and $\alpha=\l$, {respectively}. Finally, for any $x\in\mathbb R$, we set $x^+ = \max\{0,x\}$ and $x^- = -(-x)^+$.

\begin{SCfigure}
\includegraphics[width=6cm]{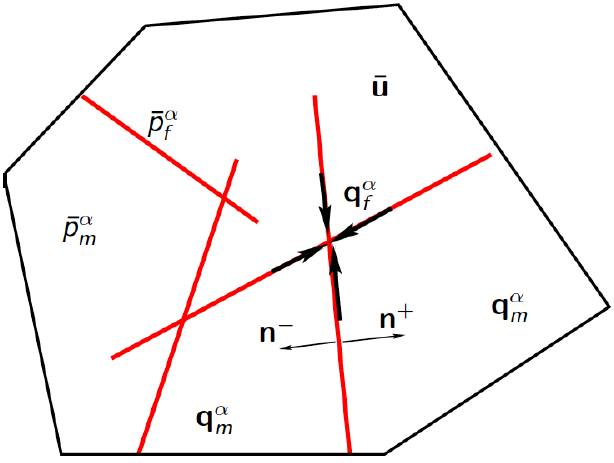}
\caption{Example of a 2D domain $\Omega$  with its fracture network $\Gamma$, the unit normal vectors $\n^\pm$ to $\Gamma$, the phase pressures $\bar p^\alpha_m$ in the matrix and $\bar p^\alpha_f$ in the fracture network, the displacement vector {field} $\bar \bu$, the matrix Darcy velocities $\q^\alpha_m$ and the fracture tangential Darcy velocities $\q^\alpha_f$ integrated along the fracture width. }
\end{SCfigure}
\ \
\begin{SCfigure}
\includegraphics[scale=.7]{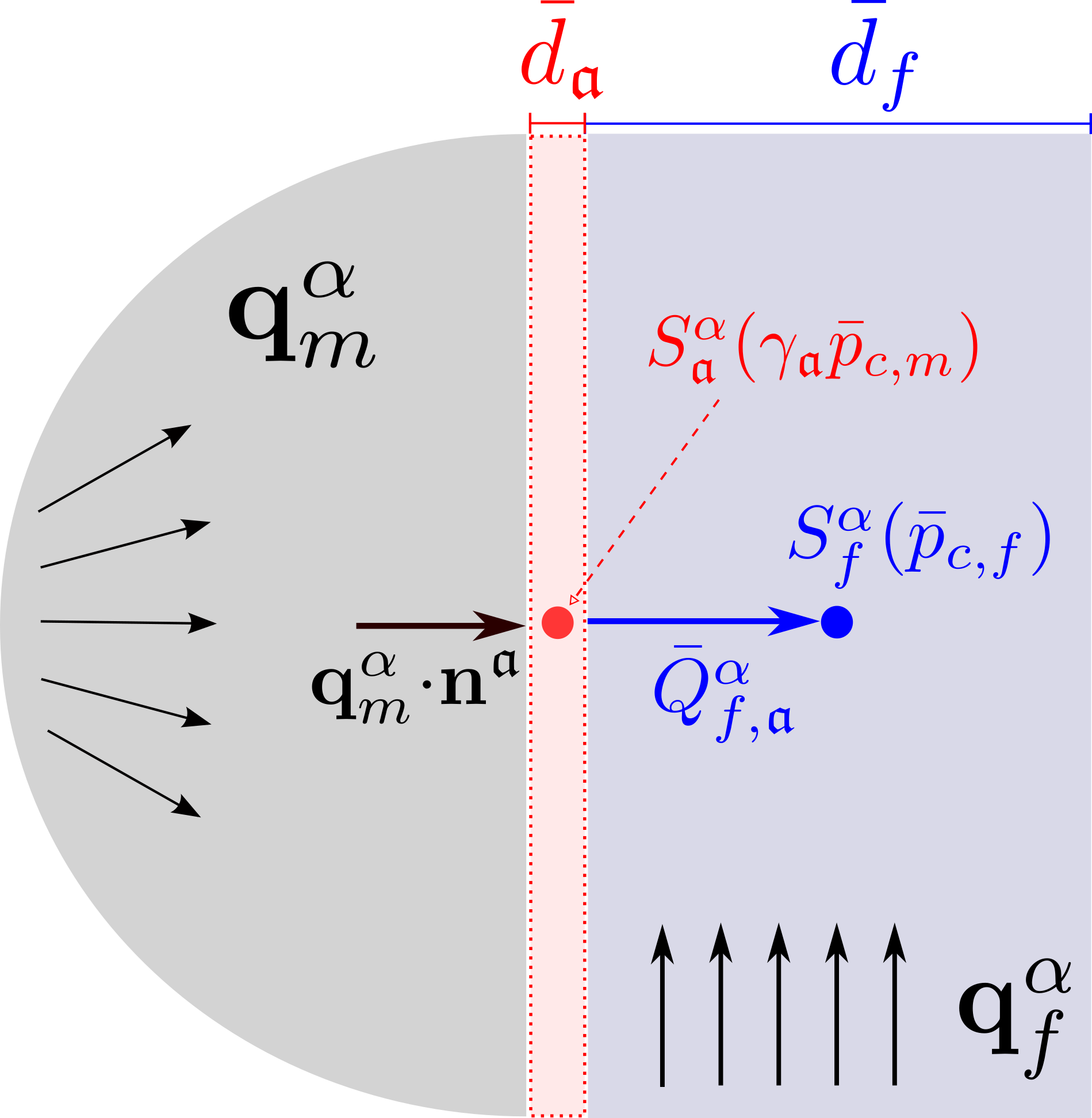}
\caption{{Illustration of the flux transmission condition between matrix and fracture, including a layer of damaged rock of thickness $d_\aa$, $\aa =~\pm$. It can be seen as an upwind two-point-like approximation of $\bar Q^\a_{f ,\aa}$. The arrows show the positive orientation of the normal fluxes $\q_m^\a\cdot\n^\aa$ (inward with respect to the damaged layer) and $\bar Q^\a_{f,\aa}$ (outward with respect to the damaged layer)}.}
\end{SCfigure}

The PDEs model {reads}: find the phase pressures $\bar p^\alpha_\nu$, $\nu\in\{m,f\}$, $\alpha\in\{\g,\l\}$, and the displacement vector {field} $\bar \bu$, such that $\bar p_{c,\nu} = \bar p_\nu^\g - \bar p_\nu^\l$,  and for $\alpha\in\{\g,\l\}$, 
\begin{equation}
\label{eq_edp_hydromeca} 
\left\{
\begin{array}{lll}
&\partial_t\left(  \bar\phi_m S^{\a}_m(\bar p_{c,m}) \right) + \div \left(  \q^\a_m \right) = h_m^\a & \mbox{ on } (0,T)\times \Omega\setminus\overline\Gamma,\\[1ex]
& \q^\a_m = \dsp - \eta^\a_m ( S^{\a}_m(\bar p_{c,m})) \K_m  \nabla \bar p^\a_m & \mbox{ on } (0,T)\times \Omega\setminus\overline\Gamma,\\[1ex]
&\partial_t\left(  \bar d_f S^{\a}_f(\bar p_{c,f}) \right)
+ \div_\tau (  \q^\a_f ) -   Q^\a_{f,+} -   Q^\a_{f,-}  = h_f^\a & \mbox{ on } (0,T)\times \Gamma,\\[1ex]
& \q^\a_f = \dsp -  \eta^\a_f (S^{\a}_f(\bar p_{c,f})) ({1\over 12} \bar d_f^3) \nabla_\tau \bar p_f^\a  & \mbox{ on } (0,T)\times \Gamma,\\[1ex]
& -\div \(\bbsig(\bar\bu) - b ~ \bar p^E_m{\mathbb I}\)= \mathbf{f} & \mbox{ on } (0,T)\times \Omega\setminus\overline\Gamma\\[1ex]
& \bbsig(\bar\bu) =  2\mu ~\bbeps(\bar\bu) + \lambda ~\div(\bar\bu) ~\mathbb{I} & \mbox{ on } (0,T)\times \Omega\setminus\overline\Gamma,  
\end{array}
\right.
\end{equation}
with the coupling conditions 
\begin{equation}
\label{closure_laws}
\left\{
\begin{array}{lll}
  & \partial_t \bar\phi_m = \dsp b~\div \partial_t \bar\bu + \frac{1}{M} \partial_t \bar p^E_m & \mbox{ on } (0,T)\times \Omega\setminus\overline\Gamma,\\  [2ex]
  & \q_m^\a\cdot\n^\aa - Q_{f,\aa}^\a = \bar d_\aa \bar \phi_\aa \partial_t S^\a_\aa(\gamma_\aa \bar p_{c,m}) &   \mbox{ on } (0,T)\times \Gamma, \aa=\pm,\\[1ex]
  & \bar Q_{f,\aa}^\a =  \eta_\aa(S^\a_\aa(\gamma_\aa \bar p_{c,m})) T_f \jump{\bar p^\a}_\aa^+  -
  \eta_f(S^\a_f(\bar p_{c,f})) T_f \jump{\bar p^\a}_\aa^-  &   \mbox{ on } (0,T)\times \Gamma, \aa=\pm, \\[1ex]
& {(\bbsig(\bar\bu) - b ~ \bar p^E_m \mathbb I)\n^\pm = - \bar p^E_f \n^\pm}  & \mbox{ on } (0,T)\times \Gamma,\\[1ex]
& \bar d_f = -\jump{\bar\bu}   & \mbox{ on } (0,T)\times \Gamma,  
\end{array}
\right.
\end{equation}
the initial conditions
$$
\bar p^\a_\nu|_{t=0} = \bar p^{\a}_{0,\nu},  \quad \bar \phi_m|_{t=0}= \bar\phi_m^0, 
$$
and normal flux conservation for $\q^\a_f$ at fracture intersections not located at the boundary $\partial\Omega$. 
Above, the equivalent pressure $\bar p^E_\nu$, $\nu\in \{m,f\}$ is defined, following~\cite{coussy}, by
$$
\bar p^E_{\nu} = \dsp \sum_{\alpha \in \{\g,\l\} } \bar p^\alpha_\nu~S^\alpha_\nu(\bar p_{c,\nu})  - U_{\nu}(\bar p_{c,\nu}), 
$$
where
\begin{equation}\label{eq.capillary_energy}
U_{\rm rt}(\bar p_c) =  \int_0^{\bar p_c} q  \left( S^{\g}_{{\rm rt}} \right)'(q) \d q
\end{equation}
is the capillary energy density function for each rock type ${\rm rt}\in \{m,f,\pm\}$. As already noticed in \cite{KTJ2013,JJ14,bonaldi:hal-02549111}, this is a key choice to obtain the energy estimates that are the starting point for the convergence analysis.

We make the following main assumptions on the data:
\begin{enumerate}[label=(H\arabic*),leftmargin=*]
\item For each phase $\a \in\{\g,\l\}$ and rock type ${\rm rt}\in\{m,f,\pm\}$, the mobility function $\eta^\a_{\rm rt}$ is continuous, non-decreasing, and there exist
  $0< \eta_{\rm rt,{\rm min}}^\a \leq \eta^\a_{\rm rt,{\rm max}}< +\infty$  such that 
  $\eta^\a_{\rm rt,{\rm min}} \leq \eta^\a_{\rm rt}(s) \leq \eta^\a_{\rm rt,{\rm max}}$ for all $s\in [0,1]$.
  \label{first.hyp}
\item For each rock type ${\rm rt}\in\{m,f,\pm\}$, the non-wetting phase saturation function $S^\g_{\rm rt}$ is a non-decreasing {Lipschitz}
  continuous function with values in $[0,1]$, and   $S^\l_{\rm rt} = 1 - S^\g_{\rm rt}$.
  \label{second.hyp}
  \item For $\aa=\pm$, the width $\bar d_\aa$ and porosity $\bar\phi_\aa$  of the damaged rock are strictly positive constants. 
\item $b\in [0,1]$ is the Biot coefficient, $M> 0$ is the Biot modulus, and $\lambda >0$, $\mu>0$ are the Lam\'e coefficients. These coefficients are assumed to be constant for simplicity.
\item The initial matrix porosity satisfies $\bar\phi_m^0 \in L^\infty(\O)$. 
%\item The initial fracture aperture $\bar d_f^0$ satisfies  $\bar d_f^0(\x) \geq d_0(\x)$ for a.e.\ $\x \in \Gamma$.    
\item The initial pressures are such that $(\bar p^\alpha_{0,m},\bar p^\alpha_{0,f}) \in V^0$, $p^\alpha_{0,m}\in L^\infty(\Omega)$ and $\bar p^\alpha_{0,f}\in L^\infty(\Gamma)$, $\alpha\in \{\g,\l\}$.
\item The source terms satisfy $\mathbf f \in L^2(\Omega)^d$, $h_m^\alpha\in L^2((0,T)\times \Omega),$ and $h_f^\alpha \in L^2((0,T)\times \Gamma)$.
  \item The normal fracture transmissivity $T_f\in L^\infty(\Gamma)$ is uniformly bounded from below by a strictly positive constant. 
 \item The matrix permeability tensor $\K_{m}$ is symmetric and uniformly elliptic on $\Omega$.
 \label{last.hyp}
\end{enumerate}
Let us denote by $C_c^\infty([0,T)\times\Omega\setminus\overline\Gamma)$ the space of smooth functions $\bar v:[0,T]\times (\O\setminus\overline\G)\to\R$ vanishing on $\partial\Omega$ and at $t=T$, and whose derivatives of any order admit finite limits on each side of $\Gamma$.
We will also use the boldface notation $\mathbf{C}_c^\infty([0,T]\times\Omega\setminus\overline\Gamma)$ for $C_c^\infty([0,T)\times\Omega\setminus\overline\Gamma)^d$.
  
\begin{definition}[Weak solution of the model]
  A \emph{weak solution} of the model is given by $\bar p^\a = (\bar p^\alpha_m, \bar p^\a_f) \in L^2(0,T;V^0)$, $\alpha\in \{\g,\l\}$, and $\bar \bu \in L^\infty(0,T;\U^0)$, such that, for any $\alpha\in \{\g,\l\}$, $\bar{d}_f^{\;\nf 3 2}\nabla_\tau\bar p_f^\alpha\in L^2((0,T)\times\G))^{d-1}$ and, for all $\bar \varphi^\a = (\bar \varphi^\alpha_m,\bar\varphi^a_f) \in C_c^\infty([0,T)\times\Omega\setminus\overline\Gamma)\times C_c^\infty([0,T)\times\Gamma)$ and all $\bar \bv \in \mathbf{C}_c^\infty([0,T]\times\Omega\setminus\overline\Gamma)$,
\begin{subequations}\label{eq_pb_weak}
\begin{equation}
\label{eq_var_hydro}      
\left.
\begin{array}{ll}
&\dsp \int_0^T \int_\O \(-\bar \phi_m S^{\a}_m(\bar p_{c,m}) \partial_t \bar \varphi_m^\a + \eta^\a_m ( S^{\a}_m(\bar p_{c,m})) \K_m  \nabla \bar p^\a_m \cdot \nabla \bar \varphi^\a_m\) \d\x \d t \\[2ex]
  &   + \dsp \int_0^T \int_\G \( -\bar d_f S^{\a}_f(\bar p_{c,f}) \partial_t \bar \varphi^\a_f  + \eta^\a_f ( S^{\a}_f(\bar p_{c,f})) {\bar d_f^{\;3}\over 12}  \nabla_\tau \bar p^\a_f \cdot \nabla_\tau \bar \varphi^\a_f \) \d\sigma(\x) \d t\\[2ex]
  &   + \dsp \sum_{\aa=\pm} 
  \int_0^T \int_\G  \( \bar Q^\a_{f,\aa}  \jump{\bar \varphi^\a}_\aa   
    - \bar d_\aa \bar \phi_\aa S^{\a}_\aa( \gamma_\aa \bar p_{c,m}) \partial_t \gamma_\aa \bar \varphi^\a_m\) \d\sigma(\x) \d t\\[3ex]
  & \dsp  - \sum_{\aa=\pm}  \int_\G   \bar d_\aa \bar \phi_\aa S^{\a}_\aa( \gamma_\aa {\bar p}^0_{c,m})  \gamma_\aa \bar\varphi^\a_m(0,\cdot) \d\sigma(\x)  \\[3ex] 
  &  - \dsp \int_\O \bar \phi_{m}^0 S^{\a}_m(\bar p_{c,m}^0) \bar \varphi_m^\a(0,\cdot)\d\x
  - \int_\G \bar d_f^0 S^{\a}_f(\bar p_{c,f}^0) \bar\varphi^\a_f(0,\cdot)\d\sigma(\x) \\[2ex]
  & \qquad = \dsp \int_0^T \int_\O h_m^\a \bar \varphi^\a_m \d\x \d t + \int_0^T \int_\G h_f^\a \bar \varphi^\a_f  ~\d\sigma(\x) \d t,
\end{array}
\right.
\end{equation}
\begin{equation}
  \label{eq_var_meca}
\begin{array}{ll}
  & \dsp \int_0^T \int_\O \( \bbsig(\bar \bu): \bbeps(\bar \bv) - b ~\bar p_m^E \div(\bar \bv)\) \d\x \d t + \int_0^T \int_\G \bar p_f^E ~\jump{\bar \bv}  ~\d\sigma(\x) \d t \\[2ex]
&  \qquad\dsp = \int_0^T \int_\Omega \mathbf{f}\cdot\bar\bv ~\d\x \d t, 
\end{array}
\end{equation}
with
\begin{equation}
  \label{eq_cont_closures}
\begin{array}{l}
  \bar Q^\a_{f,\aa} = T_f \[ \eta_\aa(S^\a_\aa(\gamma_\aa \bar p_{c,m}))  \jump{\bar p^\a}_\aa^+  -  \eta_f(S^\a_f(\bar p_{c,f}))  \jump{\bar p^\a}_\aa^- \],\\[1ex]
  \bar p_{c,\nu} = \bar p^\g_\nu - \bar p^\l_\nu,\\[1ex]
  \bar d_f = -\jump{\bar\bu},\\
  \bar\phi_m -\bar \phi_m^0= \dsp b~\div(\bar\bu-\bar\bu^0) + \frac{1}{M} (\bar p^E_m -\bar p_m^{E,0}),
  \end{array}
\end{equation}
\end{subequations}
where $\bar d_f^0 = -\jump{\bar \bu^0}$ and $\bar \bu^0$ is the solution of \eqref{eq_var_meca} without the time integral and using the initial equivalent pressures 
$\bar p_m^{E,0}$ and $\bar p_f^{E,0}$ obtained from the initial pressures {$\bar p^\alpha_{0,m}$ and $\bar p^\alpha_{0,f}$, $\alpha\in \{\g,\l\}$.}
\label{def:weak.sol}
\end{definition}
\begin{remark}[Regularity of the fracture aperture]
Notice that, by the Sobolev--trace embeddings \cite[Theorem 4.12]{AF03}, $\bar \bu \in L^\infty(0,T;\U^0)$ implies that $\bar d_f = -\jump{\bar\bu}\in L^\infty(0,T;L^4(\G))$. All the integrals above are thus well-defined.
\end{remark}

\section{The gradient discretization method}\label{sec:gradientscheme}

The gradient discretization (GD) for the Darcy discontinuous pressure model, introduced in \cite{DHM16}, is defined by a finite-dimensional vector space of discrete unknowns
$$
X^0_{\D_p} = X^0_{\D_p^m} \times X^0_{\D_p^f}
$$
and 
\begin{itemize}
\item  two discrete gradient linear operators on the matrix and fracture domains
$$
 \nabla_{\D_p}^m : X^0_{\D_p^m} \rightarrow L^\infty(\Omega)^d, \quad \quad 
   \nabla_{\D_p}^f : X^0_{\D_p^f} \rightarrow L^\infty(\Gamma)^{d-1},
$$ 
\item two function reconstruction linear operators on the matrix and fracture domains 
$$
\Pi_{\D_p}^m : X^0_{\D_p^m} \rightarrow L^\infty(\Omega),\quad\quad \Pi_{\D_p}^f : X^0_{\D_p^f} \rightarrow L^\infty(\Gamma),
$$
\item for $\aa=\pm$, jump reconstruction linear operators $\jump{\cdot}^\aa_{\D_p}$: $X^0_{\D_p} \rightarrow L^\infty(\Gamma)$,  and trace reconstruction linear operators $\mathbb{T}^\aa_{\D_p}$: $X^0_{\D_p^m} \rightarrow L^\infty(\Gamma)$.
\end{itemize}  
The operators $\Pi_{\D_p}^m$, $\Pi_{\D_p}^f$, $\mathbb{T}^\aa_{\D_p}$ are assumed \emph{piecewise constant} \cite[Definition 2.12]{gdm}. 
A consequence of the piecewise-constant property is the following: there is a basis $(\mathbf{e}_i)_{i\in I}$ of $X^0_{\D_p^m}$ such that, if $v=\sum_{i\in I} v_i\mathbf{e}_i$ and if, for a mapping $g:\R\to\R$ with $g(0)=0$, we define $g(v)=\sum_{i\in I}g(v_i)\mathbf{e}_i\in X^0_{\D_p^m}$ by applying $g$ component-wise, then $\Pi_{\D_p}^{m} g(v)=g(\Pi_{\D_p}^m v)$ and  $\mathbb{T}_{\D_p}^{\aa} g(v)=g(\mathbb{T}_{\D_p}^\aa v)$. Note that the basis $(\mathbf{e}_i)_{i\in I}$ is usually canonical and chosen in the design of $X^0_{\D_p^m}$. The same property holds for $X^0_{\D_p^f}$ and $\Pi_{\D_p}^f$. 
The vector space $X^0_{\D_p}$ is endowed with the following quantity, assumed to define a norm:
$$
\|v\|_{\D_p} \coloneqq \displaystyle  \|\nabla_{\D_p}^m v\|_{L^2(\Omega)^d} 
+ \|d_{0}^{\nf 3 2}\nabla_{\D_p}^f v\|_{ L^2(\Gamma)^{d-1} } + \sum_{\aa=\pm} \|\jump{v}^\aa_{\D_p}\|_{L^2(\Gamma)}.
$$

The gradient discretization for the mechanics is defined by a finite-dimensional vector space of discrete unknowns $X^0_{\D_\bu}$ and
\begin{itemize}
\item  a discrete symmetric gradient linear operator  $\bbeps_{\D_\bu} : X^0_{\D_{\bu}} \rightarrow   L^2(\Omega, \S_{d}(\R))$
where $\S_{d}(\R)$ is the vector space of real symmetric matrices of size $d$, 
\item a displacement function reconstruction linear operator  
$\Pi_{\D_\bu} : X^0_{\D_{\bu}} \rightarrow L^2(\Omega)^d$,
\item a normal jump function reconstruction linear operator 
  $\jump{\cdot}_{\D_\bu} : X^0_{\D_{\bu}} \rightarrow L^4(\Gamma)$.
\end{itemize}
Let us define the divergence operator $\div_{\D_\bu}(\cdot) = \mbox{\rm Trace}(\bbeps_{\D_\bu}(\cdot))$, the stress tensor operator 
$$
\bbsig_{\D_\bu}(\bv) =  2\mu \bbeps_{\D_\bu}(\bv) + \lambda \, \div_{\D_\bu}(\bv) \mathbb{I}, 
$$
and the fracture width $ d_{f,\D_\bu} = -\jump{\bu}_{\D_\bu}$. 
It is assumed that the following quantity defines a norm on $X^0_{\D_{\bu}}$:
\begin{equation}\label{norm.Du}
\|\bv\|_{\D_\bu} \coloneqq \|\bbeps_{\D_\bu}(\bv)\|_{L^2(\Omega,\mathcal S_d(\mathbb R))}.
\end{equation}

A spatial GD can be extended into a space-time GD by complementing it with
\begin{itemize}
\item 
a discretization $ 0 = t_0 < t_1 < \dots < t_N = T $ of the time interval $[0,T]$,
\item
interpolators $I^\nu_{\D_p} \colon V_\nu^0 \rightarrow X^0_{\D_p^\nu}$, $\nu\in\{m,f\}$, and $J^m_{\D_p} \colon L^2(\Omega)\rightarrow X^0_{\D_p^m}$ of initial conditions.
\end{itemize}
For $n\in\{0,\ldots,N\}$, we denote by $\dtn = t_{n+1}-t_n$ the time steps, and by $\Delta t = \max_{n=0,\ldots,N} \dtn$ the maximum time step. 

Spatial operators are extended into space-time operators as follows. Let $\Psi_\D$ be a spatial GDM operator defined in $X_\D^0$ with $\D=\D_\bu$, $\D_p^m$ or $\D_p^f$,  and let $w=(w_n)_{n=0}^{N}\in (X^0_{\D})^{N+1}$. Then, its space-time extension is defined by 
$$
\Psi_{\D}w(0,\cdot) = \Psi_{\D}w_0 \mbox{ and, } \forall n\in\{0,\dots,N-1\}\,,\;\forall t\in (t_n,t_{n+1}],\,\;
\Psi_{\D}w(t,\cdot) = \Psi_{\D}w_{n+1}.
$$
For convenience, the same notation is kept for the spatial and space-time operators.
Moreover, we define the discrete time derivative as follows: for $f:[0,T]\to L^1(\Omega)$ piecewise constant on the time discretization, with $f_n=f_{|(t_{n-1},t_n]}$ and $f_0=f(0)$, we set $\delta_t f (t) = \frac{f_{n+1} - f_n}{\dtn}$ for all  $t\in (t_n,t_{n+1}]$, $n\in\{0,\ldots,N-1\}$. 
  
Notice that the space of piecewise constant $X^0_{\D}$-valued functions $f$ on the time discretization together with the initial value $f_0=f(0)$ can be identified with $(X^0_{\D})^{N+1}$. The same definition of discrete derivative can thus be given for an element $w\in(X^0_{\D})^{N+1}$.
Namely, \mbox{$\delta_t w \in (X^0_{\D})^{N}$} is defined by setting, for any $n\in\{0,\ldots,N-1\}$ and $t\in (t_{n},t_{n+1}]$, $\delta_t w(t)=(\delta_t w)_{n+1} \coloneqq \frac{w_{n+1} - w_n}{\dtn}$. If $\Psi_{\D}$ is a space-time GDM operator, by linearity the following commutativity property holds: $\Psi_{\D} \delta_t w (t,\cdot) = \delta_t(\Psi_{\D} w(t,\cdot))$.

The gradient scheme for \eqref{eq_edp_hydromeca}  consists in writing the weak formulation \eqref{eq_var_hydro}--\eqref{eq_var_meca} with continuous spaces and operators substituted by their discrete counterparts, after a formal integration by part: find $p^\alpha = (p^\a_m, p^\a_f) \in (X^0_{\D_p})^{N+1}$, $\alpha\in \{\g,\l\}$, and $\bu \in (X^0_{\D_\bu})^{N+1}$, such that for all $\varphi^\alpha = (\varphi^\a_m, \varphi^\a_f) \in (X_{\D_p}^0)^{N+1}$, $\bv \in (X^0_{\D_\bu})^{N+1}$ and $\alpha\in \{\g,\l\}$,
\begin{subequations}\label{eq:GS}
\begin{align}
&\left.\begin{array}{llll}
  && \dsp \int_0^T \int_\Omega \( \delta_t \(\phi_\D \Pi_{\D_p}^m s^\alpha_m \)\Pi_{\D_p}^m \varphi_m^\alpha 
  +  \eta_m^\alpha(\Pi_{\D_p}^m s_m^\alpha) \K_m \nabla_{\D_p}^m p^\alpha_m \cdot   \nabla_{\D_p}^m \varphi^\alpha_m \) \d\x \d t\\[2ex]
  && + \dsp \int_0^T \int_\Gamma \delta_t \(d_{f,\D_\bu} \Pi_{\D_p}^f s^\alpha_f \)\Pi_{\D_p}^f \varphi^\alpha_f \d\sigma(\x)\d t\\[2ex]
  && + \dsp \int_0^T \int_\Gamma  \eta_f^\alpha(\Pi_{\D_p}^f s_f^\alpha) {d_{f,\D_\bu}^3 \over 12} \nabla_{\D_p}^f p^\alpha_f \cdot   \nabla_{\D_p}^f \varphi^\alpha_f  \d\sigma(\x) \d t\\[2ex]
 &&   + \dsp \sum_{\aa=\pm}  
  \int_0^T \int_\G  \( Q^\a_{f,\aa}  \jump{\varphi^\a}^\aa_{\D_p}   
  +   \bar d_\aa \bar \phi_\aa \delta_t \(\mathbb{T}^\aa_{\D_p} s^{\a}_\aa \) \mathbb{T}^\aa_{\D_p}\varphi^\a_m\) \d\sigma(\x) \d t \\[3ex] 
  && = \dsp \int_0^T \int_\Omega h_m^\alpha \Pi_{\D_p}^m \varphi^\alpha_m \d\x \d t + \int_0^T \int_\Gamma h_f^\alpha \Pi_{\D_p}^f \varphi^\alpha_f \d\sigma(\x)\d t,\\[3ex]
  \end{array}\right.
  \label{GD_hydro}
\\
& 
  \left.\begin{array}{llll}
    &&   \dsp \int_0^T \int_\Omega \( \bbsig_{\D_\bu}(\bu) : \bbeps_{\D_\bu}(\bv)  
    - b ~\Pi_{\D_p}^m p_m^E~  \div_{\D_\bu}(\bv)\)  \d\x \d t\\
    && \quad \quad\quad + \dsp \int_0^T \int_\Gamma \Pi_{\D_p}^f p_f^E~  \jump{\bv}_{\D_\bu} \d\sigma(\x)\d t 
    = \int_0^T \int_\Omega \mathbf{f} \cdot \Pi_{\D_\bu} \bv ~\d\x \d t, 
  \end{array}\right.
  \label{GD_meca}
\end{align}
with the closure equations, for $\nu\in \{m,f\}$ and $\aa= \pm$,
\begin{equation}
  \left\{\begin{array}{ll}
  & Q^\a_{f,\aa} = T_f \[ \eta^\a_\aa(\mathbb{T}^\aa_{\D_p} s^\a_\aa)  (\jump{p^\a}^\aa_{\D_p})^+  -
  \eta^\a_f(\Pi^f_{\D_p}s^\a_f)  (\jump{p^\a}_{\D_p}^\aa)^- \],\\[2ex]
  & p_{c,\nu} = p^\g_\nu - p^\l_\nu,  \quad s^\alpha_\nu = S^\alpha_\nu(p_{c,\nu}), \quad s^\a_\aa = S^\a_\aa(p_{c,m}),\\[2ex]
  & \dsp p_\nu^E = \sum_{\alpha\in\{\g,\l\}}  p^\alpha_\nu s^\alpha_\nu - U_\nu(p_{c,\nu}),\\[4ex]
  &  \phi_{\D} - \Pi_{\D_p}^m  \phi_m^0 = b ~\div_{\D_\bu} (\bu-\bu^0) + {1\over M} \Pi_{\D_p}^m (p_m^E-p_m^{E,0}),\\[2ex]
  & d_{f,\D_\bu} = -\jump{\bu}_{\D_\bu},\\[2ex]
  & \bbsig_{\D_\bu}(\bv) =  2\mu \bbeps_{\D_\bu}(\bv) + \lambda \, \div_{\D_\bu}(\bv) \mathbb{I}.  
  \end{array}\right.
  \label{GD_closures}
\end{equation}
\end{subequations}
The initial conditions are given by $p^\a_{0,\nu} = I^\nu_{\D_p} \bar p^\a_{0,\nu}$ ($\a\in\{\g,\l\}$, $\nu\in \{m,f\}$), $\phi_m^0 = J_{\D_p}^m \bar \phi^0$, and the initial displacement $\bu^0$ is the solution in $X^0_{\D_\bu}$ of
\eqref{GD_meca} without the time variable and with the equivalent pressures obtained from the initial pressures $(p^\alpha_0)_{\alpha\in\{\g,\l\}}$.

\subsection{Properties of gradient discretizations}\label{properties_gd}

Let $(\D_p^l)_{l\in\N}$ and $(\D_\bu^l)_{l\in\N}$ be sequences of GDs. We state here the assumptions on these sequences which ensure that the solutions to the corresponding schemes converge. Most of these assumptions are adaptation of classical GDM assumptions \cite{gdm}, except for the chain-rule and product rule used in Subsection \ref{subsec:compactness} to obtain compactness properties; we note that all these assumptions hold for standard discretizations used in porous media flows.

Following  \cite{DHM16}, the spatial GD of the Darcy flow
$$
\D_p=\(X_{\D_p}^0,\nabla_{\D_p}^\nu,\Pi_{\D_p}^\nu,\mathbb{T}^\aa_{\D_p}, \jump{\cdot}^\aa_{\D_p}; \nu\in\{m,f\},\aa=\pm\), 
$$
is assumed to satisfy the following coercivity, consistency, limit-conformity and compactness properties.

\noindent {\bf Coercivity of $\D_p$}. Let $C_{\D_p} > 0$ be defined by 
\begin{equation}\label{def_CDdarcy}
C_{\D_p} = \max_{0 \neq v=(v_m,v_f) \in X_{\D_p}^0} {\|\Pi^m_{\D_p} v_m\|_{L^2(\Omega)} + \|\Pi^f_{\D_p} v_f \|_{L^2(\Gamma)} + \sum_{\aa=\pm} \|\mathbb{T}^\aa_{\D_p} v_m\|_{L^2(\Gamma)} \over \|v \|_{\D_p}}. 
\end{equation}
Then, 
a sequence of spatial GDs $(\D_p^l)_{l\in{\mathbb N}}$
 is said to be {\emph{coercive}} if there exists $\overline C_p >0$ 
such that $C_{\D^l_p} \leq \overline C_p$ for all $l\in \N$.

\noindent {\bf Consistency of $\D_p$}.
%Let $r>8$ be given, and for all $w = (w_m,w_f) \in C_c^\infty(\Omega\setminus\overline\Gamma)\times C^\infty_c(\Gamma)$ and $v=(v_m,v_f) \in X_{\D_p}^0$ let us define 
Let $r>8$ be given, and for all $(w_m,v_m) \in C_c^\infty(\Omega\setminus\overline\Gamma) \times X_{\D_p^m}^0$ and all $(w_f,v_f) \in C^\infty_c(\Gamma) \times X_{\D_p^f}^0$, let us define 
\begin{equation}\label{def_sDdarcy}
%\left.\begin{array}{r@{\,\,}c@{\,\,}ll}
\begin{aligned}
S_{\D_p^m}(w_m,v_m) = \ &  
\|\nabla^m_{\D_p} v_m  - \nabla w_m\|_{L^2(\Omega)} + \|\Pi^m_{\D_p} v_m  -w_m\|_{L^2(\Omega)} \\
&
+ \sum_{\aa=\pm} \( \|\mathbb{T}^\aa_{\D_p} v_m  - \gamma_\aa w_m\|_{L^2(\Gamma)} 
+ \|\jump{(v_m,0)}^\aa_{\D_p}  - \jump{(w_m,0)}_\aa\|_{L^2(\Gamma)}\),\\[2ex]
S_{\D_p^f}(w_f,v_f) = \ & \|\nabla^f_{\D_p} v_f  - \nabla_\tau w_f\|_{L^r(\Gamma)} + \|\Pi^f_{\D_p} v_f  -w_f\|_{L^r(\Gamma)} %\\
%&
 + \!\! \sum_{\aa=\pm} \|\jump{(0,v_f)}^\aa_{\D_p}  - \jump{(0,w_f)}_\aa\|_{L^2(\Gamma)} ,
%S_{\D_p}(w,v) &=&  
%\|\nabla^m_{\D_p} v_m  - \nabla w_m\|_{L^2(\Omega)} 
%+ \|\nabla^f_{\D_p} v_f  -\nabla_\tau w_f\|_{L^r(\Gamma)} \\
%&+& \|\Pi^m_{\D_p} v_m  -w_m\|_{L^2(\Omega)} + \|\Pi^f_{\D_p} v_f - w_f\|_{L^r(\Gamma)}\\
%&+&  \dsp \sum_{\aa=\pm} \( \|\mathbb{T}^\aa_{\D_p} v_m  - \gamma_\aa w_m\|_{L^2(\Gamma)} + \|\jump{v}^\aa_{\D_p}  - \jump{w}_\aa\|_{L^2(\Gamma)}\)
\end{aligned}
\end{equation}
and 
${\cal S}_{\D_p^\nu}(w_\nu) = \min_{v_\nu \in X_{\D_p^\nu}^0} S_{\D_p^\nu}(w_\nu,v_\nu)$, 
$\nu\in\{m,f\}$.
Then, a sequence of spatial GDs $(\D_p^l)_{l\in{\mathbb N}}$ is said to be {\emph{consistent}} if for all $w_\nu\in V^0_\nu$ one has
$
\lim_{l \rightarrow +\infty} {\cal S}_{\D_p^{\nu,l}}(w_\nu) = 0,
$
$\nu\in\{m,f\}$.
Moreover, if $(\D_p^l)_{l\in\N}$ is a sequence of \emph{space-time} GDs, then it is said to be consistent if the underlying sequence of spatial GDs is consistent as above and if, for any $\varphi=(\varphi_m,\varphi_f)\in V^0$ and $\psi\in L^2(\Omega)$,
as $l\to+\infty$,
\begin{equation}\label{eq:cons.stGD}
\begin{aligned}
& \Delta\!t^l\to 0\,,\\
& \|\Pi^m_{\D_p^l} I^m_{\D_p^l} \varphi_m-\varphi_m\|_{L^2(\Omega)}+
\|\Pi^f_{\D_p^l} I^f_{\D_p^l} \varphi_f-\varphi_f\|_{L^2(\Gamma)}
\!+\! \sum_{\aa=\pm} \( \|\mathbb{T}^\aa_{\D_p} I^m_{\D_p^l} \varphi_m  - \gamma_\aa \varphi_m\|_{L^2(\Gamma)}\) \to 0,\\[1ex]
&
\|\Pi^m_{\D_p^l} J_{\D_p^l}^m \psi-\psi\|_{L^2(\Omega)}\to 0.
 \end{aligned}
\end{equation}

\begin{remark}[Consistency]
In \cite{DHM16}, the consistency is only considered for $r=2$. We have here to adopt a slightly stronger assumption to deal with the coupling and non-linearity involving the fracture aperture $d_f$. Note that, under standard mesh regularity assumptions, this stronger consistency property is still satisfied for all classical GDs.
\end{remark}

\noindent {\bf Limit-conformity of $\D_p$}. For all $\q = (\q_m, \q_f) \in C^\infty(\Omega\setminus\overline\Gamma)^d \times C^\infty(\Gamma)^{d-1}$, $\varphi_\aa \in C^\infty(\Gamma)$, and $v=(v_m,v_f)\in X_{\D_p}^0$, let us define 
\begin{equation}\label{def_wDdarcy}
  \left.\begin{array}{r@{\,\,}c@{\,\,}l}
  W_{\D_p}(\q, \varphi_\aa,v) &=&\dsp 
  \int_\Omega \( \q_m \cdot \nabla^m_{\D_p} v_m +\ \Pi^m_{\D_p}v_m ~ \div(\q_m) \)\d\x \\
  && \dsp + \int_\Gamma \( \q_f \cdot \nabla^f_{\D_p} v_f +\ \Pi^f_{\D_p}v_f ~ \div_\tau(\q_f) \)\d\sigma(\x)\\
  && \dsp - \sum_{\aa=\pm} \int_\Gamma \mathbb{T}^\aa_{\D_p} v_m \q_m \cdot \n_\aa  \d\sigma(\x)\\
  && + \dsp \sum_{\aa=\pm} \int_\Gamma \varphi_\aa \(\mathbb{T}^\aa_{\D_p} v_m - \Pi^f_{\D_p}v_f - \jump{v_m}^\aa_{\D_p}\) \d\sigma(\x),  
  \end{array}\right.
\end{equation}
and 
${\cal W}_{\D_p}(\q, \varphi_\aa) = \max_{0\neq v \in X_{\D_p}^0}
\frac{1}{\|v \|_{\D_p}}|W_{\D_p}(\q, \varphi_\aa,v)|$. 
Then, a sequence of spatial GDs $(\D_p^l)_{l\in{\mathbb N}}$ 
is said to be {\emph{limit-conforming}} if for all $\q \in C^\infty(\Omega\setminus\overline\Gamma)^d \times C_c^\infty(\Gamma)^{d-1}$ and $\varphi_\aa \in C^\infty(\Gamma)$ one has 
$\lim_{l \rightarrow +\infty} {\cal W}_{\D_p^l}(\q, \varphi_\aa) = 0$. Here $C_c^\infty(\Gamma)^{d-1}$ denotes the space of functions whose restriction to each $\Gamma_i$ {is} in $C^\infty(\Gamma_i)^{d-1}$ tangent to $\Gamma_i$, compactly supported away from the tips, and satisfying normal flux conservation at fracture intersections not located at the boundary $\partial\Omega$. 

\noindent {\bf (Local) compactness} of $\D_p$. A sequence of spatial GDs  
$(\D_p^l)_{l\in{\mathbb N}}$ is said to 
be {\emph{locally compact}} if for all sequences $(v^l)_{l\in {\mathbb N}}{\in (X^0_{\D_p^l})_{l\in\N}}$
such that $\sup_{l\in\N}\|v^l\|_{\D_p^l}<+\infty$ and all compact sets $K_m\subset \Omega$ and $K_f\subset\Gamma$, such that $K_f$ is disjoint from the intersections $(\overline{\Gamma}_i\cap\overline{\Gamma}_j)_{i\not=j}$, the sequences $(\Pi^m_{\D_p^l}v^l)_{l\in\N}$ and  $(\Pi^f_{\D_p^l}v^l)_{l\in \N}$ are relatively compact in $L^2(K_m)$ and $L^2(K_f)$, respectively.

\begin{remark}[Local compactness through estimates of space translates]\label{rem:compact.equivalent}
For $K_m,K_f$ as above, set
\begin{equation*}
   \begin{array}{r@{\,\,}c@{\,\,}l}
  T_{\D_p^l,K_m,K_f}(\xi,\eta)=\max_{v=(v_m,v_f) \in X^0_{\D_p^l}\backslash\{0\}} \frac{1}{\|v\|_{\D_p^l}} \(
  && \dsp \|\Pi^m_{\D_p^l}v_m(\cdot+\xi)-\Pi^m_{\D_p^l}v_m\|_{L^2(K_m)}\\
  && +\dsp \sum_{i\in I}\|\Pi^f_{\D_p^l}v_f(\cdot+\eta_i)-\Pi^f_{\D_p^l}v_f\|_{L^2(K_f\cap\Gamma_i)} \\
  && +\dsp \sum_{i\in I} \sum_{\aa=\pm}\|\mathbb{T}^\aa_{\D_p^l}v_m(\cdot+\eta_i)-\mathbb{T}^\aa_{\D_p^l}v_m\|_{L^2(K_f\cap\Gamma_i)} \)
  \end{array}
\end{equation*}
where $\xi\in\R^d$, $\eta=(\eta_i)_{i\in I}$ with $\eta_i$ tangent to $\Gamma_i$; for $\xi$ and $\eta$ small enough,
this expression is well defined since $K_m$ and $K_f$ are compact in $\Omega$ and $\Gamma$, respectively. Following \cite[Lemma 2.21]{gdm}, An equivalent formulation of the local compactness property is: for all $K_m,K_f$ as above,
\begin{equation*}
\lim_{\xi,\eta\to 0}\sup_{l\in\N}T_{\D_p^l,K_m,K_f}(\xi,\eta)=0.
\end{equation*}
\end{remark}

\begin{remark}[Usual compactness property for GDs]
The standard compactness property for GD is not \emph{local} but \emph{global}, that is, on the entire domain not any of its compact subsets (see, e.g., \cite[Definition 2.8]{gdm} and also below for $\D_\bu$). Two reasons pushed us to consider here the weaker notion of local compactness: firstly, for standard GDs, the global compactness does not seem obvious to establish (or even true) in the fractures, because of the weight $d_0$ in the norm $\|{\cdot}\|_{\D_p}$, which prevents us from estimating the translates of the reconstructed function by the gradient near the fracture tips; secondly, we will only use compactness on saturations, which are uniformly bounded by 1 and for which local and global compactness are therefore equivalent.

In the following, for brevity we refer to the local compactness of $(\D_p^l)_{l\in\N}$ simply as the \emph{compactness} of this sequence of GDs.
\end{remark}

\noindent\textbf{Bounds on reconstruction operators of $(\D_p^l)_{l\in\N}$}.
\begin{itemize}[leftmargin=*]
\item \emph{Chain rule estimate}. For any Lipschitz-continuous function $F:\R\to\R$, there is $C_F\ge 0$ such that, for all $l\in\N$, and any $v_m\in X_{\D_p^{m,l}}^0$,  $$\|\nabla_{\D_p^l}^m F(v_m)\|_{L^2(\Omega)}\le C_F \|\nabla_{\D_p^l}^m v_m\|_{L^2(\Omega)}.$$

\item \emph{Product rule estimate}. There exists $C_P \ge 0$ such that, for any $l\in\N$ and any $u^l_m,v^l_m\in X_{\D_p^{m,l}}^0$, it holds
$$
\|\nabla^m_{\D_p} (u^l_m v^l_m) \|_{L^2(\Omega)} \leq C_P\( |u^l_m|_{\infty} \|\nabla^m_{\D_p} v^l_m \|_{L^2(\Omega)} + |v^l_m|_{\infty} \|\nabla^m_{\D_p} u^l_m \|_{L^2(\Omega)} \),
$$
where $|w_m|_{\infty}\coloneqq\max_{i\in I_m}|w_i|$ whenever $w_m=\sum_{i\in I_m}w_i\mathbf{e}_i$ with $(\mathbf{e}_i)_{i\in I_m}$ the canonical basis of $X^0_{\D_p^{m,l}}$.

\item \emph{Bound on the jump operator}. For any $l\in\mathbb N$ and any $v^l=(v^l_m,v^l_f) \in X_{\D_p^{l}}^0$, there is $C\ge 0$ such that
$$\| \jump{v^l}_\aa \|_{L^\infty(\Gamma)} \leq C \( |v^l_m|_{\infty} + |v^l_f|_{\infty} \),$$
where $|v_f|_{\infty}\coloneqq\max_{i\in I_f}|v_i|$ whenever $v_f=\sum_{i\in I_f}v_i\mathbf{e}_i$ with $(\mathbf{e}_i)_{i\in I_f}$ the canonical basis of $X^0_{\D_p^{f,l}}$.
\end{itemize}

\noindent {\bf Coercivity of $(\D_\bu^l)_{l\in\N}$}.
Let $C_{\D_\bu} > 0$ be defined by 
\begin{equation}\label{def_CDmeca}
C_{\D_\bu} = \max_{\mathbf 0 \neq \bv \in X_{\D_\bu^l}^0} {\|\Pi_{\D_\bu^l} \bv\|_{L^2(\Omega)} + \|\jump{\bv}_{\D_\bu^l}\|_{L^4(\Gamma)} \over \|\bv \|_{\D_\bu^l}}.
\end{equation}
Then, 
the sequence of spatial GDs $(\D_\bu^l)_{l\in{\mathbb N}}$
 is said to be {\emph{coercive}} if there exists $\overline C_{\bu} >0$ 
such that $C_{\D^l_\bu} \leq \overline C_{\bu}$ for all $l\in \N$.

\noindent {\bf Consistency of $(\D_\bu^l)_{l\in\N}$}. For all $\textbf{w}\in \U^0$, it holds $\lim_{l \rightarrow +\infty} {\cal S}_{\D_\bu^l}(\textbf{w}) = 0$ where
\begin{equation}\label{eq:def.SDu}
{\cal S}_{\D_\bu^l}(\textbf{w})=  {\min_{\bv \in X_{\D_\bu^l}^0}}\Big[
\|\bbeps_{\D_\bu^l}(\bv)  - \bbeps(\textbf{w})\|_{L^2(\Omega,\S_{d}(\R))} 
+ \|\Pi_{\D_\bu^l} \bv  -\textbf{w}\|_{L^2(\Omega)} + \left\|\jump{\bv}_{\D_\bu^l}  - \jump{\textbf{w}}\right\|_{L^4(\Gamma)}\Big].
\end{equation}

\noindent {\bf Limit-conformity of $(\D_\bu^l)_{l\in\N}$}.
Let $C_\Gamma^\infty(\Omega\setminus\overline\Gamma,\S_{d}(\R))$ denote the vector space of smooth functions $\bbsig : \Omega\setminus\overline\Gamma \to \S_{d}(\R)$ whose derivatives of any order admit finite limits on each side of $\Gamma$, and such that $\bbsig^+(\x) \n^+ + \bbsig^-(\x) \n^- = {\mathbf 0}$ and $(\bbsig^+(\x) \n^+) {\times} \n^+ = {\mathbf 0}$ for 
a.e.\ $\x\in \Gamma$. For all $\bbsig \in C_\Gamma^\infty(\Omega\setminus\overline\Gamma,\S_{d}(\R))$, it holds
$\lim_{l \rightarrow +\infty} {\cal W}_{\D_\bu^l}(\bbsig) = 0$ where
\begin{equation*}
 {{\cal W}_{\D_\bu^l}(\bbsig)}{} ={\max_{\mathbf 0 \neq \bv \in X_{\D_\bu^l}^0}\frac{1}{\|\bv \|_{\D_\bu^l}}}\left[
  \int_\Omega \( \bbsig : \bbeps_{\D_\bu^l}(\bv) + \Pi_{\D_\bu^l}\bv \cdot \div(\bbsig) \)\d\x 
 - \int_\Gamma  (\bbsig \n^+)\cdot\n^+ \jump{\bv}_{\D_\bu^l} \d\sigma(\x)\right].
\end{equation*}

\noindent {\bf Compactness of $(\D_\bu^l)_{l\in\N}$}. For any sequence $(\bv^l)_{l\in {\mathbb N}}{\in (X^0_{\D_\bu^l})_{l\in\N}}$ such that {$\sup_{l\in\N}\|\bv^l\|_{\D_\bu^l}<+\infty$}, the sequences $(\Pi_{\D_\bu^l}\bv^l)_{l\in\N}$ and $(\jump{\bv^l}_{\D_\bu^l})_{l\in \N}$ 
are relatively compact in $L^2(\Omega)^d$ and in $L^s(\Gamma)$ for all $s<4$, respectively.

\begin{remark}[Compactness through estimates of space translates]\label{rem:compact.equivalent.Du}
Similarly to Remark \ref{rem:compact.equivalent} (see also \cite[Lemma 2.21]{gdm}), the compactness of $(\D_\bu^l)_{l\in\N}$ is equivalent to
\begin{equation*}
\lim_{\xi,\eta\to 0}\sup_{l\in\N}T_{\D_\bu^l,s}(\xi,\eta)=0\quad\forall s<4,
\end{equation*}
where
\begin{equation*}
T_{\D_\bu^l,s}(\xi,\eta)=\max_{\bv\in X^0_{\D_\bu^l}\backslash\{0\}}\frac{\|\Pi_{\D_\bu^l}\bv(\cdot+\xi)-\Pi_{\D_\bu^l}\bv\|_{L^2(\Omega)}+\sum_{i\in I}\left\|\jump{\bv^l}_{\D_\bu^l}(\cdot+\eta_i)-\jump{\bv^l}_{\D_\bu^l}\right\|_{L^s(\Gamma_i)}}{\|\bv\|_{\D_\bu^l}},
\end{equation*}
with $\xi\in\R^d$, $\eta=(\eta_i)_{i\in I}$ with $\eta_i$ tangent to $\Gamma_i$, and the functions extended by $0$ outside their respective domain $\Omega$ or $\Gamma$.
\end{remark}

\section{Convergence analysis}\label{sec:convergence}
The main theoretical result of this work is the following convergence theorem. 

\begin{theorem}\label{thm.conv}
Let $(\D_p^l)_{l\in\N}$, $(\D_\bu^l)_{l\in\N}$, $\{(t^l_n)_{n=0}^{N^l}\}_{l\in\N}$, be sequences of space time GDs assumed to satisfy the properties described in Section~\ref{properties_gd}.
Let $\phi_{m,{\rm min}} >0$ and assume that, for each $l\in \N$, the gradient scheme \eqref{GD_hydro}--\eqref{GD_meca} has a solution $p^\alpha_l=(p^\a_{m,l},p^\a_{f,l})\in (X_{\D_p^l}^0)^{N^l+1}$, $\alpha\in \{\g,\l\}$, $\bu^l \in (X_{\D_\bu^l}^0)^{N^l+1}$ such that
  \begin{itemize}
\item[(i)] $d_{f,\D_\bu^l}(t,\x) \geq d_0(\x)$ for a.e.\ $(t,\x) \in (0,T)\times \Gamma$, 
\item[(ii)] $\phi_{\D^l}(t,\x) \geq \phi_{m,{\rm min}}$ for a.e.\ $(t,\x) \in (0,T)\times \Omega$.
  \end{itemize}
  Then, there exist $\bar p^\alpha=(\bar p^\a_m, \bar p^\a_f)\in L^2(0,T;V^0)$, $\alpha\in \{\g,\l\}$, $\bar \bu\in L^\infty(0,T;\U^0)$  and $\bar Q^\a_{f,\aa} \in L^2(0,T;L^2(\Gamma))$  satisfying the weak formulations \eqref{eq_var_hydro}-\eqref{eq_var_meca} such that for $\alpha\in \{\g,\l\}$ and up to a subsequence 
\begin{equation*}
%\left\{
\begin{array}{lll}
  & \Pi^m_{\D_p^l} p^\alpha_{m,l} \weakto \bar p^\alpha_m & \mbox{weakly in } L^2(0,T;L^2(\Omega)), \\[1ex]
  & \Pi^f_{\D_p^l} p^\alpha_{f,l} \weakto \bar p^\alpha_f & \mbox{weakly in } L^2(0,T;L^2(\Gamma)), \\[1ex]
  & \Pi_{\D_\bu^l} \bu^l \weakto \bar \bu & \mbox{weakly-$\star$ in } L^\infty(0,T;L^2(\Omega)^d),\\[1ex]
  & \phi_{\D^l} \weakto \bar \phi_m & \mbox{weakly-$\star$ in } L^\infty(0,T;L^2(\Omega)),\\[1ex]
  & d_{f,\D_\bu^l} \rightarrow \bar d_f & \mbox{in } L^\infty(0,T;L^p(\Gamma)) \mbox{ for } 2 \leq p < 4,\\[1ex]
  & \Pi^m_{D_p^l} S^\alpha_m(p_{c,m}^l) \rightarrow S^\alpha_m(\bar p_{c,m}) & \mbox{in } L^2(0,T;L^2(\Omega)) ,\\[1ex]
  & \Pi^f_{D_p^l} S^\alpha_f(p_{c,f}^l) \rightarrow S^\alpha_f(\bar p_{c,f}) & \mbox{in } L^2(0,T;L^2(\Gamma)),\\[1ex]
  & Q^\a_{f,\aa} \weakto \bar Q^\a_{f,\aa}  & \mbox{weakly in } L^2(0,T;L^2(\Gamma)).
\end{array}
%\right.
\end{equation*}
where $\bar \phi_m = \bar \phi_m^0 + \dsp b~\div(\bar\bu-\bar\bu^0) + \frac{1}{M} (\bar p^E_m -\bar p_m^{E,0})$,
$\bar d_f = -\jump{\bar\bu}$, and $\bar p_c = \bar p^\g - \bar p^\l$. 
\end{theorem}
We first present in Subsections \ref{subsec:apriori} and \ref{subsec:compactness} a sequence of intermediate results that will be useful for the proof of Theorem \ref{thm.conv} detailed in Subsection~\ref{subsec:convergence}.

\begin{remark}[Limit interface fluxes]\label{rem:noQf}
The theorem states that the limit functions satisfy all but the first closure equations in \eqref{eq_cont_closures}. It does not, however, identify the limit interface fluxes $\bar Q_{f,\aa}^\alpha$, $\alpha\in\{\g,\l\}$. This identification requires the limit functions to satisfy an energy equality, which is known under some assumption on the limit fracture width $\bar d_f$. See the discussion in Subsection~\ref{identification} for more details.
\end{remark}

\subsection{Energy estimates}\label{subsec:apriori}
Using the phase pressures and velocity (time derivative of the displacement field) as test functions, the following \emph{a priori} estimates can be inferred.
\begin{lemma}[A priori estimates]   
  \label{lemma_apriori}
  Let $p^\alpha=(p^\a_m,p^\a_f),\bu$ be a solution to problem~\eqref{eq:GS} such that 
  \begin{itemize}
\item[(i)] $d_{f,\D_\bu}(t,\x) \geq d_0(\x)$ for a.e.\ $(t,\x) \in (0,T)\times \Gamma$, 
\item[(ii)] $\phi_{\D}(t,\x) \geq  \phi_{m,\rm min}$ for a.e.\ $(t,\x) \in (0,T)\times \Omega$,
where $\phi_{m,\rm min}>0$ is a constant.
  \end{itemize}
Under hypotheses \emph{\ref{first.hyp}--\ref{last.hyp}}, there exists a real number $C > 0$ depending on the data, the coercivity constants $C_{\D_p}$, $C_{\D_\bu}$, and $\phi_{m,\rm min}$, such that the following estimates hold: 
\begin{equation}\label{apriori.est}
\begin{aligned}
\|\grad_{\D_p}^m p_m^\alpha\|_{L^2((0,T)\times\Omega)} \le C, 
&\quad& \| d_{f,\D_\bu}^{\nf 3 2} \grad_{\D_p}^f p_f^\alpha\|_{L^2((0,T)\times\Gamma)} \le C, \\
\| \jump{p^\a}^\aa_{\D_p} \|_{L^2((0,T)\times\Gamma)} \le C,
&\quad&  \|U_\aa(\mathbb{T}^\aa_{\D_p}p_{c,m})\|_{L^\infty(0,T;L^1(\Gamma))} \le C,\\
\|U_m(\Pi^m_{\D_p}p_{c,m})\|_{L^\infty(0,T;L^1(\Omega))} \le C,
&\quad&
\|d_0 U_f(\Pi^f_{\D_p}p_{c,f})\|_{L^\infty(0,T;L^1(\Gamma))} \le C,\\
\|\Pi^m_{\D_p} p^E_m\|_{L^\infty(0,T;L^2(\Omega))} \le C, 
&\quad&
\| \bbeps_{\D_\bu}(\bu)\|_{L^\infty(0,T;L^2(\Omega,\S_d(\R)))} \le C,\\
\|d_{f,\D_\bu}\|_{L^\infty(0,T;L^4(\Gamma))}\le C.
\end{aligned}
\end{equation}
\end{lemma}

{
\begin{remark}[Existence of the discrete solution]
Since these a priori estimates are obtained assuming lower bounds on the fracture aperture and porosity, the existence of the discrete solution cannot be deduced from these estimates and will be assumed in the following convergence analysis. As noticed in the introduction, these assumptions on lower bounds of the fracture aperture and porosity are mandatory since the model itself does not account for possible contact of fracture walls nor a nonlinear behavior of pore volume contraction, and thus cannot yield such lower bounds. The analysis of a model with contact is a topic for future work.
\end{remark}
}

\begin{proof}
For a piecewise constant function $v$ 
on~$[0,T]$ with $v(t)= v_{n+1}$ for all $t \in (t_n,t_{n+1}]$, $n\in\{0,\ldots,N-1\}$, and the initial value $v(0) = v_0$, we define the piecewise constant function $\hat v$ such that $\hat v(t) = v_n$ for all $t \in (t_n,t_{n+1}]$. We notice the following expression for the discrete derivative of the product of two such functions: 
\begin{equation}\label{der.prod}
\delta_t(u v)(t) = \hat u(t)\delta_t v(t) + v(t)\delta_t u(t).
\end{equation}
In \eqref{GD_hydro}, upon choosing $\varphi^\alpha = p^\alpha$ we obtain
$T_1+T_2+T_3+T_4 +\sum_{\aa=\pm} (T_1^\aa + T^\aa_2)= T_5 + T_6$, with
\begin{equation}\label{first.eq}
 \hspace{-.1cm}
\begin{array}{llll}
 T_1 = \dsp \int_0^T \int_\Omega \delta_t \(\phi_\D \Pi_{\D_p}^m s^\alpha_m \)\Pi_{\D_p}^m p^\alpha_m \d\x \d t, \quad & T_2 = \dsp \int_0^T \int_\Omega \eta_m^\alpha(\Pi_{\D_p}^m s_m^\alpha) \K_m \nabla_{\D_p}^m p^\alpha_m \cdot   \nabla_{\D_p}^m p^\alpha_m \d\x \d t,\\[2ex]
 T_3 = \dsp \int_0^T \int_\Gamma \delta_t \(d_{f,\D_\bu} \Pi_{\D_p}^f s^\alpha_f \)\Pi_{\D_p}^f p^\alpha_f \d\sigma(\x)\d t, \quad & T_4 = \dsp \int_0^T \int_\Gamma  \eta_f^\alpha(\Pi_{\D_p}^f s_f^\alpha) {d_{f,\D_\bu}^3 \over 12} \nabla_{\D_p}^f p^\alpha_f \cdot   \nabla_{\D_p}^f p^\alpha_f  \d\sigma(\x) \d t,\\[2ex]
 \dsp T^\aa_1 =  \int_0^T \int_\G   \bar d_\aa \bar \phi_\aa \delta_t \(\mathbb{T}^\aa_{\D_p} s^{\a}_\aa \) \mathbb{T}^\aa_{\D_p} p^\a_m \d\sigma(\x) \d t,  \quad& 
 \dsp T^\aa_2 =   \int_0^T \int_\G  \eta^\a_{\aa,f} T_f (\jump{p^\a}^\aa_{\D_p})^2 \d\sigma(\x) dt, \\[2ex]
   T_5 = \dsp \int_0^T \int_\Omega h_m^\alpha \Pi_{\D_p}^m p^\alpha \d\x \d t,
  \quad & T_6 = \dsp \int_0^T \int_\Gamma h_f^\alpha \Pi_{\D_p}^f p^\alpha \d\sigma(\x)\d t,
 \end{array}
\end{equation}
where  $\eta^\a_{\aa,f}=\eta^\a_\aa(\mathbb{T}^\aa_{\D_p} s^\a_\aa)$ if $\jump{p^\a}^\aa_{\D_p} \geq 0$ and $\eta^\a_{\aa,f}=\eta^\a_f(\Pi^f_{\D_p}s^\a_f)$ otherwise.  
First, we focus on the matrix and fracture accumulation terms $T_1$ and $T_3$, respectively. Using~\eqref{der.prod} and the piecewise constant function reconstruction property of $\Pi^\rt_{\D_p}$, $\rt\in\{m,f\}$, we can write 
$$
\begin{aligned}
\delta_t (\phi_\D S^\alpha_m(\Pi^m_{\D_p} p_{c,m})) & = 
\hat\phi_\D \delta_t S^\alpha_m(\Pi^m_{\D_p} p_{c,m}) + S^\alpha_m(\Pi^m_{\D_p} p_{c,m})
\delta_t \phi_\D,\\
\delta_t (d_{f,\D_\bu} S^\alpha_f(\Pi^f_{\D_p} p_{c,f})) & = 
\hat d_{f,\D_\bu} \delta_t S^\alpha_f(\Pi^f_{\D_p} p_{c,f}) + S^\alpha_f(\Pi^f_{\D_p} p_{c,f})
\delta_t d_{f,\D_\bu}.
\end{aligned}
$$
Summing on $\alpha\in\{\l,\g\}$, we obtain
$$\begin{aligned}
\sum_\alpha (T_1+T_3) 
= & \sum_\alpha \( \int_0^T \int_\Omega \hat\phi_\D \Pi^m_{\D_p} p^\alpha_m\,
 \delta_t S^\alpha_m(\Pi^m_{\D_p} p_{c,m}) \d\x \d t +
 \int_0^T \int_\Omega S^\alpha_m(\Pi^m_{\D_p} p_{c,m}) \Pi^m_{\D_p} p^\alpha_m \,\delta_t \phi_\D \d\x \d t \\
 &\hspace{-.3cm}+ \int_0^T \int_\Gamma \hat d_{f,\D_\bu} \Pi^f_{\D_p} p^\alpha_f\,
 \delta_t S^\alpha_f(\Pi^f_{\D_p} p_{c,f}) \d\sigma(\x) \d t+
 \int_0^T \int_\Gamma S^\alpha_f(\Pi^f_{\D_p} p_{c,f}) \Pi^f_{\D_p} p^\alpha_f\,
 \delta_t d_{f,\D_\bu} \d\sigma(\x) \d t \).
 \end{aligned}
$$
Now, for $\rt\in\{m,f\}$,
\begin{equation}\label{eq.der_cap_en}
\sum_\alpha \Pi^\rt_{\D_p} p^\alpha_\rt\,
 \delta_t S^\alpha_\rt(\Pi^\rt_{\D_p} p_{c,\rt}) = \Pi^\rt_{\D_p} p_{c,\rt}\,
 \delta_t S^\g_\rt(\Pi^\rt_{\D_p} p_{c,\rt}) \ge \delta_t U_\rt(\Pi^\rt_{\D_p} p_{c,\rt}).
 \end{equation}
Indeed, for $n\in\{0,\dots,N-1\}$, by the definition~\eqref{eq.capillary_energy}
of the capillary energy $U_\rt$ and letting \mbox{$\pi^n_{c,\rt}=\Pi^\rt_{\D_p}p^n_{c,\rt}$}, we have
$$
\begin{aligned}
\pi^{n+1}_{c,\rt}(S_{\rt}^\g(\pi^{n+1}_{c,\rt})-S_{\rt}^\g(\pi^n_{c,\rt})) & = U_\rt(\pi^{n+1}_{c,\rt}) - U_\rt(\pi^n_{c,\rt}) + \int_{\pi^n_{c,\rt}}^{\pi^{n+1}_{c,\rt}} (S_\rt^\g(q) - S_\rt^\g(\pi^n_{c,\rt}))\d q \\
& \ge U_\rt(\pi^{n+1}_{c,\rt}) - U_\rt(\pi^n_{c,\rt}),
\end{aligned}$$
where the last inequality holds since $S_\rt^\g$ is a non-decreasing function. Thus, we obtain
$$\begin{aligned}
\sum_\alpha {}&(T_1+T_3) 
\ge  \int_0^T \int_\Omega \hat\phi_\D \delta_t U_m(\Pi^m_{\D_p} p_{c,m}) \d\x \d t
+\int_0^T \int_\Gamma \hat d_{f,\D_\bu} \delta_t U_f(\Pi^f_{\D_p} p_{c,f}) \d\sigma(\x)\d t \\
& +\sum_\alpha\( \int_0^T \int_\Omega S^\alpha_m(\Pi^m_{\D_p} p_{c,m}) \Pi^m_{\D_p} p^\alpha_m \delta_t \phi_\D \d\x \d t + \int_0^T \int_\Gamma S^\alpha_f(\Pi^f_{\D_p} p_{c,f}) \Pi^f_{\D_p} p^\alpha_f
 \delta_t d_{f,\D_\bu} \d\sigma(\x)\d t \).
 \end{aligned}
$$
Applying again~\eqref{der.prod}, we have 
$$
\begin{aligned}
\hat\phi_\D \delta_t U_m(\Pi^m_{\D_p} p_{c,m}) & = 
\delta_t(\phi_\D U_m(\Pi^m_{\D_p} p_{c,m})) - U_m(\Pi^m_{\D_p} p_{c,m}) \delta_t\phi_\D,\\
\hat d_{f,\D_\bu} \delta_t U_f(\Pi^f_{\D_p} p_{c,f}) & = \delta_t(d_{f,\D_\bu} U_f(\Pi^f_{\D_p} p_{c,f}))
- U_f(\Pi^f_{\D_p} p_{c,f}) \delta_t d_{f,\D_\bu}.
\end{aligned}
$$
In the light of the closure equations~\eqref{GD_closures}, this allows us to infer that
\begin{equation}\label{lhs_accumulation}
\begin{aligned}
\sum_\alpha (T_1+T_3) \ge 
& \int_0^T\int_\Omega \delta_t(\phi_\D U_m(\Pi^m_{\D_p} p_{c,m})) 
\d\x \d t + \int_0^T \int_\Gamma \delta_t (d_{f,\D_\bu} U_f(\Pi^f_{\D_p} p_{c,f})) \d\sigma(\mathbf x) \d t \\
& + \int_0^T \int_\Omega \frac{1}{2M} \delta_t\(\Pi^m_{\D_p} p^E_m\)^2 \d\x \d t + 
\int_0^T \int_\Omega b\, \Pi^m_{\D_p} p^E_m \, \div_{\D_\bu} (\delta_t \bu) \d\x \d t \\
& - \int_0^T \int_\Gamma \Pi^f_{\D_p} p^E_f\, \jump{\delta_t \bu}_{\D_\bu} \d\sigma(\x) \d t,
\end{aligned}
\end{equation}
where we have used the fact that, for $v$ piecewise constant on~$[0,T]$,
\begin{equation}\label{eq:vdtv}
v\delta_t v \ge \delta_t\left(\frac{v^2}{2}\right).
\end{equation}

Using, as in \eqref{eq.der_cap_en}, the relation
$$
\sum_\alpha \mathbb{T}^\aa_{\D_p} p^\alpha_m\,
 \delta_t S^\alpha_\aa(\mathbb{T}^\aa_{\D_p} p_{c,m}) = \mathbb{T}^\aa_{\D_p} p_{c,m}\,
 \delta_t S^\g_\aa(\mathbb{T}^\aa_{\D_p} p_{c,m}) \ge \delta_t U_\aa(\mathbb{T}^\aa_{\D_p} p_{c,m}),
$$
we obtain
\begin{equation}
 \label{acc_interface}
\sum_\alpha T^\aa_1 \geq  \int_0^T \int_\Gamma {\bar d_\aa \bar \phi_\aa} \delta_t (U_\aa(\mathbb{T}^\aa_{\D_p} p_{c,m})) \d\sigma(\mathbf x) \d t. 
\end{equation}
    
Then, taking into account assumptions~\ref{first.hyp}--\ref{last.hyp} and (i) in the lemma, there exists a real number $C > 0$ depending only on the data such that
\begin{equation}\label{lhs_diffusion}
\begin{aligned}
\sum_\alpha (T_2+T_4 + \sum_{\aa=\pm} T_2^\aa) \ge C & \( \int_0^T \int_\Omega \sum_\alpha |\grad_{\D_p}^m p^\alpha_m|^2 \d\x \d t 
+  \int_0^T \int_\Gamma \sum_\alpha |d_{f,\D_\bu}^{\nf 3 2} \grad^f_{\D_p} p^\alpha_f|^2 \d\sigma(\x) \d t \\
& 
+\int_0^T \int_\Gamma \sum_\alpha |\jump{p^\a}^\aa_{\D_p}|^2 \d\sigma(\x) \d t
\).
\end{aligned}
\end{equation}

On the other hand, upon choosing $\bv = \delta_t\bu$ in \eqref{GD_meca}, we get
$T_7+T_8+T_9=T_{10}$, with
\begin{equation}\label{second.eq}
 \hspace{-.1cm}
 \begin{array}{llll}
 T_7 = \dsp \int_0^T \int_\Omega \bbsig_{\D_\bu}(\bu) : \bbeps_{\D_\bu}(\delta_t \bu) \d\x \d t,  \quad 
 & \ \, T_8 = - \dsp \int_0^T \int_\Omega b\, \Pi_{\D_p}^m p_m^E\,  \div_{\D_\bu}(\delta_t \bu) \d\x \d t\,\\[2ex]
   T_9 = \dsp \int_0^T \int_\Gamma \Pi_{\D_p}^f p_f^E \, \jump{\delta_t \bu}_{\D_\bu} \d\sigma(\x) \d t, \quad 
   & T_{10} = \dsp \int_0^T \int_\Omega \mathbf{f} \cdot \Pi_{\D_\bu} (\delta_t \bu) \d\x \d t.
 \end{array}
\end{equation}
Using \eqref{eq:vdtv} and developing the definition of $\bbsig_{\D_\bu}$, we see that
\begin{equation}\label{est_elastic_energy}
T_7 \ge  \int_0^T \int_\Omega  \delta_t \( \frac{1}{2} \bbsig_{\D_\bu}(\bu) : \bbeps_{\D_\bu}(\bu) \) \d\x \d t,
\end{equation}
so that, all in all, taking into account
$$
\sum_\alpha \(T_1+T_2+T_3+T_4 + \sum_{\aa=\pm} (T_1^\aa+ T_2^\aa) \) + T_7 + T_8 + T_9 = \sum_\alpha(T_5 + T_6) + T_{10},
$$
and inequalities \eqref{lhs_accumulation}, \eqref{lhs_diffusion} and \eqref{est_elastic_energy},
we obtain the following estimate
for the solutions of \eqref{eq:GS}: there is a real number $C > 0$ depending on the data such that
\begin{equation}\label{GD_energy_estimate}
\begin{aligned}
&\int_0^T \int_\Omega \delta_t(\phi_\D U_m(\Pi^m_{\D_p} p_{c,m})) \, \d\mathbf x \d t
  + \int_0^T \int_\Gamma \delta_t (d_{f,\D_\bu} U_f(\Pi^f_{\D_p} p_{c,f})) \, \d\sigma(\mathbf x) \d t \\
  & + \sum_{\aa=\pm} \int_0^T \int_\Gamma \delta_t U_\aa(\mathbb{T}^\aa_{\D_p} p_{c,m}) \, \d\sigma(\mathbf x) \d t\\
& + \int_0^T \int_\Omega \delta_t\left(\frac{1}{2}\bbsigma_{\D_\bu}(\bu):\bbeps_{\D_\bu}(\bu) 
+ \frac{1}{2M}(\Pi^m_{\D_p}p^E_m)^2\right)\, \d\x \d t  \\
&+\sum_{\alpha} \int_0^T\int_\Omega |\grad^m_{\D_p}p^\alpha_m|^2 \, \d\mathbf x \d t + 
\sum_{\alpha} \int_0^T\int_\Gamma | d_{f,\D_\bu}^{\nf 3 2} \grad^f_{\D_p}p^\alpha_f |^2 \, \d\sigma(\mathbf x) \d t \\
& +\sum_\alpha \sum_{\aa=\pm}\int_0^T \int_\Gamma  |\jump{p^\a}^\aa_{\D_p}|^2 \d\sigma(\x) \d t\\
& \le C\left( \int_0^T \int_\Omega \mathbf f \cdot \delta_t \Pi_{\D_\bu} \bu \, \d\mathbf x \d t +
\sum_{\alpha} \int_0^T \int_\Omega h_m^\alpha \Pi^m_{\D_p}p^\alpha_m\,\d\mathbf x \d t %\right. %\\ 
%& \qquad 
+  \sum_{\alpha} \int_0^T \int_\Gamma h_f^\alpha \Pi^f_{\D_p}p^\alpha_f\,\d\sigma(\mathbf x) \d t \right).
\end{aligned}
\end{equation}
Now, we have
$$
\begin{alignedat}{1}
 \int_0^T \int_\Omega \mathbf f \cdot \delta_t \Pi_{\D_\bu} \bu \, \d\x \d t
   & = \int_\Omega \mathbf f \cdot (\Pi_{\D_\bu}\bu(T) - \mathbf f\cdot \Pi_{\D_\bu} \bu(0)) \d\x \\
    & \le  C_{\D_\bu} \|\mathbf f\|_{L^2(\Omega)}( \|\bbeps_{\D_\bu}(\bu)(T)\|_{L^2(\Omega,\S_d(\R))} + \|\bbeps_{\D_\bu}(\bu)(0)\|_{L^2(\Omega,\S_d(\R))}),
\end{alignedat}
$$
$$
\begin{alignedat}{1}
& \sum_{\alpha} \( \int_0^T \int_\Omega h_m^\alpha \Pi^m_{\D_p}p^\alpha_m\,d\mathbf x \d t  +  \int_0^T \int_\Gamma h_f^\alpha \Pi^f_{\D_p}p_f^\alpha\,d\sigma(\mathbf x) \d t \)
\\
& \le C_{\D_p} \sum_\alpha \(\|h_m^\alpha\|_{L^2((0,T)\times\Omega)} + 
\|h_f^\alpha\|_{L^2((0,T)\times\Gamma)}\)\\
& \quad\quad \times \(\|\grad^m_{\D_p} p^\alpha_m\|_{L^2(0,T;L^2(\Omega))}
+\| d_{f,\D_\bu}^{\nf 3 2} \grad^f_{\D_p} p^\alpha_f\|_{L^2(0,T;L^2(\Gamma))} +  \sum_{\aa=\pm} \| \jump{p^\a}^\aa_{\D_p}\|_{L^2(0,T;L^2(\Gamma))} \),
\end{alignedat}
$$
where we have used the coercivity properties of the two gradient discretizations along
with the Cauchy--Schwarz inequality and $d_0\le d_{f,\D_\bu}$. Using Young's inequality in the last two estimates as well as hypotheses \ref{first.hyp}--\ref{last.hyp} and (ii) in the lemma, and using telescopic sums on the terms involving $\delta_t$, it is then possible to infer from~\eqref{GD_energy_estimate} the existence of a real number $C > 0$ depending on the data and on $\phi_{m,\rm{min}}$ such that
$$
\begin{multlined}
  \| U_m(\Pi^m_{\D_p} p_{c,m})(T)\|_{L^1(\Omega)} + \| d_0 U_f(\Pi^f_{\D_p} p_{c,f})(T)\|_{L^1(\Gamma)} \\
  + \sum_{\aa=\pm} \| U_\aa(\mathbb{T}^\aa_{\D_p} p_{c,m})(T)\|_{L^1(\Gamma)}
+ \|(\Pi^m_{\D_p} p^E_m)(T)\|_{L^2(\Omega)}^2 
+ \|\bbeps_{\D_\bu}(\bu)(T)\|^2_{L^2(\Omega,\S_d(\R))}\\
+ \sum_\alpha \( \|\grad^m_{\D_p}p^\alpha_m\|_{L^2(0,T;L^2(\Omega))}^2 +
\| d_{f,\D_\bu}^{\nf 3 2} \grad^f_{\D_p}p^\alpha_f\|_{L^2(0,T;L^2(\Gamma))}^2 +  \sum_{\aa=\pm} \| \jump{p^\a}^\aa_{\D_p}\|^2_{L^2(0,T;L^2(\Gamma))} \) \\
\le C \( \| \mathbf f \|_{L^2(\Omega)}^2 + 
\sum_\alpha \( \|h_m^\alpha\|_{L^2((0,T)\times\Omega)}^2 + 
\|h_f^\alpha\|_{L^2((0,T)\times\Gamma)}^2 \)  \\
+  \| U_m(\Pi^m_{\D_p} p_{c,m})(0)\|_{L^1(\Omega)} + \| d_{f,\D_\bu}(0) U_f(\Pi^f_{\D_p} p_{c,f})(0)\|_{L^1(\Gamma)}\\
+   \sum_{\aa=\pm} \| U_\aa(\mathbb{T}^\aa_{\D_p} p_{c,m})(0)\|_{L^1(\Gamma)}
+ \|(\Pi^m_{\D_p} p^E_m)(0)\|_{L^2(\Omega)}^2 
+ \|(\Pi^f_{\D_p} p^E_f)(0)\|_{L^2(\Gamma)}^2 \)
\end{multlined}
$$
The above inequality,along with the fact that $T$ can be replaced by any $t\in (0,T]$ in the left-hand side, and in view of \eqref{eq:cons.stGD}--\eqref{conv_pE0}--\eqref{conv_u0}, yields the a priori estimates~\eqref{apriori.est} on $p^\alpha_\nu$, $p_{c,\nu}$, $p_m^E$ and $\bu$. The estimate on $d_{f,\D_\bu}$ follows from its definition and from the definition \eqref{def_CDmeca} of $C_{\D_\bu}$.
\end{proof}

\subsection{Compactness properties}\label{subsec:compactness}

Throughout the analysis, we write $a\lesssim b$ for $a\le Cb$ with constant $C$ depending only on the coercivity constants $C_{\D_p}$, $C_{\D_\bu}$ of the considered GDs, and on the physical parameters.

\subsubsection{Estimates on time translates}

\begin{proposition}\label{prop_timetranslates}
  Let $\D_p$, $\D_\bu$, $(t_n)_{n=0}^N$ be given space time GDs and $\phi_{m,{\rm min}} > 0$. It is assumed that the gradient scheme \eqref{GD_hydro}--\eqref{GD_meca} has a solution $p^\alpha=(p^\a_m,p^\a_f)\in (X_{\D_p}^0)^{N+1}$, $\alpha\in \{\g,\l\}$, $\bu \in (X_{\D_\bu}^0)^{N+1}$ such that
  $\phi_{\D}(t,\x) \geq \phi_{m,{\rm min}}$ for a.e.\ $(t,\x) \in (0,T)\times \Omega$ and $d_{f,\D_\bu}(t,\x) \geq d_0(\x)$ for a.e.\ $(t,\x) \in (0,T)\times \Gamma$. 
Let $\tau,\tau'  \in (0,T)$ and, for $s\in (0,T]$, denote by $n_s$ the natural number such that $s\in (t_{n_s},t_{n_s+1}]$. For any $\varphi =(\varphi_m,\varphi_f) \in X^0_{\D_p}$, it holds
  \begin{equation}
    \label{est_timetranslates}
\begin{array}{ll}
& \Big| 
\dsp \< [\phi_\D \Pi^m_{\D_p} s^\alpha_m](\tau) - [\phi_\D \Pi^m_{\D_p} s^\alpha_m](\tau'), \Pi^m_{\D_p} \varphi_m \>_{L^2(\Omega)} \\\\
& \qquad +\,\,   
\dsp \< [d_{f,\D_{\bu}} \Pi^f_{\D_p} s^\alpha_f](\tau) - [d_{f,\D_{\bu}} \Pi^f_{\D_p}s^\alpha_f](\tau'), \Pi^f_{\D_p} \varphi_f \>_{L^2(\Gamma)}  \\\\
& \qquad + \,\, \dsp \sum_{\aa=\pm} \dsp \< \bar d_{\aa}\bar\phi_\aa  \[\mathbb{T}^\aa_{\D_p} s^\alpha_\aa(\tau) -   \mathbb{T}^\aa_{\D_p} s^\alpha_\aa(\tau')\], \mathbb{T}^\aa_{\D_p} \varphi_m \>_{L^2(\Gamma)} 
\Big|\\\\
& \lsim  \dsp 
\sum_{n = n_{\tau}+1}^{n_{\tau'}}
 \dtn  \left( 
 \xi^{(1),\alpha,n+1}_m  \| \nabla^m_{\D_p} \varphi_m \|_{L^2(\Omega)} +
 \xi^{(1),\alpha,n+1}_f \| \nabla^f_{\D_p} \varphi_f \|_{L^8(\Gamma)}  \right.\\
& \qquad\qquad\qquad\qquad + \,\, \xi^{(2),\alpha,n+1}_m  \| \Pi^m_{\D_p} \varphi_m \|_{L^2(\Omega)} +
 \xi^{(2),\alpha,n+1}_f \| \Pi^f_{\D_p} \varphi_f \|_{L^2(\Gamma)}  %\\
%& \qquad\qquad\qquad\qquad
\left. 
+ \,\dsp \sum_{\aa=\pm}  \xi^{(1),\alpha,n+1}_\aa \| \jump{\varphi}^\aa_{\D_p}\|_{L^2(\Gamma)} 
 \right),
\end{array}
\end{equation}
with 
$$
\sum^{N-1}_{n = 0} \dtn \[ \sum_{\rt=m,f,\pm} \left( \xi^{(1),\alpha,n+1}_{\rm rt} \right)^2  + \sum_{\rt=m,f} \left( \xi^{(2),\alpha,n+1}_{\rm rt} \right)^2  \]  \lsim 1,  
$$
and  
\begin{align*}
  & \xi^{(1),\alpha,n+1}_m  = \| \nabla^m_{\D_p} p^{\alpha,n+1}_{m} \|_{L^2(\Omega)},
  \quad\quad   \xi^{(1),\alpha,n+1}_f = \| (d^{n+1}_{f,\D_{\bu}})^{\nf 3 2} \nabla^f_{\D_p} p^{\alpha,n+1}_f \|_{L^2(\Gamma)}  \| d^{n+1}_{f,\D_{\bu}} \|^{\nf 3 2}_{L^4(\Gamma)},\\
  & \xi^{(2),\alpha,n+1}_m  = \Big\| {1\over \dtn}\int_{t_n}^{t_{n+1}} h^{\alpha}_m(t,\cdot)\d t \Big\|_{L^2(\Omega)}, \quad\quad 
  \xi^{(2),\alpha,n+1}_f  = \Big\| {1\over \dtn}\int_{t_n}^{t_{n+1}} h^{\alpha}_f(t,\cdot)\d t \Big\|_{L^2(\Gamma)},\\
  & \xi^{(1),\alpha,n+1}_\aa = \| \jump{p^{\a,n+1}}^\aa_{\D_p}\|_{L^2(\Gamma)}.  
\end{align*}
\end{proposition}
\begin{proof} 
For any $\varphi \in X^0_{\D_p}$, writing the difference of piecewise-constant functions at times $\tau$ and $\tau'$ as the sum of their jumps between these two times, one has 
\begin{align}
\Big| 
{}& \< [\phi_\D  \Pi^m_{\D_p} s^\alpha_m](\tau) - [\phi_\D \Pi^m_{\D_p} s^\alpha_m](\tau'), \Pi^m_{\D_p} \varphi_m \>_{L^2(\Omega)} \nonumber\\
& +
 \< [d_{f,\D_{\bu}} \Pi^f_{\D_p} s^\alpha_f](\tau) - [d_{f,\D_{\bu}} \Pi^f_{\D_p}s^\alpha_f](\tau'), \Pi^f_{\D_p} \varphi_f \>_{L^2(\Gamma)}  \nonumber\\
 &+ \dsp \sum_{\aa=\pm} \dsp \< \bar d_{\aa}\bar\phi_\aa \[\mathbb{T}^\aa_{\D_p} s^\alpha_\aa(\tau) -  \mathbb{T}^\aa_{\D_p} s^\alpha_\aa(\tau')\], \mathbb{T}^\aa_{\D_p} \varphi_m \>_{L^2(\Gamma)} \Big| \nonumber\\
 \leq{}& 
\sum_{n = n_{\tau}+1}^{n_{\tau'}}
  \dtn \Big| \< \delta_t [\phi_\D \Pi^m_{\D_p} s^\alpha_m](t_{n+1}), \Pi^m_{\D_p} \varphi_m \>_{L^2(\Omega)} \nonumber\\ 
  &  \dsp  \quad + \,\, \<  \delta_t  [d_{f,\D_{\bu}} \Pi^f_{\D_p} s^\alpha_f](t_{n+1}), \Pi^f_{\D_p} \varphi_f \>_{L^2(\Gamma)} 
    + \,\sum_{\aa=\pm}  \< \bar d_{\aa}\bar\phi_\aa\delta_t [ \mathbb{T}^\aa_{\D_p} s^\alpha_\aa](t_{n+1}), \mathbb{T}^\aa_{\D_p} \varphi_m \>_{L^2(\Gamma)} 
\Big|.
\label{est_tt1}
\end{align}
From the gradient scheme discrete variational equation \eqref{GD_hydro}, we deduce that 
\begin{align}
\Big| 
{}& \< \delta_t  [\phi_\D \Pi^m_{\D_p} s^\alpha_m](t_{n+1}), \Pi^m_{\D_p} \varphi_m \>_{L^2(\Omega)} + 
 \dsp \< \delta_t  [d_{f,\D_{\bu}} \Pi^f_{\D_p} s^\alpha_f](t_{n+1}), \Pi^f_{\D_p} \varphi_f \>_{L^2(\Gamma)} \nonumber\\
 \dsp & + \,\sum_{\aa=\pm}  \< \bar d_{\aa}\bar\phi_\aa \delta_t [ \mathbb{T}^\aa_{\D_p} s^\alpha_\aa](t_{n+1}), \mathbb{T}^\aa_{\D_p} \varphi_m \>_{L^2(\Gamma)} 
 \Big| \nonumber\\
 \lsim{}&
 \| \nabla^m_{\D_p} p^{\a,n+1}_m \|_{L^2(\Omega)} ~ \| \nabla^m_{\D_p} \varphi_m \|_{L^2(\Omega)} 
 +
 \| (d^{n+1}_{f,\D_{\bu}})^{\nf 3 2} \nabla^f_{\D_p} p^{\a,n+1}_f \|_{L^2(\Gamma)} ~ \| (d^{n+1}_{f,\D_{\bu}})^{\nf 3 2} \nabla^f_{\D_p} \varphi_f \|_{L^2(\Gamma)}\nonumber\\
 & + \dsp \sum_{\aa=\pm}\| \jump{p^{\a,n+1}}^\aa_{\D_p}\|_{L^2(\Gamma)} \| ~\jump{\varphi}^\aa_{\D_p}\|_{L^2(\Gamma)} \nonumber\\
 & + \Big\|{1\over \dtn} \int_{t_n}^{t_{n+1}} h^{\alpha}_m(t,\cdot)\d t \Big\|_{L^2(\Omega)} ~\| \Pi^m_{\D_p} \varphi_m \|_{L^2(\Omega)} + 
 \Big\| {1\over \dtn} \int_{t_n}^{t_{n+1}} h^{\alpha}_f(t,\cdot)\d t \Big\|_{L^2(\Gamma)} ~  \|\Pi^f_{\D_p} \varphi_f \|_{L^2(\Gamma)}\nonumber\\
 \lsim{}&
\xi^{(1),\alpha,n+1}_m  \| \nabla^m_{\D_p} \varphi_m \|_{L^2(\Omega)} +
 \xi^{(1),\alpha,n+1}_f \| \nabla^f_{\D_p} \varphi_f \|_{L^8(\Gamma)} \nonumber\\
& + \, \xi^{(2),\alpha,n+1}_m  \| \Pi^m_{\D_p} \varphi_m \|_{L^2(\Omega)} +
 \xi^{(2),\alpha,n+1}_f \| \Pi^f_{\D_p} \varphi_f \|_{L^2(\Gamma)}   + \,\dsp \sum_{\aa=\pm}  \xi^{(1),\alpha,n+1}_\aa \| \jump{\varphi}^\aa_{\D_p}\|_{L^2(\Gamma)}, 
\label{est_tt2}
 \end{align}
 where the term $\| (d^{n+1}_{f,\D_{\bu}})^{\nf 3 2} \nabla^f_{\D_p} \varphi \|_{L^2(\Gamma)}$ has been estimated using the generalized H\"older inequality with exponents $(8,8/3)$, which satisfy $\frac{1}{8}+\frac{3}{8}=\frac{1}{2}$.
The result follows from \eqref{est_tt1}, \eqref{est_tt2}, the a priori estimates of Lemma \ref{lemma_apriori}, and from the assumptions $h_m^\alpha \in L^2((0,T)\times\Omega)$,
$h_f^\alpha \in L^2((0,T)\times\Gamma)$. 
\end{proof}

\begin{remark}
  Summing the estimate \eqref{est_timetranslates} on $\alpha\in \{\g,\l\}$, and using the fact that the two-phase saturations add up to $1$ in each medium, we obtain the following time translate estimates on $\phi_\D$ and $d_{f,\D_\bu}$:
  \begin{equation}
    \label{est_timetranslates_phi_df}
\begin{array}{ll}
& \Big| 
\dsp \< \phi_\D(\tau) - \phi_\D(\tau'), \Pi^m_{\D_p} \varphi_m \>_{L^2(\Omega)} + 
\< d_{f,\D_{\bu}}(\tau) - d_{f,\D_{\bu}}(\tau'), \Pi^f_{\D_p} \varphi_f \>_{L^2(\Gamma)} \Big| \\\\
& \lsim  \dsp 
\sum_{\alpha\in \{\g,\l\}}\sum_{n = n_{\tau}+1}^{n_{\tau'}}
\dtn  \left(
 \xi^{(1),\alpha,n+1}_m  \| \nabla^m_{\D_p} \varphi_m \|_{L^2(\Omega)} +
 \xi^{(1),\alpha,n+1}_f \| \nabla^f_{\D_p} \varphi_f \|_{L^8(\Gamma)}  \right.\\
& \qquad\qquad\qquad\qquad + \,\, \xi^{(2),\alpha,n+1}_m  \| \Pi^m_{\D_p} \varphi_m \|_{L^2(\Omega)} +
 \xi^{(2),\alpha,n+1}_f \| \Pi^f_{\D_p} \varphi_f \|_{L^2(\Gamma)} % \\
%& \qquad\qquad\qquad\qquad
\left. + \,\dsp \sum_{\aa=\pm}  \xi^{(1),\alpha,n+1}_\aa \| \jump{\varphi}^\aa_{\D_p}\|_{L^2(\Gamma)} 
 \right). 
\end{array}
\end{equation}
\end{remark}

\subsubsection{Compactness properties of $\Pi^m_{\D_p} s^\alpha_m$ and $\mathbb{T}^\aa_{\D_p} s^\a_\aa$}
\label{subsec:compactnessmat}

\begin{proposition}
  \label{prop_compactness_Sm}
  Let $(\D_p^l)_{l\in\N}$, $(\D_\bu^l)_{l\in\N}$, $\{(t^l_n)_{n=0}^{N^l}\}_{l\in\N}$ be sequences of space time GDs assumed to satisfy the coercivity and compactness properties, and such that $\lim_{l\rightarrow +\infty} \Delta t^l =0$.
  Let $\phi_{m,{\rm min}} > 0$ and assume that, for each $l\in \N$, the gradient scheme \eqref{GD_hydro}--\eqref{GD_meca} has a solution $p^\alpha_l=(p^\a_{m,l},p^\a_{f,l})\in (X_{\D_p^l}^0)^{N^l+1}$, $\alpha\in \{\g,\l\}$, $\bu^l \in (X_{\D_\bu^l}^0)^{N^l+1}$ such that
  $\phi_{\D^l}(t,\x) \geq \phi_{m,{\rm min}}$ for a.e.\ $(t,\x) \in (0,T)\times \Omega$ and $d_{f,\D^l_\bu}(t,\x) \geq d_0(\x)$ for a.e.\ $(t,\x) \in (0,T)\times \Gamma$. 
  Then, the sequences $(\Pi^m_{\D_p^l}s_m^{\alpha,l})_{l\in \N}$ and $(\mathbb{T}^\aa_{\D_p^l} s^{\a,l}_\aa)_{l\in\N}$,  with $s_m^{\alpha,l}= S_m^\alpha(p_{c,m}^l)$ and $s_\aa^{\alpha,l}= S_\aa^\alpha(p_{c,m}^l)$, are relatively compact in $L^2((0,T)\times \Omega)$ and $L^2((0,T)\times \Gamma)$, respectively. 
\end{proposition}

\begin{proof}
  The superscript $l\in \N$ will be dropped in the proof and all hidden constants in the following estimates are independent of $l$.
  Setting
  $$
  F^\a(p) = S^\a_m(p) + \sum_{\aa=\pm} S^\a_\aa(p), 
  $$
  it results from hypothesis \ref{second.hyp} that
  $$
  \(S^\a_{\rt}(p) - S^\a_{\rt}(p)\)^2 \leq \(S^\a_{\rt}(p) - S^\a_{\rt}(p)\)\(F^\a(p) - F^\a(q)\), 
  $$
  for ${\rm rt}\in\{m,\pm\}$. Using that $\phi_{\D}(t,\x) \geq \phi_{m,{\rm min}}$ for a.e.\ $(t,\x) \in (0,T)\times \Omega$  and noting that $\Pi^m_{\D_p} s_m^{\alpha}=S_m^\alpha(\Pi^m_{\D_p} p_{c,m})\in [0,1]$ and 
$\mathbb{T}^\aa_{\D_p} s_\aa^{\alpha}=S_\aa^\alpha(\mathbb{T}^\aa_{\D_p} p_{c,m})\in [0,1]$,
  we obtain
  $$
\begin{array}{l}
  \dsp \int_0^T \| \Pi^m_{\D_p} s^\alpha_m(\cdot + \tau,\cdot ) - \Pi^m_{\D_p} s^\alpha_m\|^2_{L^2(\Omega)} \d t
  + \sum_{\aa=\pm}  \int_0^T \| \mathbb{T}^\aa_{\D_p} s_\aa^{\alpha}(\cdot + \tau,\cdot )-  \mathbb{T}^\aa_{\D_p} s_\aa^{\alpha}\|_{L^2(\Gamma)}^2 \d t\\
  \lsim \dsp \tau + \int^{T-\tau}_0 \int_\O \phi_\D \( \Pi^m_{\D_p} s^\alpha_m(\cdot + \tau,\cdot ) - \Pi^m_{\D_p} s^\alpha_m \)
  \( F^\a(\Pi^m_{\D_p} p_{c,m}(\cdot + \tau,\cdot )) - F^\a(\Pi^m_{\D_p} p_{c,m}) \) \d\x \d t \\[3ex]
  + \dsp \sum_{\aa=\pm} \int^{T-\tau}_0 \int_\Gamma \bar d_\aa \bar \phi_\aa \( \mathbb{T}^\aa_{\D_p} s^\alpha_\aa(\cdot + \tau,\cdot ) - \mathbb{T}^\aa_{\D_p} s^\alpha_\aa \)
   \( F^\a(\mathbb{T}^\aa_{\D_p} p_{c,m}(\cdot + \tau,\cdot )) - F^\a(\mathbb{T}^\aa_{\D_p} p_{c,m}) \) \d\x \d t \\[3ex]
  =  \tau + T_1 + T_2, 
\end{array}
$$
where
\begin{align*}
  & \dsp T_1 = \int^{T-\tau}_0 \Big| \<  [\phi_\D \Pi^m_{\D_p} s^\alpha_m]( t + \tau ) - [\phi_\D \Pi^m_{\D_p} s^\alpha_m](t),  \Pi^m_{\D_p} \zeta_m^\alpha(t) \>_{L^2(\Omega)} \\
  & \quad\quad\quad\quad\quad\quad\quad\quad + \sum_{\aa=\pm} \<  \bar d_\aa\bar \phi_\aa \[ \mathbb{T}^\aa_{\D_p} s^\alpha_\aa( t + \tau ) - \mathbb{T}^\aa_{\D_p} s^\alpha_\aa(t)\],  \mathbb{T}^\aa_{\D_p} \zeta_m^\alpha(t) \>_{L^2(\Gamma)} \Big| \d t, \\
  & \dsp T_2 = \int^{T-\tau}_0 \Big| \<  \phi_\D( t + \tau ) - \phi_\D(t),  \Pi^m_{\D_p} \chi_m^\alpha(t) \>_{L^2(\Omega)} \Big| \d t,
\end{align*}
with 
$\zeta_m^\a(t) =  \( F^\a(p_{c,m}(t + \tau )) - F^\a(p_{c,m}(t)) \)$ and $\chi_m^\a(t) =  \zeta_m^\a(t)  ~s^\alpha_m(t+\tau)$. 
Let us set $\zeta^\a(t) = (\zeta^\a_m(t),0) \in X_{\D_p}^0$. In view of the estimates \eqref{est_timetranslates} for $\varphi =\zeta^\a(t)$,  we have 
$$
\begin{aligned}
  T_1 &\lsim \dsp \int^{T-\tau}_0 \sum_{n = n_t+1}^{n_{(t+\tau)}} \dtn
  \( \xi^{(1),\alpha,n+1}_m \| \nabla^m_{\D_p} \zeta^\alpha_m (t)\|_{L^2(\Omega)}
  +  \xi^{(2),\alpha,n+1}_m \| \Pi^m_{\D_p} \zeta^\alpha_m (t)\|_{L^2(\Omega)}\\
  & \quad\quad\quad\quad\quad\quad\quad\quad\quad\quad + \sum_{\aa=\pm} \xi^{(1),\alpha,n+1}_\aa \| \jump{\zeta^\alpha (t)}^\aa_{\D_p} \|_{L^2(\Gamma)}
  \) 
  ~\d t \\
&\lsim
 \dsp \int^{T-\tau}_0 \sum_{n = n_t+1}^{n_{(t+\tau)}} \dtn  \( (\xi^{(1),\alpha,n+1}_m)^2 + (\xi^{(2),\alpha,n+1}_m)^2 + \sum_{\aa=\pm} (\xi^{(1),\alpha,n+1}_\aa)^2\\
& \quad\quad\quad\quad\quad\quad\quad\quad\quad\quad +\, \| \nabla^m_{\D_p} \zeta^\alpha_m (t)\|^2_{L^2(\Omega)} + \| \Pi^m_{\D_p} \zeta^\alpha_m (t)\|_{L^2(\Omega)}^2 + \sum_{\aa=\pm} \| \jump{\zeta^\alpha(t)}^\aa_{\D_p} \|^2_{L^2(\Gamma)}
 \)~\d t.  
\end{aligned}
$$
From Proposition \ref{prop_timetranslates}, we have 
$$
\sum_{n = 0}^{N-1} \dtn  \((\xi^{(1),\alpha,n+1}_m)^2 + (\xi^{(2),\alpha,n+1}_m)^2 + \sum_{\aa=\pm} (\xi^{(1),\alpha,n+1}_\aa)^2
\) \lsim 1. 
$$
Using the a priori estimates of Lemma \ref{lemma_apriori}, $h_m^\alpha \in L^2((0,T)\times\Omega)$, the Lipschitz property and boundedness of $S^\alpha_m$ and $S^\a_\aa$, the chain rule  estimate on the sequence of GDs $(\D_p^l)_{l\in\N}$,  and the bound on the jump operator, we obtain that
$$
\int^{T-\tau}_0 \(\| \nabla^m_{\D_p} \zeta^\alpha_m (t)\|^2_{L^2(\Omega)} + \| \Pi^m_{\D_p} \zeta^\alpha_m (t)\|_{L^2(\Omega)}^2
+ \sum_{\aa=\pm} \| \jump{\zeta^\alpha (t)}^\aa_{\D_p} \|^2_{L^2(\Gamma)}
\)\d t \lsim 1. 
$$
%\corr{}{Il manque une hypothese sur $\jump{(.,0)}^\aa_{\D_p}$ de type Poincaré dans la coercivité? ou est ce une conséquence de la limite conformité qui serait du coup à rajouter dans les hypothèses de la prop?}[RM]. 
%\corr{}{}{[JD: je pense que c'est plus simple que ca. Toutes les composantes de $\zeta$ sont bornees puisque c'est un $F^\alpha$ (donc somme de saturations). Du coup, je pense qu'il suffit de rajouter une hypothese algebrique sur $\jump{\cdot}^\aa_{\D_p}$, du
%genre "borne par la norme $l^\infty$", comme on a fait avec la product rule. Du coup, aussi, je serai pour rassembler product rule, chain rule et borne du saut sous une seule propriete "bounds on reconstruction operators", avec 3 items dedans pour chacune de ces regles]}
We deduce from \cite[Lemma 4.1]{brenner.hilhorst2013} that $T_1 \lsim \tau + \Delta t$. 
Similarly, using the time translate estimate \eqref{est_timetranslates_phi_df} and the product rule  estimate on the sequence of GDs $(\D_p^l)_{l\in\N}$, one shows that $T_2 \lsim \tau+ \Delta t$, which provides the time translates estimates on $\Pi^m_{\D_p} s^\alpha_m$ in $L^2(0,T;L^2(\Omega))$ and on $\mathbb{T}^\aa_{\D_p} s^\alpha_\aa$ in $L^2(0,T;L^2(\Gamma))$. 

Let us consider any compact sets $K_m\subset \Omega$ and $K_f\subset\Gamma$. 
The space translates  estimates for $\Pi^m_{\D_p} s_m^\alpha$ in $L^2(0,T;L^2(K_m))$  and for $\mathbb{T}^\aa_{\D_p} s^\alpha_\aa$ in $L^2(0,T;L^2(K_f))$ derive from the a priori estimates of Lemma \ref{lemma_apriori}, the Lipschitz properties of 
$S^{\alpha}_m$ and $S^\a_\aa$, and from the local compactness property of the sequence of spatial GDs $(\D_p^l)_{l\in\N}$ (cf.~Remark \ref{rem:compact.equivalent}).
Combined with the time translate estimates above,  the Fr\'echet--Kolmogorov theorem implies that $\Pi^m_{\D_p} s_m^\alpha$  is relatively compact in $L^2(0,T;L^2(K_m))$ and that $\mathbb{T}^\aa_{\D_p} s^\alpha_\aa$ is relatively compact in $L^2(0,T;L^2(K_f))$. Since $\Pi^m_{\D_p} s_m^\alpha \in [0,1]$ and $\mathbb{T}^\aa_{\D_p} s^\alpha_\aa\in [0,1]$ it results that $\Pi^m_{\D_p} s_m^\alpha$  is relatively compact in $L^2(0,T;L^2(\Omega))$ and that $\mathbb{T}^\aa_{\D_p} s^\alpha_\aa$ is relatively compact in $L^2(0,T;L^2(\Gamma))$.

\end{proof}

\subsubsection{Uniform-in-time $L^2$-weak convergence of $d_{f,\D_\bu} \Pi^f_{\D_p} s^\alpha_f$ and $d_{f,\D_\bu}$}

\begin{proposition}
  \label{prop_uniftimeweakL2_dfsf}
  Let $(\D_p^l)_{l\in\N}$, $(\D_\bu^l)_{l\in\N}$, $\{(t^l_n)_{n=0}^{N^l}\}_{l\in\N}$
  be sequences of space time GDs assumed to satisfy the coercivity and consistency properties.
Let $\phi_{m,{\rm min}} >0$  and assume that, for each $l\in \N$, the gradient scheme \eqref{GD_hydro}--\eqref{GD_meca} has a solution $p^\alpha_l =(p^\a_{m,l},p^\a_{f,l}) \in (X_{\D_p^l}^0)^{N^l+1}$, $\alpha\in \{\g,\l\}$, $\bu^l \in (X_{\D_\bu^l}^0)^{N^l+1}$ such that
  \begin{itemize}
\item[(i)] $d_{f,\D_\bu^l}(t,\x) \geq d_0(\x)$ for a.e.\ $(t,\x) \in (0,T)\times \Gamma$, 
\item[(ii)] $\phi_{\D^l}(t,\x) \geq \phi_{m,{\rm min}}$  for a.e.\ $(t,\x) \in (0,T)\times \Omega$.
  \end{itemize}
  Then, the sequences $(d_{f,\D^l_\bu})_{l\in\N}$ and $(d_{f,\D^l_\bu} \Pi^f_{\D_p}s_f^{\alpha,l})_{l\in \N}$, with $s_f^{\alpha,l}= S_f^\alpha(p_{c,f}^l)$, converge up to a subsequence uniformly in time and weakly in $L^2(\Gamma)$, as per \cite[Definition C.14]{gdm}.
\end{proposition}

\begin{proof}
In the following, the superscript $l\in \N$ is dropped when not required for the clarity of the proof, and the hidden constants are independent of $l$.
Let $\overline{\varphi}_f \in C^\infty_c(\Gamma)$ and let $\varphi_f \in X^0_{\D_p^f}$ the element that realizes the minimum of ${\cal S}_{\D_p^f}(\bar\varphi_f,\varphi_f)$ in~\eqref{def_sDdarcy}. From Proposition \ref{prop_timetranslates} (with $\varphi_m=0$) we have 
\begin{align*}
&\Big| \< [d_{f,\D_\bu} \Pi^f_{\D_p} s^\alpha_f](\tau) - [d_{f,\D_\bu} \Pi^f_{\D_p} s^\alpha_f](\tau'), \Pi^f_{\D_p}  \varphi_f \>_{L^2(\Gamma)}\Big|\\
& \lsim   
 \max\left( \| \nabla^f_{\D_p}  \varphi_f \|_{L^8(\Gamma)}, \| \Pi^f_{\D_p}  \varphi_f \|_{L^2(\Gamma)} ,\max_{\aa=\pm}\|\jump{(0,\varphi_f)}_{\D_p}^\aa\|_{L^2(\Gamma)} 
 \right)  \\
& \qquad\times 
 \left(\sum_{n = n_\tau+1}^{n_{\tau'}}
 \dtn  \left( 
 \left( \xi^{(1),\alpha,n+1}_f \right)^2 +
 \left( \xi^{(2),\alpha,n+1}_f \right)^2
 \right) \right)^{1\over 2} 
 \times\left(\sum_{n = n_\tau+1}^{n_{\tau'}}
 \dtn \right)^{1\over 2} \\
 & \lsim  \( |\tau-\tau'|^{1\over 2} + \Delta t^{1\over 2}\). 
 \end{align*}
Notice that, since $\nabla^f_{\D_p}\varphi_f\to\nabla_\tau\bar\varphi_f$, $\Pi_{\D_p}^f\varphi_f\to\bar\varphi_f$ and $\jump{(0,\varphi_f)}_{\D_p}^\aa\to\jump{(0,\bar\varphi_f)}_\aa$ in their respective spaces, their norms are bounded, so that the maximum in the right-hand side above is well defined.
Using the estimate 
$$
\Big| \< [d_{f,\D_\bu} \Pi^f_{\D_p} s^\alpha_f](\tau) - [d_{f,\D_\bu} \Pi^f_{\D_p} s^\alpha_f](\tau'), \bar\varphi_f-\Pi^f_{\D_p}  \varphi_f \>_{L^2(\Gamma)}\Big| \lsim
\|d_{f,\D_\bu}\|_{L^\infty(0,T;L^2(\Gamma))} \| \bar\varphi_f-\Pi^f_{\D_p}  \varphi_f \|_{L^2(\Gamma)}. 
$$
and the a priori estimates of Lemma \ref{lemma_apriori} we deduce that
$$
\Big| \< [d_{f,\D_\bu} \Pi^f_{\D_p} s^\alpha_f](\tau) - [d_{f,\D_\bu} \Pi^f_{\D_p} s^\alpha_f](\tau'), \bar \varphi_f \>_{L^2(\Gamma)}\Big| \lsim
  \omega(|\tau-\tau'|) + \Delta t^{1\over 2} + \varpi_{\D_p}, 
$$
with $\lim_{h\rightarrow 0}\omega(h)=0$ and
$\varpi_{\D_p}= \| \bar\varphi_f-\Pi^f_{\D_p}  \varphi_f \|_{L^2(\Gamma)}$ a consistency error term such that
$\lim_{l\rightarrow + \infty}\varpi_{\D^l_p} = 0$.
It follows from the discontinuous Ascoli-Arzel\`a theorem \cite[Theorem C.11]{gdm} that (up to a subsequence) the sequence $d_{f,\D_\bu} \Pi^f_{\D_p} s^\alpha_f$ converges uniformly in time weakly in $L^2(\Gamma)$. Summing over $\alpha \in \{\rm nw,\rm w\}$, we also deduce the uniform-in-time $L^2(\Gamma)$-weak convergence of $d_{f,\D_\bu}$.
\end{proof}

\subsubsection{Strong convergence of $d_{f,\D_\bu}$, $d_{f,\D_\bu} \Pi^f_{\D_p} s_{f}^\alpha$, and $\Pi^f_{\D_p} s^\alpha_{f}$}

\begin{proposition}
  \label{prop_compactnessdfsf}
  Let $(\D_p^l)_{l\in\N}$, $(\D_\bu^l)_{l\in\N}$, $\{(t^l_n)_{n=0}^{N^l}\}_{l\in\N}$
be sequences of space time GDs assumed to satisfy the coercivity, consistency and compactness properties.
Let $\phi_{m,{\rm min}} >0$ and assume that, for each $l\in \N$, the gradient scheme \eqref{GD_hydro}--\eqref{GD_meca} has a solution $p^\alpha_l = (p^\a_{m,l}, p^\a_{f,l}) \in (X_{\D_p^l}^0)^{N^l+1}$, $\alpha\in \{\g,\l\}$, $\bu^l \in (X_{\D_\bu^l}^0)^{N^l+1}$ such that
  \begin{itemize}
\item[(i)] $d_{f,\D_\bu^l}(t,\x) \geq d_0(\x)$ for a.e.\ $(t,\x) \in (0,T)\times \Gamma$, 
\item[(ii)] $\phi_{\D^l}(t,\x) \geq \phi_{m,{\rm min}}$ for a.e.\ $(t,\x) \in (0,T)\times \Omega$.
  \end{itemize}
  Then, the sequence $(d_{f,\D^l_\bu})_{l\in\N}$ converges up to a subsequence in $L^\infty(0,T;L^p(\Gamma))$ for all $2 \leq p < 4$,
  and the sequences $(d_{f,\D^l_\bu} \Pi^f_{\D^l_p}s_f^{\alpha,l})_{l\in\N}$ and $(\Pi^f_{\D_p^l}s_f^{\alpha,l})_{l\in \N}$,  with $s_f^{\alpha,l}= S_f^\alpha(p_{c,f}^l)$ converge, up to a subsequence, in  $L^4(0,T;L^2(\Gamma))$. 
\end{proposition}
\begin{proof}
The proof is based on the same arguments employed in the proof of~\cite[Proposition 4.8]{bonaldi:hal-02549111}.
\end{proof}

\subsection{Convergence to a weak solution}\label{subsec:convergence}

\begin{proof}[Proof of Theorem~\ref{thm.conv}]
  The superscript $l$ will be dropped in the proof, and all convergences are up to appropriate subsequences. From Lemma \ref{lemma_apriori} and Proposition \ref{prop_compactnessdfsf}, there exist $\bar d_f\in L^\infty(0,T;L^4(\Gamma))$ 
and $\bar s_f^\alpha\in L^\infty((0,T)\times\Gamma)$ such that
\begin{equation}
  \label{conv_dfsf}
\begin{array}{lll}
  & d_{f,\D_\bu} \rightarrow \bar d_f & \mbox{in } L^\infty(0,T;L^p(\Gamma)),\, 2\leq p < 4, \\[1ex]
  & \Pi^f_{D_p} S^\alpha_f(p_{c,f}) \rightarrow \bar s^\alpha_f & \mbox{in } L^4(0,T;L^2(\Gamma)).
\end{array} 
\end{equation}
From Proposition \ref{prop_compactness_Sm}, there exist $\bar s_m^\alpha\in L^\infty((0,T)\times\Omega)$ and $\bar s_\aa^\alpha\in L^\infty((0,T)\times\Gamma)$ such that
\begin{equation}
  \label{conv_sm}
\begin{array}{lll}
  & \Pi^m_{D_p} S^\alpha_m(p_{c,m}) \rightarrow \bar s^\alpha_m & \mbox{in } L^2(0,T;L^2(\Omega)),\\
  & \mathbb{T}^\aa_{D_p} S^\alpha_\aa(p_{c,m}) \rightarrow \bar s^\alpha_\aa & \mbox{in } L^2(0,T;L^2(\Gamma)).
\end{array} 
\end{equation}
The identification of the limit \cite[Proposition~3.1]{BHMS2016}, resulting from the limit-conformity property, can easily be adapted to our definition of $V^0$, with weight $d_0^{\nf 3 2}$ and the use in the definition of limit-conformity of fracture flux functions that are compactly supported away from the tips. Using this lemma and the a priori estimates of Lemma \ref{lemma_apriori}, we obtain $\bar p^\alpha= (\bar p^\a_m,\bar p^\a_f)\in L^2(0,T;V^0)$ and  ${\bf g}^\alpha_f \in L^2(0,T;L^2(\Gamma)^{d-1})$, such that the following weak limits hold
\begin{equation}
  \label{conv_pnablap}
\left. 
\begin{array}{lll}
  & \Pi^m_{\D_p} p^\alpha_m \weakto \bar p^\alpha_m & \mbox{in } L^2(0,T;L^2(\Omega))\mbox{ weak}, \\
  & \Pi^f_{\D_p} p^\alpha_f \weakto \bar p^\alpha_f & \mbox{in } L^2(0,T;L^2(\Gamma))\mbox{ weak}, \\
  & \mathbb{T}^\aa_{\D_p} p^\a_m \weakto \gamma_\aa \bar p^\a_m & \mbox{in } L^2(0,T;L^2(\Gamma))\mbox{ weak}, \\
  & \nabla^m_{\D_p} p^\alpha_m \weakto  \nabla \bar p^\alpha_m & \mbox{in } L^2(0,T;L^2(\Omega)^d)\mbox{ weak}, \\
  & \jump{p^\a}^\aa_{\D_p}  \weakto \jump{\bar p^\a}_\aa & \mbox{in } L^2(0,T;L^2(\Gamma))\mbox{ weak}, \\
  & d_0^{\nf 3 2}\nabla^f_{\D_p} p^\alpha_f \weakto d_0^{\nf 3 2} \nabla_\tau \bar p^\alpha_f & \mbox{in } L^2(0,T;L^2(\Gamma)^{d-1})\mbox{ weak}, \\
  & d_{f,\D_\bu}^{\nf 3 2}\nabla^f_{\D_p} p^\alpha_f \weakto {\bf g}^\alpha_f & \mbox{in } L^2(0,T;L^2(\Gamma)^{d-1})\mbox{ weak}. 
\end{array}
\right.
\end{equation}
Let $\boldsymbol{\varphi} \in C_c^0((0,T)\times\Gamma)^{d-1}$ whose support is contained in $(0,T)\times K$, with $K$ compact set not containing the tips of $\Gamma$. We have
$$
\int_{0}^T\int_\Gamma d_{f,\D_\bu}^{\nf 3 2}\nabla^f_{\D_p} p^\alpha_f \cdot \boldsymbol{\varphi} ~\d\sigma(\x)\d t \rightarrow \int_{0}^T\int_\Gamma {\bf g}^\alpha_f \cdot \boldsymbol{\varphi}~ \d\sigma(\x)\d t. 
$$
On the other hand, it results from  \eqref{conv_pnablap} and the fact that $d_0$ is bounded away from $0$ on $K$ (because $d_0$ is continuous and does not vanish outside the tips of $\Gamma$) that $\nabla^f_{\D_p} p^\alpha_f \weakto \nabla_\tau \bar p^\alpha_f$ in $L^2(0,T;L^2(K)^{d-1})$. 
Combined with the convergence $d_{f,\D_\bu}^{\nf 3 2} \boldsymbol{\varphi} \rightarrow (\bar d_f)^{\nf 3 2}\boldsymbol{\varphi}$ in $L^{\infty}(0,T;L^2(\Gamma)^{d-1})$ given by \eqref{conv_dfsf}, we infer that 
\begin{equation*}
\int_{0}^T\int_\Gamma d_{f,\D_\bu}^{\nf 3 2}\nabla^f_{\D_p} p^\alpha_f \cdot \boldsymbol{\varphi}~ \d\sigma(\x)\d t
 \rightarrow 
\int_{0}^T\int_\Gamma (\bar d_f)^{\nf 3 2} \nabla_\tau \bar p^\alpha_f \cdot \boldsymbol{\varphi} ~\d\sigma(\x)\d t. 
\end{equation*}
This shows that ${\bf g}^\alpha_f = (\bar d_f)^{\nf 3 2} \nabla_\tau \bar p^\alpha_f$ on $(0,T)\times\Gamma$.

Combining the strong convergence of $\Pi^m_{D_p} S^\alpha_m(p_{c,m})= S^\alpha_m(\Pi^m_{D_p} p_{c,m})$, 
the weak convergence of  $\Pi^m_{D_p} p_{c,m}$, it results from the Minty trick (see, e.g.,~\cite[Lemma 2.6]{EGHM13}) that $\bar s_m^\alpha = S_m^\alpha(\bar p_{c,m})$ with $(\bar p_{c,m},\bar p_{c,f}) = \bar p_c = (\bar p_m^\g-\bar p_m^\l,\bar p_f^\g-\bar p_f^\l) = \bar p^\g-\bar p^\l$. Using the same arguments, we also have $\bar s_f^\alpha = S_f^\alpha(\bar p_{c,f})$ and $\bar s_\aa^\alpha = S_\aa^\alpha(\gamma_\aa \bar p_{c,m})$. 

From the a priori estimates of Lemma \ref{lemma_apriori} and the limit-conformity property of the sequence of GDs $(\D_\bu^l)_{l\in \N}$ (see \cite[Lemma A.3]{bonaldi:hal-02549111}), there exists $\bar\bu\in L^\infty(0,T;\U^0)$, such that
\begin{equation}
   \label{conv_unablau}
\left.
\begin{array}{lll}  
  & \Pi_{\D_\bu} \bu \weakto \bar \bu & \mbox{in } L^\infty(0,T;L^2(\Omega)^d) \mbox{ weak $\star$},\\
  & \bbeps_{\D_\bu} (\bu) \weakto \bbeps(\bar\bu) & \mbox{in } L^\infty(0,T;L^2(\Omega,\S_d(\R))) \mbox{ weak $\star$},\\
  & \div_{\D_\bu} \bu \weakto \div(\bar\bu) & \mbox{in } L^\infty(0,T;L^2(\Omega)) \mbox{ weak $\star$},\\
  & d_{f,\D_\bu} = -\jump{\bu}_{\D_\bu} \weakto -\jump{\bar\bu} & \mbox{in } L^\infty(0,T;L^2(\Gamma)) \mbox{ weak $\star$},
\end{array}
\right.
\end{equation}
from which we deduce that $\bar d_f = -\jump{\bar\bu}$ and that $\bbsig_{\D_u}(\bu)$ converges to $\bbsig(\bar \bu)$ in $L^\infty(0,T;L^2(\Omega,\S_d(\R)))$ weak $\star$. 

From the a priori estimates and the closure equations
\eqref{GD_closures}, there exist $\bar \phi_m \in L^\infty(0,T;L^2(\Omega)$ and $\bar p^E_m \in L^\infty(0,T;L^2(\Omega)$ such that  
\begin{equation}
   \label{conv_phipEm}
\left.
\begin{array}{lll}  
  & \phi_\D \weakto \bar \phi_m & \mbox{in } L^\infty(0,T;L^2(\Omega)) \mbox{ weak $\star$},\\
  & \Pi^m_{\D_p} p^E_m \weakto \bar p^E_m & \mbox{in } L^\infty(0,T;L^2(\Omega)) \mbox{ weak $\star$}. 
\end{array}
\right.
\end{equation}
Since $0\le  U_{\rm rt}(p) = \int_0^p q (S^{\g}_{\rm rt})'(q) \d q \leq 2 |p|$ for ${\rm rt} \in \{m,f\}$,
it results from the a priori estimates of Lemma \ref{lemma_apriori} that there exist $\bar p^E_f \in L^2(0,T;L^2(\Gamma))$, $\bar U_f \in L^2(0,T;L^2(\Gamma))$ and $\bar U_m \in L^2(0,T;L^2(\Omega))$ such that
\begin{equation}
   \label{conv_pEfUmf}
\left.
\begin{array}{lll}  
  & \Pi^f_{\D_p} p^E_f \weakto \bar p^E_f & \mbox{in } L^2(0,T;L^2(\Gamma)) \mbox{ weak},\\
  & \Pi^f_{\D_p} U_f(p_{c,f}) \weakto \bar U_f & \mbox{in } L^2(0,T;L^2(\Gamma)) \mbox{ weak},\\
  & \Pi^m_{\D_p} U_m(p_{c,m}) \weakto \bar U_m & \mbox{in } L^2(0,T;L^2(\Omega))\mbox{ weak}.   
\end{array}
\right.
\end{equation}
For ${\rm rt}\in \{\g,\l\}$, it is shown in \cite{DHM16}, following ideas from \cite{DE15}, that $U_{\rm rt}(p) = B_{\rm rt}(S^\g_{\rm rt}(p))$ where $s\in [0,1]\mapsto B_{\rm rt}(s)\in (-\infty,+\infty]$ is a convex lower semi-continuous function  with finite limits at $s=0$ and $s=1$ (note that $B_{\rm rt}$ is therefore continuous).
Since $\Pi^{m}_{\D_p} s^\g_{m}$ converges strongly in $L^2((0,T)\times\Omega)$  to $S^\g_m(\bar p_{c,m})$, it converges a.e.\ in $(0,T)\times\Omega$.
It  results that $B_{m}(\Pi^{m}_{\D_p} s^\g_{m})$  converges a.e.\ in $(0,T)\times\Omega$ to $B_{m}(S^\g_{m}(\bar p_{c,m}))$, and hence that 
$\bar U_{m} =  B_{m}(S^\g_{m}(\bar p_{c,m})) = U_{m}(\bar p_{c,m})$. Similarly, $\bar U_{f} =  B_{f}(S^\g_{f}(\bar p_{c,f})) = U_{f}(\bar p_{c,f})$. We deduce, using the strong convergences of the saturations and the weak convergences of the pressures, that
$$
\bar p^E_m = \sum_{\alpha\in \{\g,\l\}} \bar p^\alpha S^\alpha_m(\bar p_{c,m}) - U_m(\bar p_{c,m}) \quad \mbox{ and } \quad 
\bar p^E_f = \sum_{\alpha\in \{\g,\l\}} \bar p^\alpha S^\alpha_f(\bar p_{c,f}) - U_f(\bar p_{c,f}). 
$$
Using the estimate
$$
 \left|U_{\rm rt}(p_2)-U_{\rm rt}(p_1)\right|  = \left|\int^{p_2}_{p_1} q (S^{\g}_{\rm rt})'(q) \d q\right| \leq |p_2-p_1| + |p_2 S^{\g}_{\rm rt}(p_2) - p_1 S^{\g}_{\rm rt}(p_1)|,
$$
the Lipschitz property of $S^{\g}_{\rm rt}$, $\bar p^\alpha_0 = (\bar p^\a_{0,m},\bar p^\a_{0,f})  \in V^0\cap (L^\infty(\Omega)\times L^\infty(\Gamma))$, $\alpha\in \{\g,\l\}$, and the consistency of the sequence of GDs $(\D_p^l)_{l\in \N}$, we deduce that
  \begin{equation}
    \label{conv_pE0}
\left.
\begin{array}{lll}  
  & \Pi^m_{\D_p} p^{E,0}_m \rightarrow \bar p^{E,0}_m & \mbox{in } L^2(\Omega),\\
  & \Pi^f_{\D_p} p^{E,0}_f \rightarrow \bar p^{E,0}_f & \mbox{in } L^2(\Gamma). 
\end{array}
\right.
  \end{equation}
  Then, from \cite[Proposition A.4]{bonaldi:hal-02549111} it holds that 
  \begin{equation}
    \label{conv_u0}
\left.
\begin{array}{lll}  
  & \div_{\D_\bu}(\bu^0) \rightarrow \div(\bar\bu^0)  & \mbox{in } L^2(\Omega),\\
  & \jump{\bu^0}_{\D_\bu} \rightarrow \jump{\bar\bu^0} = - \bar d_f^0  & \mbox{in } L^2(\Gamma).    
\end{array}
\right.
\end{equation}
It results from  \eqref{conv_unablau}, \eqref{conv_phipEm}, \eqref{conv_pE0}, \eqref{conv_u0} and the definition of $\phi_\D$ that   
$$
\bar \phi_m = \bar \phi_m^0 +  b~\div(\bar\bu-\bar\bu^0) + \frac{1}{M} (\bar p^E_m -\bar p_m^{E,0}). 
$$

Let us now prove that the functions $\bar p^\alpha$, $\alpha\in \{\g,\l\}$, and $\bar\bu$ satisfy the variational formulation \eqref{eq_var_hydro}--\eqref{eq_var_meca} by passing to the limit in the gradient scheme \eqref{eq:GS}.

For $\theta\in C^\infty_c([0,T))$ and $\psi=(\psi_m,\psi_f)\in C_c^\infty(\Omega\setminus\overline\Gamma)\times C^\infty_c(\Gamma)$ let us set, with $P_{\D_p}\psi= (P^m_{\D_p}\psi_m, P^f_{\D_p}\psi_f)\in X^0_{\D_p}$ and $P^\nu_{\D_p}\psi_\nu$ realising the minimum of $\S_{\D_p^\nu}(\psi_\nu)$,
$$
\varphi =  (\varphi^1,\ldots,\varphi^{N}) \in (X^0_{\D_p})^{N} \mbox{ with } \varphi^i = (\varphi^i_m,\varphi^i_f)=\theta(t_{i-1}) (P_{\D_p}\psi).
$$
Let us set $\varphi_\nu = (\varphi^1_\nu,\ldots,\varphi^{N}_\nu), \nu\in\{m,f\}$. From the consistency properties of $(\D_p^l)_{l\in \N}$ with given $r> 8$, we deduce that
\begin{equation}
   \label{conv_thetapDpsi}
\left.
\begin{array}{lllll}
  & \Pi^m_{\D_p} P^m_{\D_p}\psi_m  \rightarrow \psi_m & \mbox{in } L^2(\Omega),\quad &
   \Pi^f_{\D_p}  P^f_{\D_p}\psi_f \rightarrow \psi_f & \mbox{in } L^2(\Gamma),\\  
  & \Pi^m_{\D_p}\varphi_m  \rightarrow \theta\psi_m & \mbox{in } L^\infty(0,T;L^2(\Omega)),\quad
  & \Pi^f_{\D_p} \varphi_f \rightarrow \theta\psi_f & \mbox{in } L^\infty(0,T;L^2(\Gamma)),\\
  & \nabla^m_{\D_p}\varphi_m  \rightarrow \theta \nabla \psi_m & \mbox{in } L^\infty(0,T;L^2(\Omega)^d),\quad
  & \nabla^f_{\D_p} \varphi_f \rightarrow \theta\nabla_\tau \psi_f & \mbox{in } L^\infty(0,T;L^r(\Gamma)^{d-1}),\\
   &  \mathbb{T}^\aa_{\D_p} \varphi_m \rightarrow \theta \gamma_\aa \psi_m & \mbox{in } L^\infty(0,T;L^2(\Gamma)),\quad
   & \jump{\varphi}^\aa_{\D_p} \rightarrow  \theta\jump{\psi}_\aa   & \mbox{in } L^\infty(0,T;L^2(\Gamma)).
\end{array}
\right.
\end{equation}

Setting
$$
\left.\begin{array}{lll}
 && T_1 = \dsp \int_0^T \int_\Omega \delta_t \(\phi_\D \Pi_{\D_p}^m s^\alpha_m \)\Pi_{\D_p}^m \varphi_m ~\d\x \d t \\[2ex]
 && T_2=  \dsp \int_0^T \int_\Omega \eta_m^\alpha(\Pi_{\D_p}^m s_m^\alpha) \K_m \nabla_{\D_p}^m p^\alpha_m \cdot   \nabla_{\D_p}^m \varphi_m  ~\d\x \d t\\[2ex]
 && T_3 = \dsp \int_0^T \int_\Gamma \delta_t \(d_{f,\D_\bu} \Pi_{\D_p}^f s^\alpha_f \)\Pi_{\D_p}^f \varphi_f ~\d\sigma(\x)\d t\\[2ex]
  && T_4= \dsp \int_0^T \int_\Gamma  \eta_f^\alpha(\Pi_{\D_p}^f s_f^\alpha) {d_{f,\D_\bu}^3 \over 12} \nabla_{\D_p}^f p^\alpha_f \cdot   \nabla_{\D_p}^f \varphi_f  ~\d\sigma(\x) \d t\\[2ex]
  && T_5 = \dsp \int_0^T \int_\Omega h_m^\alpha \Pi_{\D_p}^m \varphi_m ~d\x \d t + \int_0^T \int_\Gamma h_f^\alpha \Pi_{\D_p}^f \varphi_f ~\d\sigma(\x)\d t, \\[2ex]
 && T_1^\aa = \dsp \int_0^T \int_\Gamma \bar d_\aa\bar \phi_\aa \delta_t \( \mathbb{T}_{\D_p}^\aa s^\alpha_\aa \)\mathbb{T}_{\D_p}^\aa \varphi_m ~\d\sigma(\x)\d t\\[2ex]
  && T_2^\aa= \dsp \int_0^T \int_\Gamma   Q^\a_{f,\aa}  \jump{\varphi}_{\D_p}^\aa  ~\d\sigma(\x) \d t, 
  \end{array}\right.
$$
the gradient scheme variational formulation \eqref{GD_hydro} states that
$$
T_1 + T_2 + T_3 + T_4 + \sum_{\aa=\pm} (T_1^\aa + T^\aa_2)= T_5. 
$$
For $\omega\in C^\infty_c([0,T))$ and a smooth function ${\bf w}:  \Omega\setminus\overline\Gamma \rightarrow \R^d$ vanishing on $\partial\Omega$ and admitting finite limits on each side of $\Gamma$, let us set 
$$
\bv = (\bv^1,\ldots,\bv^{N}) \in (X^0_{\D_\bu})^{N} \mbox{ with } \bv^i = \omega(t_{i-1}) (P_{\D_\bu}{\bf w})
$$
where $P_{\D_\bu}\mathbf{w}$ realises the minimum in the definition \eqref{eq:def.SDu} of $\mathcal S_{\D_\bu}(\textbf{w})$.
From the consistency properties of $(\D_\bu^l)_{l\in \N}$, we deduce that 
\begin{equation}
   \label{conv_omegaw}
\left.
\begin{array}{lll}  
  & \Pi_{\D_\bu}\bv  \rightarrow \omega\psi & \mbox{in } L^\infty(0,T;L^2(\Omega)^d),\\ 
  & \bbeps_{\D_\bu}(\bv)  \rightarrow \omega \bbeps({\bf w}) & \mbox{in } L^\infty(0,T;L^2(\Omega,\S_d(\R))),\\
  & \jump{\bv}_{\D_\bu}  \rightarrow \omega \jump{{\bf w}} & \mbox{in } L^\infty(0,T;L^2(\Gamma)). 
\end{array}
\right.
\end{equation}
Setting 
$$
  \left.\begin{array}{lll}
    && T_6 =  \dsp \int_0^T \int_\Omega \( \bbsig_{\D_u}(\bu) : \bbeps_{\D_\bu}(\bv)
    -   b(\Pi_{\D_p}^m p_m^E)  \div_{\D_\bu}(\bv)\)d\x \d t,\\[2ex]
    && T_7 =   \dsp \int_0^T \int_\Gamma (\Pi_{\D_p}^f p_f^E)  \jump{\bv}_{\D_\bu} \d\sigma(\x)\d t,\\[2ex]
    && T_8 = \dsp \int_0^T \int_\Omega \mathbf{f} \cdot \Pi_{\D_\bu} \bv ~\d\x \d t.  
  \end{array}\right.
$$
the gradient scheme variational formulation \eqref{GD_meca} states that
$$
T_6 + T_7 = T_8. 
$$
Using a discrete integration by part \cite[Section D.1.7]{gdm}, we have $T_1 = T_{11} + T_{12}$ with 
\begin{align*}
& T_{11} = \dsp - \int_0^T \int_\Omega  \phi_\D (\Pi_{\D_p}^m s^\alpha_m) (\Pi_{\D_p}^m P^m_{\D_p} \psi_m) \theta'(t) ~\d\x \d t, \\
& T_{12} = - \int_\Omega  (\Pi^m_{\D_p} J_{\D_p}^m \bar \phi^0)  (\Pi_{\D_p}^m S^\alpha_m(I^m_{\D_p}\bar p^\alpha_{m,0}))  (\Pi_{\D_p}^m P^m_{\D_p}\psi) \theta(0) ~\d\x.  
\end{align*}
Using \eqref{conv_thetapDpsi} and \eqref{conv_phipEm}, and that $\Pi_{\D_p}^m s^\alpha_m \in [0,1]$ converges to $S^\alpha_m(\bar p_{c,m})$ a.e.\ in $(0,T)\times\Omega$ (this follows from \eqref{conv_sm}), it holds that
$$
T_{11} \rightarrow - \int_0^T \int_\Omega  \bar \phi_m S^\alpha_m(\bar p_{c,m})  \psi_m \theta'(t) ~\d\x \d t. 
$$  
Using  \eqref{conv_thetapDpsi}, that $\Pi^m_{\D_p} J_{\D_p}^m \bar \phi^0$ converges in $L^2(\Omega)$ to $\bar \phi^0$ and that 
$\Pi_{\D_p}^m S^\alpha_m(I^m_{\D_p}\bar p^\alpha_{m,0}) \in [0,1]$ converges a.e.\ in $\Omega$ to $S^\alpha_m(\bar p^\alpha_{m,0})$, we deduce that 
$$
T_{12} \rightarrow - \int_\Omega  \bar \phi^0  S^\alpha_m(\bar p^\alpha_{m,0})  \psi_m \theta(0) ~\d\x.  
$$
Writing $T_3 = T_{31} + T_{32}$ with 
\begin{align*}
& T_{31} = \dsp - \int_0^T \int_\Gamma  d_{f,\D_\bu} (\Pi_{\D_p}^f s^\alpha_f) (\Pi_{\D_p}^f P^f_{\D_p} \psi_f) \theta'(t) ~\d\sigma(\x) \d t, \\
& T_{32} = \int_\Gamma  \jump{\bu^0}_{\D_\bu}  (\Pi_{\D_p}^f S^\alpha_f( I^f_{\D_p} \bar p^\alpha_{f,0}))  (\Pi_{\D_p}^f P^f_{\D_p}\psi_f) \theta(0) ~\d\sigma(\x),   
\end{align*}
we obtain, using similar arguments and \eqref{conv_u0}, that
$$
T_{31} \rightarrow \dsp - \int_0^T \int_\Gamma  \bar d_f S^\alpha_f(\bar p_{c,f}) \psi_f \theta'(t) ~\d\sigma(\x) \d t,
$$
and
$$
T_{32} \rightarrow - \int_\Gamma  \bar d_f^0  S^\alpha_f(\bar p^\alpha_{f,0}) \psi_f \theta(0) ~\d\sigma(\x).  
$$
Writing $T_1^\aa = T_{11}^\aa + T_{12}^\aa$ with
\begin{align*}
& T_{11}^\aa = \dsp - \int_0^T \int_\Gamma  \bar d_\aa\bar \phi_\aa (\mathbb{T}_{\D_p}^\aa s^\alpha_\aa) (\mathbb{T}_{\D_p}^\aa P^m_{\D_p} \psi_m) \theta'(t) ~\d\sigma(\x) \d t, \\
& T_{12}^\aa = - \dsp \int_\Gamma  \bar d_\aa\bar \phi_\aa (\mathbb{T}_{\D_p}^\aa S^\alpha_\aa( I^m_{\D_p} \bar p^\alpha_{m,0}))  (\mathbb{T}_{\D_p}^\aa P^m_{\D_p}\psi_m) \theta(0) ~\d\sigma(\x),   
\end{align*}
we also obtain that
$$
T_{11}^\aa \rightarrow \dsp - \int_0^T \int_\Gamma  \bar d_\aa\bar \phi_\aa S^\alpha_\aa(\gamma_\aa \bar p_{c,m}) \gamma_\aa \psi_m \theta'(t) ~\d\sigma(\x) \d t,
$$
and
$$
T_{12}^\aa \rightarrow - \int_\Gamma  \bar d_\aa\bar \phi_\aa S^\alpha_\aa(\gamma_\aa \bar p^\alpha_{m,0}) \gamma_\aa \psi_m \theta(0) ~\d\sigma(\x).  
$$

Using that $0\le \eta^\alpha_m(\Pi_{\D_p}^m s_m^\alpha)\leq \eta^\alpha_{\rm m,max}$, the continuity of $\eta_m^\alpha$,
the convergence of $\Pi_{\D_p}^m s_m^\alpha$ a.e.\ in $(0,T)\times \Omega$ to $S^\alpha_m(\bar p_{c,m})$, 
\eqref{conv_pnablap} and \eqref{conv_thetapDpsi}, it holds that
$$
T_2 \rightarrow \dsp \int_0^T \int_\Omega \eta_m^\alpha(S_m^\alpha(\bar p_{c,m})) \K_m \nabla \bar p^\alpha_m \cdot   \theta \nabla \psi_m  ~\d\x \d t
$$
The convergence 
$$
T_4 \rightarrow   \int_0^T \int_\Gamma \eta^\a_f ( S^{\a}_f(\bar p_{c,f})) {\bar d_f^{\;3}\over 12}  \nabla_\tau \bar p^\a_f \cdot \theta \nabla_\tau \psi_f  \d\sigma(\x) \d t
$$
is established using $0\le \eta^\alpha_f(\Pi_{\D_p}^f s_f^\alpha)\leq \eta^\alpha_{\rm f,max}$, the continuity of $\eta_f^\alpha$, the convergence of $\Pi_{\D_p}^f s_f^\alpha$ a.e.\ in $(0,T)\times \Gamma$ to $S^\alpha_f(\bar p_{c,f})$, combined with the weak convergence of $d_{f,\D_\bu}^{\nf 3 2} \nabla_{\D_p}^f p^\alpha_f$ to $\bar d_{f}^{\nf 3 2} \nabla_\tau \bar p^\alpha_f$ in $L^2((0,T)\times\Gamma)^{d-1}$, the strong convergence of $d_{f,\D_\bu}^{\nf 3 2}$ to $\bar d_{f}^{\nf 3 2}$ in
$L^s((0,T)\times\Gamma)$ for all $2\leq s < {8\over 3}$ (resulting from \eqref{conv_dfsf}), and the strong convergence \eqref{conv_thetapDpsi} of $\nabla_{\D_p}^f \varphi_f$ to $\theta \nabla_\tau \psi_f$ in $L^\infty(0,T;L^r(\Gamma))$ with $r>8$. 

From \eqref{conv_thetapDpsi} we readily obtain the convergence
$$
T_5 \rightarrow \dsp \int_0^T \int_\O h_m^\a  ~\theta \psi_m  ~\d\x \d t + \int_0^T \int_\G h_f^\a ~\theta \psi_f~  \d\sigma(\x) \d t.
$$

The convergence
$$
T_2^\aa \rightarrow \int_0^T \int_\Gamma \bar Q^\a_{f,\aa}  \theta \jump{\psi}_\aa \d\x \d t, 
$$
results from the weak convergence of  $Q^\a_{f,\aa}$ to $\bar Q^\a_{f,\aa}$ in $L^2(0,T;L^2(\Gamma))$ combined with the strong convergence of $\jump{\varphi}^\aa_{\D_p}$ to $\theta\jump{\psi}_\aa$ in  $L^2(0,T;L^2(\Gamma))$. 

The following convergences of $T_6$, $T_7$, $T_8$
$$
\begin{aligned}
T_6 & \rightarrow \dsp \int_0^T \int_\O \( \bbsig(\bar \bu): \bbeps({\bf w})\omega - b \bar p_m^E \div( {\bf w}) \omega\) ~\d\x \d t,\\
T_7 & \rightarrow \int_0^T \int_\G \bar p_f^E ~\jump{ {\bf w}}\omega  ~\d\sigma(\x) \d t, \\
T_8 & \rightarrow \int_0^T \int_\Omega \mathbf{f}\cdot {\bf w} \omega~ \d\x \d t  
\end{aligned}
$$
classically result from the strong convergences \eqref{conv_omegaw} combined with the weak convergences \eqref{conv_unablau}.

Using the above limits in $T_1+T_2+T_3+T_4+ \sum_{\aa=\pm} (T_1^\aa + T^\aa_2)=T_5$ and $T_6+T_7=T_8$ concludes the proof that $\bar p^\alpha$, $\alpha\in \{\g,\l\}$, and $\bar\bu$ satisfy the variational formulation \eqref{eq_var_hydro}--\eqref{eq_var_meca}.
\end{proof}

\subsection{Identification of the limit interface fluxes}
\label{identification}

As mentioned in Remark \ref{rem:noQf}, the proof above does not identify the limit fluxes $\bar Q_{f,\aa}^\a$ of 
$$
Q^\a_{f,\aa}= T_f \left[ \eta^\a_\aa(S^\a_\aa(\mathbb{T}^\aa_{\D_p} p_{c,m}))  (\jump{p^\a}^\aa_{\D_p})^+  -
  \eta^\a_f(S^\a_f(\Pi^f_{\D_p}p_{c,f}))  (\jump{p^\a}_{\D_p}^\aa)^- \right]
$$
as $T_f \[ \eta_\aa(S^\a_\aa(\gamma_\aa \bar p_{c,m}))  \jump{\bar p^\a}_\aa^+  -  \eta_f(S^\a_f(\bar p_{c,f}))  \jump{\bar p^\a}_\aa^- \]$. The reason is that although the saturations $S^\a_\aa(\mathbb{T}^\aa_{\D_p} p_{c,m})$ and $S^\a_f(\Pi^f_{\D_p}p_{c,f})$ converge strongly, the pressure jumps $\jump{p^\a}^\aa_{\D_p}$ only converge weakly, which challenges the identification of the limits of their positive and negative parts (non-linear functions of the pressure jumps).

The expression of $Q^\a_{f,\aa}$ is however monotonic in terms of the pressure jumps $\jump{p^\a}^\aa_{\D_p}$, a feature that was used in \cite[Section 4.3]{DHM16} to identify the limit of these matrix--fracture fluxes in absence of mechanical deformations. The argument used there relies on a Minty technique (see, e.g., \cite[Section D.5]{gdm}). A key ingredient to this argument relies on being able to establish a energy equality for the limit of the approximations, which is done using the limit functions themselves as test functions in the weak equations \eqref{eq_var_hydro}--\eqref{eq_var_meca} they satisfy and using fine integrating-by-parts in time results from \cite{DE15}.

The caveat here is that space of test functions for \eqref{eq_var_hydro}, which is $C_c^\infty([0,T)\times\Omega\setminus\overline\Gamma)\times C_c^\infty([0,T)\times\Gamma)$, is not obviously dense in the space of trial functions $L^2(0,T;V^0)$, in which the limit pressures are found. Hence, it is not clear that we can indeed use these limit pressures as test functions in \eqref{eq_var_hydro}. The density issue comes from the fact that we would need to find smooth functions $(\varphi_{f,k}^\a)_{k\ge 1}$ such that $\bar d_f^{3/2}\nabla_\tau \varphi_{f,k}^\a\to \bar d_f^{3/2}\nabla_\tau p_{f}^\a$ in $L^2((0,T)\times \Gamma)$; in other words, we would like smooth functions to be dense in the weighted space $H^1_{\bar d_f}(\Gamma)$. 

Such a density result has been established in \cite{GWGM2015}, but under an additional assumption on the weight. Specifically:
$$
\begin{aligned}
&\bar d_f\mbox{ is bounded in time, smooth in space away from the fracture tips,}\\
&\mbox{and, near the tips, $\bar d_f(t,x,y)\approx x^{\frac12+\epsilon}f(y)$ with $\epsilon>0$}
\end{aligned}
$$
(above, $x$ is the distance to the tip, $y$ is the coordinate parallel to the tip, and $f$ is smooth). Under this assumption, the arguments of \cite{DHM16} can be reproduced and the limit fluxes $\bar Q_{f,\aa}^\a$ can be shown to satisfy the first equation in \eqref{eq_cont_closures}.

\section{Numerical experiments}
\label{sec:numerical.example}
The objective of this numerical section is to compare the discontinuous pressure poro-mechanical model investigated in this work with the continuous pressure poro-mechanical model presented in \cite{bonaldi:hal-02549111}. 
Two test cases are considered. The first one already described in  \cite{bonaldi:hal-02549111} considers the injection of gas in a cross-shaped fracture network coupled with the matrix domain initially liquid saturated. The second test case models the desaturation by suction at the interface between a ventilation tunnel and a low-permeability fractured porous medium. 

\begin{SCfigure}
\centering
\includegraphics[scale=1.2]{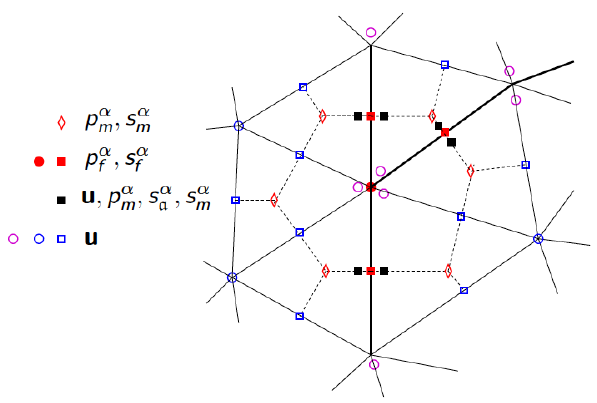}
\caption{{Example of admissible triangular mesh with three fracture edges in bold. The dot lines joining each cell center to the center of each of its edges are assumed orthogonal to the edge. The discrete unknowns of the discontinuous pressure model are presented for the two-phase flow and the mechanics. Note that the discontinuities of the pressures, of the saturations and of the displacement are captured at matrix fracture interfaces. Note also that additional nodal unknowns are defined at intersections of at least three fractures}.}
\label{mesh_tpfa_p2}
\end{SCfigure}

For both test cases, the flow part of system \eqref{eq_edp_hydromeca} is discretized in space by a Two-Point Flux Approximation (TPFA) cell-centered finite volume scheme with additional face unknowns at matrix fracture interfaces \cite{gem.aghili}. The mechanical part of \eqref{eq_edp_hydromeca} is discretized using second-order finite elements ($\P_2$) for the displacement field in the matrix \cite{daim.et.al,jeannin.et.al}, adding supplementary unknowns on the fracture faces to account for the discontinuities. The computational domain $\Omega$ is decomposed using \emph{admissible} triangular meshes for the TPFA scheme (cf.~\cite[Section~3.1.2]{finite.vol}) as illustrated in Figure \ref{mesh_tpfa_p2}. 

{For the TPFA scheme, the GD operators $\Pi_{\D_p}^m$ and $\Pi_{\D_p}^f$ are respectively cell-wise and fracture face-wise constant. It results that the porosity and the fracture aperture defined by the closure laws \eqref{GD_closures} will be projected in the matrix and fracture accumulation terms (cf.~first two terms in~\eqref{GD_hydro}) to cell-wise and face-wise constant spaces respectively. For simplicity, this face-wise constant projection of the fracture aperture is also used in the fracture conductivity.} 

Let $n\in\N^\star$ denote the time step index. The time stepping is adaptive, defined as
$$\dtn = \min \{ {\varrho}\dtnmun , \Delta t^{\max} \},$$
where $\dtzero = 0.001$ days is the initial time step, $\Delta t^{\max} = 10$ days in the first test case and 10 years in the second one, and ${\varrho} = 1.1$. At each time step, the flow unknowns
are computed by a Newton-Raphson algorithm. At each Newton-Raphson iteration, the Jacobian matrix is computed analytically and the linear system is solved using a GMRes iterative solver. The time step is reduced by a factor 2 whenever the Newton-Raphson algorithm does not converge within 50 iterations, with the stopping criteria defined by the relative residual norm lower than $10^{-5}$ or a maximum normalized variation of the primary unknowns lower than
$10^{-4}$. On the other hand, given the matrix and fracture equivalent pressures $p_m^E$ and $p_f^E$, the displacement field $\bu$ is computed using the direct solver MA48 (see~\cite{ma48}).
The coupled nonlinear system is then solved at each time step using a Newton-Krylov acceleration \cite{nitsol} of the fixed point algorithm which, for a given displacement field, solves the two-phase Darcy flow problem, 
then computes the new displacement field given the new equivalent pressures (see~\cite{ecmor} for more details). The stopping criterion is fixed to $10^{-5}$ on the relative displacement field increment.
{This Newton-Krylov algorithm is compared in \cite{ecmor} to the fixed stress algorithm \cite{KTJ11} extended to DFM models in \cite{GKW16}. It is shown to solve the robustness issue of fixed stress algorithms w.r.t. to small initial time steps in the case of incompressible fluids.}

For this section, we introduce the following notation for the \emph{total stress}:
\begin{equation}
\bbsig^T = \bbsig^0 + \bbsig(\bu) - b p_m^E \mathbb I,
\end{equation}
where $\bbsig^0$ is a possible \emph{pre-stress} state \cite[Section~4.2.4]{onate}.

\subsection{Gas injection in a cross-shaped fracture network}

The data set of the continuous pressure model is the one described in \cite{bonaldi:hal-02549111}. We recall it briefly here. We consider the square $\Omega = (0,L)^2$ lying in the $xy$-plane, with $L=100\,\text{m}$, containing a cross-shaped fracture network $\Gamma$ made up of four fractures, each one of length $\frac{L}{8}$ intersecting at $(\frac{L}{2},\frac{L}{2})$ and aligned with the coordinate axes. The matrix and fracture network have the following mobility laws: $\eta_{m}^\alpha(s^\alpha) = \frac{(s^\alpha)^2}{\mu^\alpha}$, $\eta_f^\alpha(s^\alpha) = \frac{s^\alpha}{\mu^\alpha}$, $\alpha \in \{\l,\g\}$, where $\mu^\l = 10^{-3}\,\rm{Pa{\cdot}s}$ and $\mu^\g = 1.851{\cdot}10^{-5}\,\rm{Pa{\cdot}s}$ are the dynamic viscosities of the wetting and non-wetting phases, respectively. Notice that $\eta_m^\alpha$ and $\eta_f^\alpha$ do not satisfy the assumptions of our analysis, as they are not bounded below by a strictly positive number; nevertheless, the results of the numerical experiments are not affected by this circumstance. The saturation--capillary pressure relation is Corey's law:
\begin{equation}\label{corey}
s_{\rm rt}^\g = S_{\rm rt}^\g (p_c) = \max\left(1 - \exp\left(-\frac{p_c}{R_{\rm rt}}\right),0\right),
\quad {\rm rt}\in \{m,f\},
\end{equation}
with $R_m=10^4\,\rm{Pa}$ and $R_f = 10\,\rm{Pa}$. The matrix is homogeneous and isotropic,
i.e.~$\mathbb K_m = \Lambda_m \mathbb I$, characterized by a permeability $\Lambda_m = 3{\cdot}10^{-15}\,\rm m^2$, an initial porosity $\phi_m^0 = 0.2$, effective Lam\'e parameters $\lambda = 833\,{\rm MPa}$, $\mu = 1250\,{\rm MPa}$, Biot's coefficient $b=1-\frac{K_{\rm dr}}{K_{\rm s}}\simeq0.81$, and Biot's modulus $M=18.4\,{\rm GPa}$. The pre-stress state is assumed null: $\bbsig^0=0$.  The domain is assumed to be clamped all over its boundary, i.e.~$\bf u = \bf 0$ on $(0,T)\times\del\Omega$; for the flows, we impose a wetting saturation $s_m^\l = 1$ on the upper side of the boundary $(0,T)\times((0,L)\times\{L\}) $, whereas the remaining part of the boundary is considered impervious (${\bf q}_m^\alpha\cdot\bf n = 0, \alpha\in\{\g,\l\}$).
%The system is subject to an initial pressure $p_0^\l = 10^5\,{\rm Pa}$ and an initial saturation $s_{0,{\rm rt}}^\l=1$, ${\rm rt}\in \{m,f\}$, for the wetting phase, which in turn results in an initial pressure $p_{0,{\rm rt}}^\g = p_{0,{\rm rt}}^\l + P_{c,{\rm rt}}(s_{0,{\rm rt}}^\g) =  10^5\,{\rm Pa}$ for the wetting phase. 
The system is subject to the initial conditions $p_{0}^\g = p_{0}^\l= 10^5\,{\rm Pa}$, which in turn results in an initial  saturation $s_{0,{\rm rt}}^\g=0$,  ${\rm rt}\in \{m,f\}$. 
The final time is set to $T=1000\,{\rm days}=8.64{\cdot}10^7\,\rm{s}$. The system is excited by the following source term, representing injection of non-wetting fluid at the center of the fracture network:
$$h_f^\g(t,\x) = \frac{g(\x)}{\dsp \int_\Gamma g(\x)\, \d\sigma(\x)}\frac{V_{\rm por}}{5T},\quad (t,\x)\in
(0,T)\times \Gamma,$$
where $V_{\rm por} = \dsp \int_\Omega \phi^0_m(\x)\,\d\x$ is the initial porous volume and 
$g(\x) = e^{-\beta |(\x-\x_0)/L|^2},\ \x_0 = (\frac{L}{2},\frac{L}{2})$, with $\beta=1000$ and $|{\cdot}|$ the Euclidean norm.
The remaining source terms $h^\l_f$ and $h^\alpha_m$, $\alpha\in\{\l,\g\}$, are all set to zero.

To define the discontinuous pressure model, we consider additionally the normal fracture transmissivity $T_f = 10^{-8}$ m. 

\begin{figure}
\begin{center}
  \includegraphics[width=.6\textwidth]{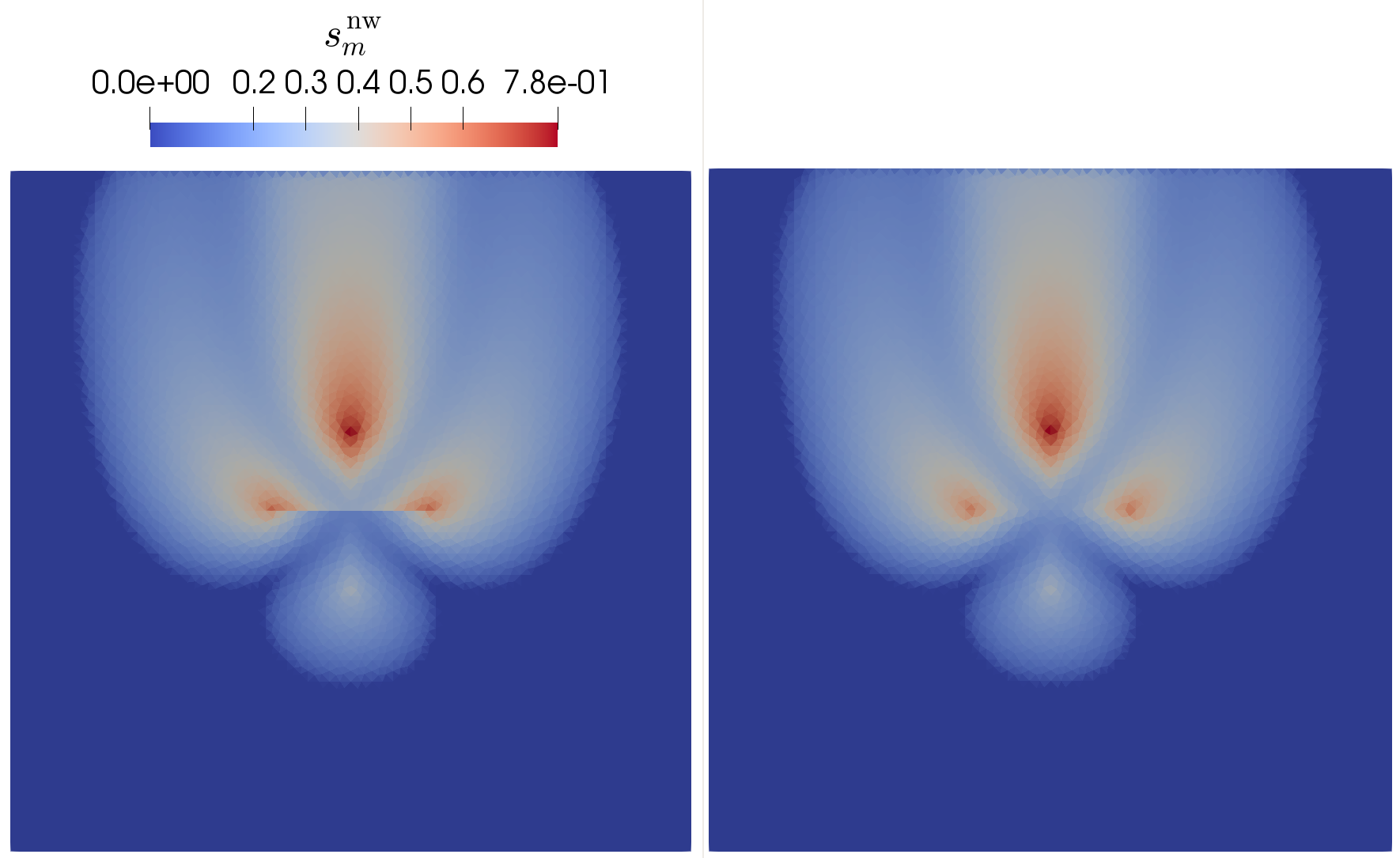}\\[5ex]
  \includegraphics[width=0.35\textwidth]{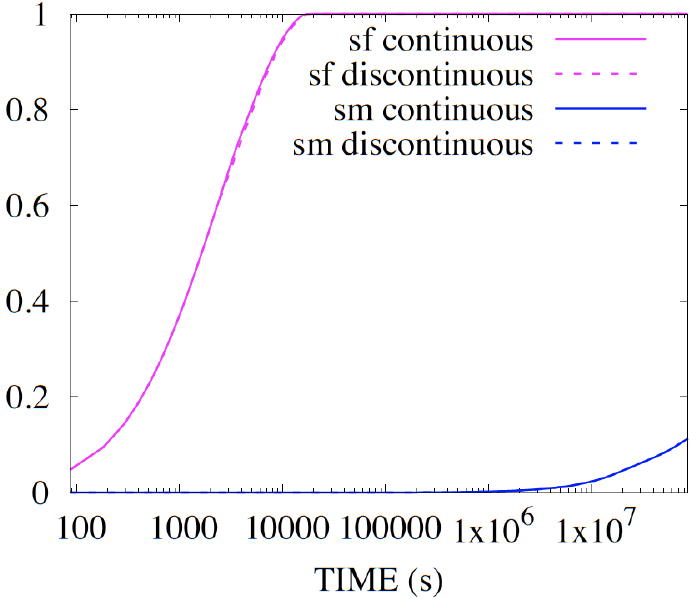}
  \includegraphics[width=0.4\textwidth]{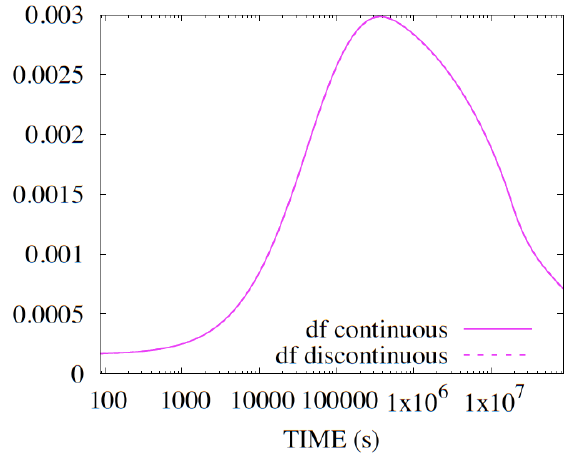}  
\caption{Top: matrix non-wetting phase saturation at final time obtained by the discontinuous pressure (left) and continuous pressure (right) models. Bottom left: mean non-wetting phase saturations in the matrix and in the fracture network as a function of time for both the continuous and discontinuous pressure models. Bottom right: mean  fracture apertures (m) as a function of time for both the continuous and discontinuous pressure models. }
\label{fig_gasinjection}
\end{center}
\end{figure}

\begin{figure}
\begin{center}
\includegraphics[width=.6\textwidth]{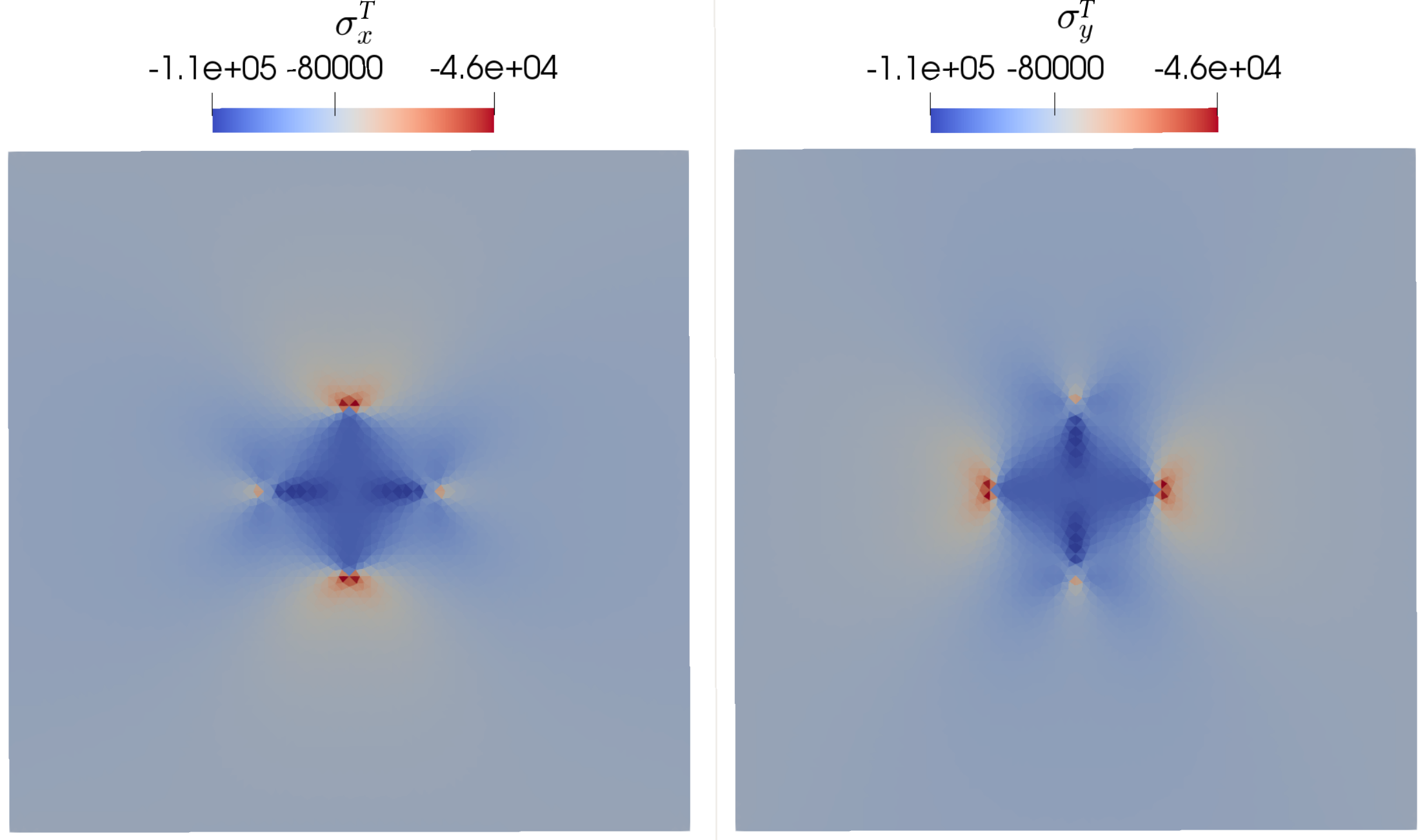}
\includegraphics[width=.3\textwidth]{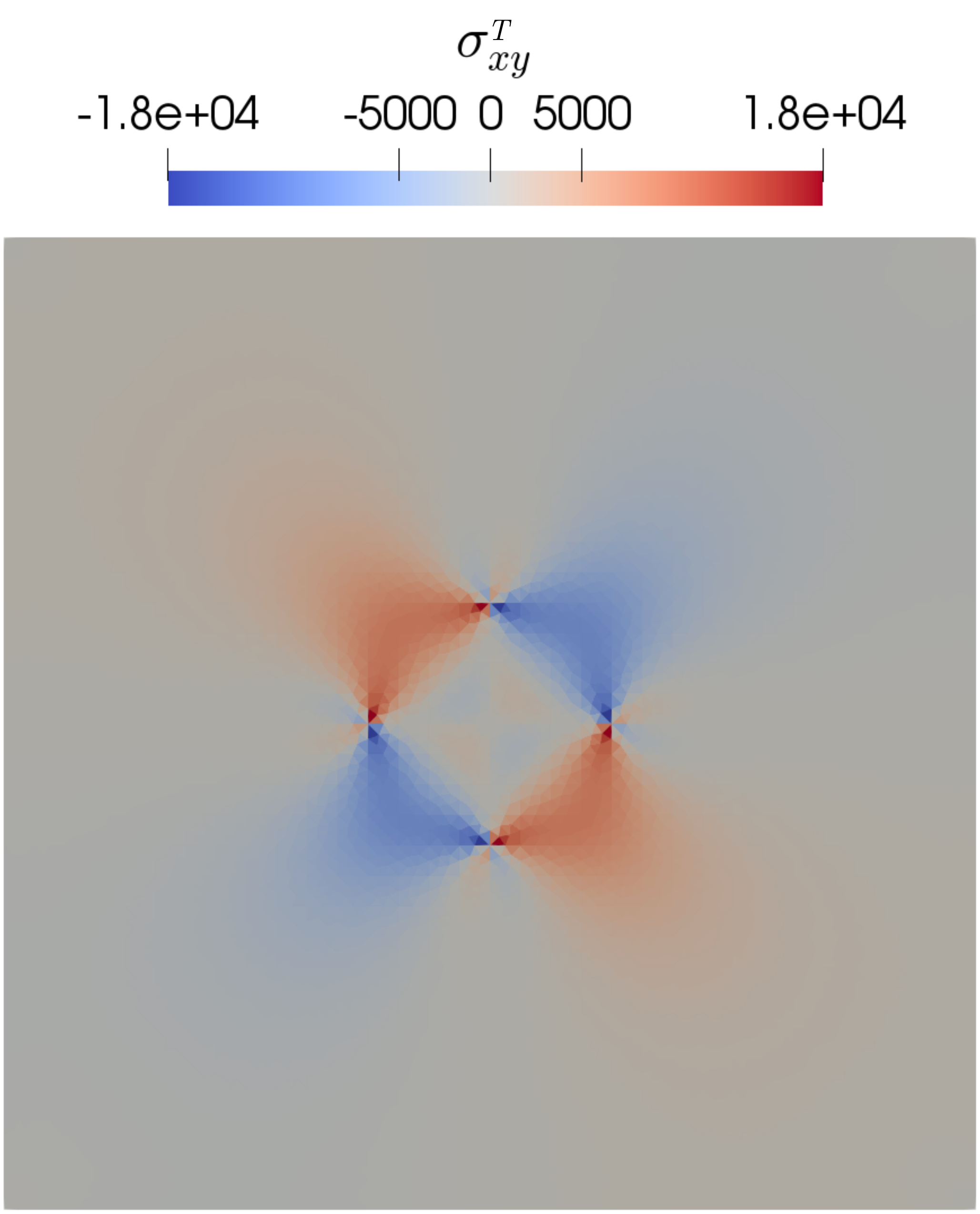}
\caption{Axial normal stresses $\sigma_{x}^T$ and $\sigma_{y}^T$ (Pa) and shear stress $\sigma_{xy}^T$ (Pa) at final time for the discontinuous pressure model.}
\label{plan_stress}
\end{center}
\end{figure}

From Figure \ref{fig_gasinjection}, it is clear that both the continuous and discontinuous pressure models provides roughly the same solutions. Nevertheless, the continuous pressure model provides a rather smoothed non-wetting phase saturation at matrix fracture interfaces while the discontinuous pressure model is more accurate as discussed in \cite{gem.aghili}. Also, Figure~\ref{plan_stress} shows the axial total stresses $\sigma_x^T$, $\sigma_y^T$, and the shear total stress $\sigma_{xy}^T$ in the matrix at the final time, for the discontinuous pressure model; from the mechanics viewpoint, in this case the difference between the two models is not so remarkable. As expected, stresses are concentrated in the neighborhood of fracture tips.
To conclude this subsection, we give an insight into the performance of our method in Table~\ref{perfs_plan}, where
%\newpage
\begin{itemize}
\item NbCells is the number of mesh cells,
\item N$_{\Delta t}$ is the number of successful time steps,
\item N$_{\text{Chops}}$ is the number of time step chops,
\item N$_{\text{Newton}}$ is the total number of Newton-Raphson iterations,
\item N$_{\text{GMRes}}$ is the total number of GMRes iterations,
\item N$_{\text{NK}}$ is the total number of Newton-Krylov iterations,
\item CPU[s] is the total computational time in seconds.
\end{itemize}
%%%
\begin{table}
\centering
{ 
  {
    \begin{tabular}{|c|c|c|}
\hline & Discontinuous pressure & Continuous pressure \\ \hline
NbCells & 14336 & 14336 \\ \hline
 {N}$_{\Delta t}$  & 187  & 187 \\ \hline
 N$_{\text{Chops}}$ & 0 & 0 \\ \hline
   N$_{\text{Newton}}$  & 2525 & 2294 \\ \hline
      N$_{\text{GMRes}}$  & 47733 & 32841 \\ \hline
    N$_{\text{NK}}$ & 1711 & 1618 \\ \hline
     CPU[s] & 307.4 & 246.7 \\ \hline
%     224 & 153  & 3298  & 13716 & 2351 \\ \hline
%     896 & 153  & 1857  & 9751 & 1308 \\ \hline
%     3584 & 153  & 1839  & 10792 & 1312 \\ \hline
%     14336 & 153  & 1976  & 14062 & 1324 \\ \hline
%     57344 & 153  & 2247  & 20330 & 1322 \\ \hline
%     229376 & 153  & 2772  & 34799 & 1329 \\ \hline
    \end{tabular}
    }
  }
    \caption{Performance of the method for the plane problem, in terms of the number of mesh elements, the number of successful time steps, the number of time step chops, the total number of Newton-Raphson iterations, the total number of GMRes iterations, the total number of Newton-Krylov iterations, and the total computational time.}
      \label{perfs_plan}
\end{table}

\subsection{Desaturation by suction of a low-permeability fractured porous medium}
\begin{figure}
\centering
\includegraphics[keepaspectratio=true,scale=.75]{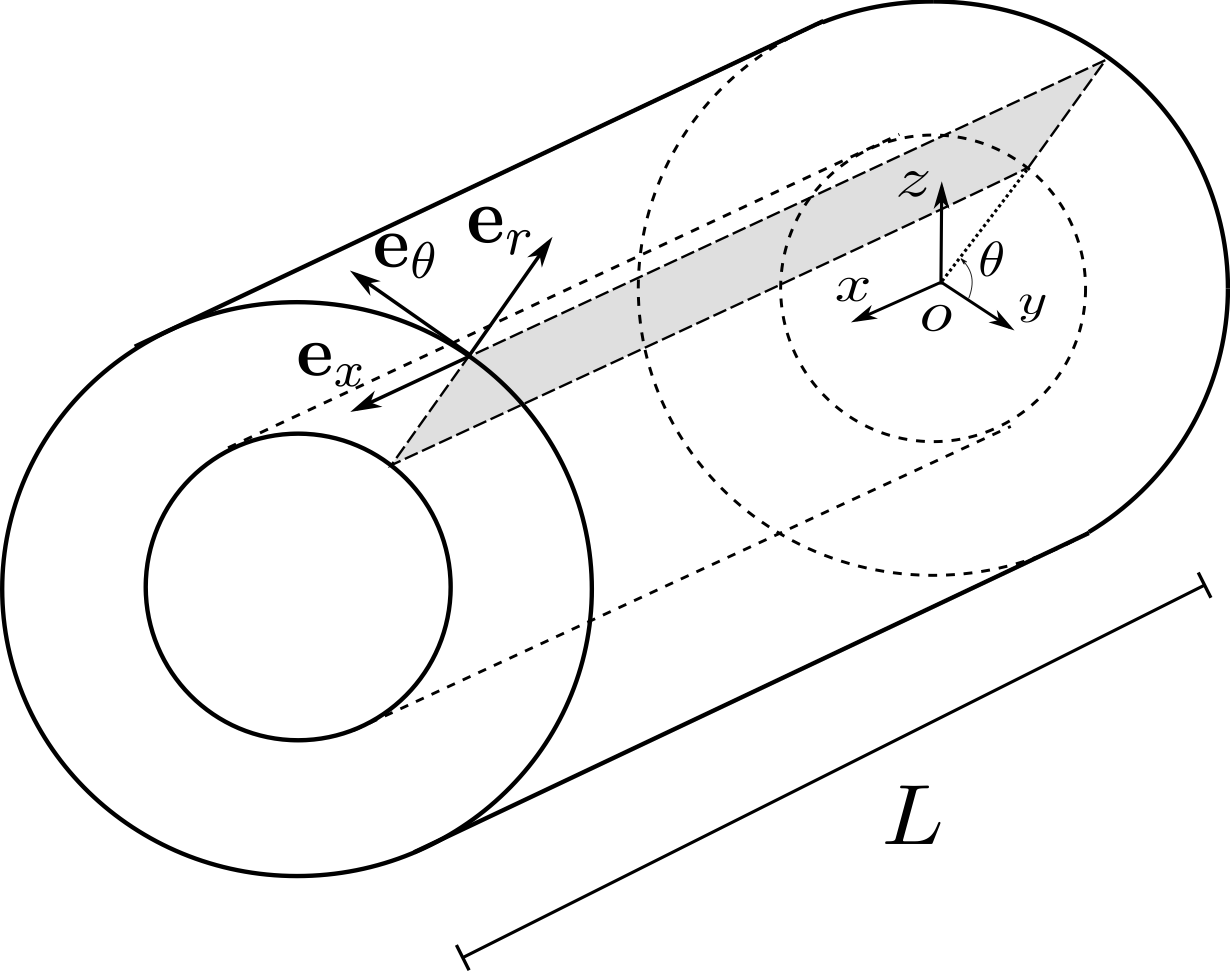}
\caption{Hollow cylinder of length $L$ and internal and external radii $R_{\rm int}$ and $R_{\rm ext}$, respectively. The diametral section is highlighted in gray, the fracture network is not shown for simplicity.}
\label{cilindro}
\end{figure}

\begin{figure}
\centering
\includegraphics[keepaspectratio=true,scale=.5]{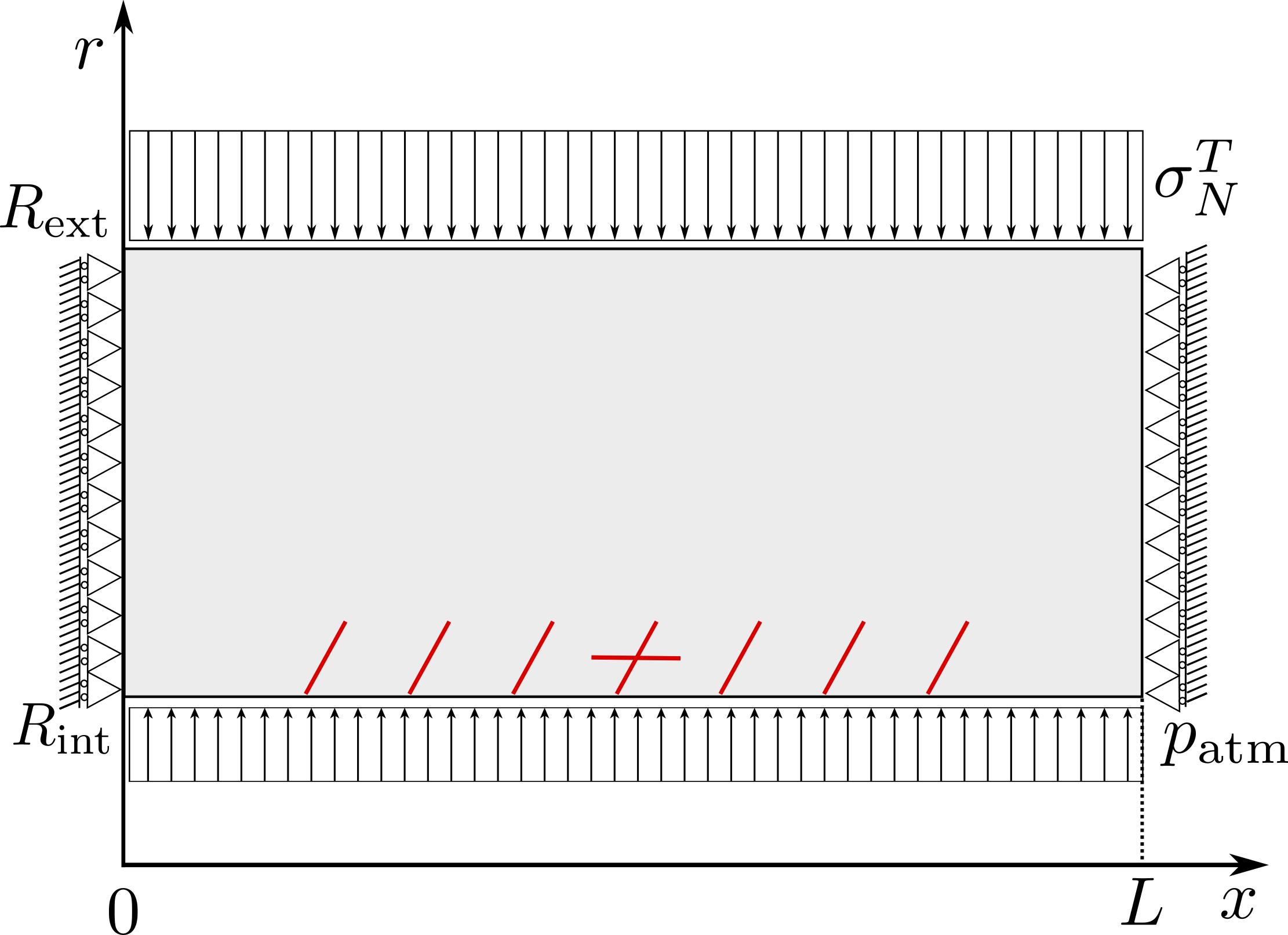}
\caption{Diametral section (not scaled) of the cylinder, with mechanical boundary conditions: zero normal displacement on the two sides $x=0$ et $x=L$ and uniform pressures $\sigma_N^T$ on the surface $r=R_{\rm ext}$ and $p_{\rm atm}$ on the surface $r=R_{\rm int}$. The fracture network is highlighted in red.}
\label{sezione}
\end{figure}

In this test case, we consider a hollow cylinder (Figure~\ref{cilindro}) made up of a low-permeability porous medium, containing an \emph{axisymmetric} fracture network, subject to axisymmetric loads -- uniform pressures exerted on the internal and external surfaces. Using cylindrical coordinates $(x,r,\theta)$, the problem can therefore be reduced to a two-dimensional formulation on the diametral section of the medium, shown in Figure~\ref{sezione} along with the fracture network, and the displacement field only consists of its axial and radial components:
$$
\begin{gathered}
\bu(x,r,\theta) = u_x(x,r) \mathbf e_x + u_r(x,r) \mathbf e_r(\theta),\quad
\mathbf e_r(\theta) = (\cos\theta) \mathbf e_y + (\sin\theta) \mathbf e_z,\\
0\le x \le L, \quad R_{\mathrm{int}} \le r \le R_{\mathrm{ext}},\quad 0 \le \theta \le 2\pi,
\end{gathered}
$$
where we have dropped time dependence for simplicity, and taken into account the system of cylindrical coordinates in Figure~\ref{cilindro}, denoting by $\mathbf e_x$ the axial unit vector, by $\mathbf e_r = \mathbf e_r(\theta)$ the radial unit vector, and by $\mathbf e_\theta$ the orthoradial unit vector. The final time for this simulation is set to $T=200$ years. The geometry is characterized by the following data set: length $L=10\,\text{m}$, internal and external radii $R_\mathrm{int} = 5\,\mathrm{m}$, $R_\mathrm{ext} = 35\,\mathrm{m}$; two consecutive fractures are spaced by 1.25\,m. The matrix is characterized by the Lamé parameters $\lambda = 1.5\,\text{GPa}$, $\mu = 2\,\text{GPa}$, by a permeability $\Lambda_m = 5\cdot10^{-20} \,\mathrm m^2$, Biot's coefficient and modulus $b=1$ and $M=1\,\text{GPa}$ respectively, and by an initial porosity $\phi_m = 0.15$. The normal transmissibility of fractures is  $T_f = 10^{-9}$ m and the initial fracture aperture is set to $10^{-2}$~m. 
The matrix relative permeabilities of 
the liquid and gas phases are defined 
by the following Van Genuchten laws:
\begin{eqnarray}
\label{eq_VGKrl}
k_{r,m}^\l(s^\l) = 
\left\{\begin{array}{r@{\,\,}c@{\,\,}ll}
&0 &\mbox{if}& s^\l < S_{lr},\\
&1 &\mbox{if}& s^\l > 1-S_{gr},\\
&\sqrt{\bar s^\l} \(1- (1-(\bar s^\l)^{1/q})^q\)^2 &\mbox{if}&  S_{lr}\leq s^\l \leq 1-S_{gr},  
\end{array}\right.
\end{eqnarray}
\begin{eqnarray}
\label{eq_VGKrg}
k_{r,m}^\g(s^\g) = 
\left\{\begin{array}{r@{\,\,}c@{\,\,}ll}
&0 &\mbox{if}& s^\g < S_{gr},\\
&1 &\mbox{if}& s^\g > 1-S_{lr},\\
&\sqrt{1-\bar s^\l} \(1- (\bar s^\l)^{1/q}\)^{2q} &\mbox{if}&  S_{gr}\leq s^\g \leq 1-S_{lr}, 
\end{array}\right.
\end{eqnarray}
with $$
\bar s^\l = {s^\l-S_{lr} \over 1-S_{lr}-S_{gr} },
$$
and the parameter $q=0.328$, the residual liquid and gas saturations 
$S_{lr}=0.40$ and $S_{gr}=0$; in the fractures, we take $k_{r,f}^\alpha(s) = s$ for both phases. The phase mobilities are then $\eta_m^\a(s^\a) = k_{r,m}^\a(s^\a)/\mu^\a$ and $\eta_f^\a(s^\a) = k_{r,f}^\a(s^\a)/\mu^\a$, $\a\in\{\l,\g\}$ both in the matrix and in the fractures, with the same viscosities as in the previous test case. Again, $\eta_m^\alpha$ and $\eta_f^\alpha$ are not bounded below by a strictly positive number, but this does not have an influence on the numerical results.
The saturation--capillary pressure relation is again Corey's law, as in~\eqref{corey}, with $R_m = 2{\cdot}10^8$\,Pa and $R_f = 10^2$\,Pa.
Moreover, the medium is supposed to be pre-stressed with the following pre-stress state:
\begin{equation*}
\bbsig^0 = \sigma_x^0 \, \mathbf e_x \otimes \mathbf e_x + \sigma^0_r \, {\mathbf e}_r \otimes \mathbf e_r +
\sigma^0_\theta \, \mathbf e_\theta \otimes \mathbf e_\theta,
\qquad \sigma_x^0 = 16\,\text{MPa}, \ \sigma_r^0 = \sigma_\theta^0 = 12\,\text{MPa}.
\end{equation*}
%(cf.~Figure~\ref{cilindro} for the representation of the unit vectors $\mathbf e_x$, $\mathbf e_r$, $\mathbf e_\theta$ corresponding to the chosen system of cylindrical coordinates).
%\paragraph{Jeu de données}
%\vspace{-.4cm}
%\begin{itemize}[leftmargin=*,itemsep=0pt]
%\item Géométrie: $L = 10\,\mathrm{m}$, $R_\mathrm{int} = 5\,\mathrm{m}$, $R_\mathrm{ext} = 35\,\mathrm{m}$, fractures espacées de $1.25\,\mathrm{m}$.
%\item Modules de Lamé de la matrice: $\lambda = 1.5\,\text{GPa}$, $\mu = 2\,\text{GPa}$.
%\item Précontraintes: 
%\item Porosité (initiale) matrice: $\phi_m = 0.15$. Porosité fractures: $\phi_f = 1$ (fractures ouvertes).
%\item Epaisseur initiale: $d_f = 10^{-2}$ m.
%\item Perméabilité matrice: $\Lambda_m = 5\cdot10^{-20} \,\mathrm m^2$, perméabilité fractures $d_f^2/12$.
%\item Transmissibilité normale des fractures:  
%\end{itemize}

Full saturation of the liquid phase is assumed at the initial state, both in the matrix and in the fracture network, with an initial uniform pressure $p_0^\l=p_0^\g = 4\,\mathrm{MPa}$.
%Cette valeur correspond à celle de la \emph{pression hydrostatique} 
%$$p_{\rm atm} + \rho_{\ell} g\zeta$$ dans la matrice, en notant $p_{\rm atm}$ la pression atmosphérique (1 bar), $g$ l'accélération de gravité ($\approx 10\,\rm m/s^2$) et $\zeta$ la profondeur du site de stockage ($\approx 400~\rm m$).
%\it%em Le déplacement initial, et par conséquent l'épaisseur initiale du réseau de fractures, est donné par les pressions initiales en résolvant la mécanique à l'instant $t=0$.
%\end{itemize}
%\paragraph{Conditions aux limites}
%\begin{itemize}[leftmargin=*]
%\item \emph{Mécanique}.

Concerning flow boundary conditions, the porous medium is assumed impervious (vanishing fluxes) on the lateral boundaries corresponding to $x=0$ and $x=L$.
On the inner surface $r=R_{\rm int}$, a given gas saturation is imposed: $s^\g_m = 0.35$ on the matrix side and $s^\g_f = 1 - 10^{-8}$ at fracture nodes, and atmospheric pressure $p_{\rm atm} = 10^5\,\text{Pa}$ everywhere. On the outer surface $r=R_{\rm ext}$, a liquid saturation $s_m^\l = 1$ and pressure $p^\l = 4\,\text{MPa}$ are imposed.

As for the mechanical boundary conditions, we impose a vanishing axial displacement $u_x$ on the lateral boundaries corresponding to $x=0$ and $x=L$. Moreover, on the same boundaries, the tangential stress is set to zero. %En revanche, remarquons que l'imposition de conditions aux limites sur le déplacement serait en fait «~interdite~» par la présence de la précontrainte. En effet, comme 
On the other hand, external surface loads $\mathbf g$ (uniform pressures) are applied on the inner and outer surfaces:
$$\mathbf g =\begin{cases}\begin{alignedat}{2} -\sigma_N^T \mathbf n, &\ \ \sigma_N^T > 0,&\ \ &\text{if }r=R_{\rm ext},\\
-p_{\rm atm}\mathbf n,  &\ \ p_{\rm atm} > 0,& \ \ &\text{if }r=R_{\rm int},\end{alignedat}\end{cases}$$
where $\mathbf n = \mathbf e_r$ for $r=R_{\rm ext}$ and $\mathbf n = -\mathbf e_r$ for $r=R_{\rm int}$.
%. En particulier, la valeur de $\sigma_N$ est donnée par la \emph{pression lithostatique} (poids des terres):
%$$
%\sigma_N = p_{\rm atm} + \rho_{\rm b} g\zeta;
%$$
%la densité de «~bulk~» (mélange roche-liquide) $\rho_{\rm b}$ est donnée par la moyenne pondérée avec la porosité de la roche $\phi_0$ de la manière suivante:
%$$\rho_{\rm b} = (1-\phi_0)\rho_{\rm r} + \phi_0\rho_{\ell},$$
%où $\rho_{\rm r}$ est la densité de la roche et $\rho_{\ell}$ la densité du liquide. 
We consider $\sigma_N^T = 10.95\,\text{MPa}$ as the numerical value for the uniform pressure on the outer surface.

\begin{figure}
\begin{center}
\includegraphics[scale=0.7]{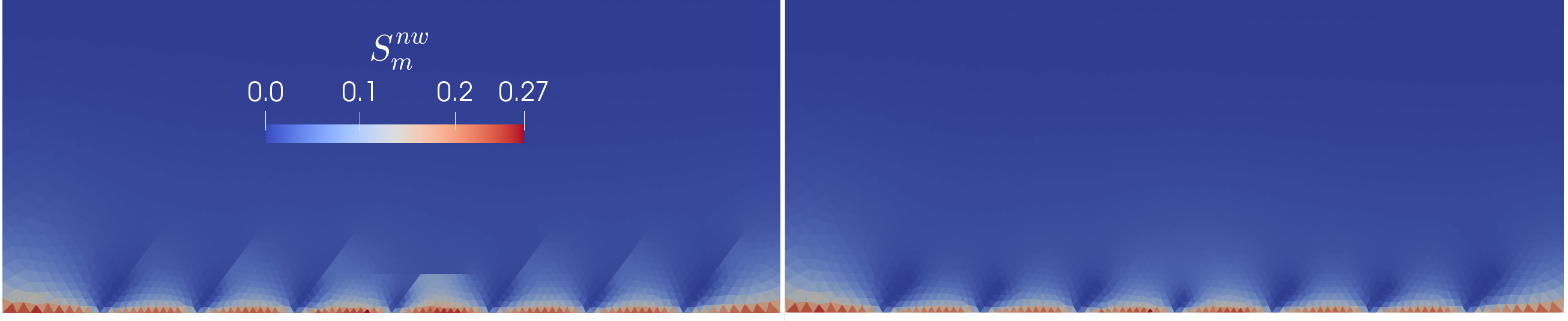}
%\\[5ex]
%\includegraphics[scale=.15]{TestAndra/Snw-cont-disc-coupe.png}
\caption{Zoom on the matrix non-wetting phase saturations at final time for the discontinuous (left) and continuous (right) pressure models.}
% Bottom: cut at $x=5.5$ m of the non-wetting phase saturations at final time for both models. }
\label{fig_Snw_tf_andra}
\end{center}
\end{figure}

\begin{figure}
\begin{center}
\includegraphics[scale=.7]{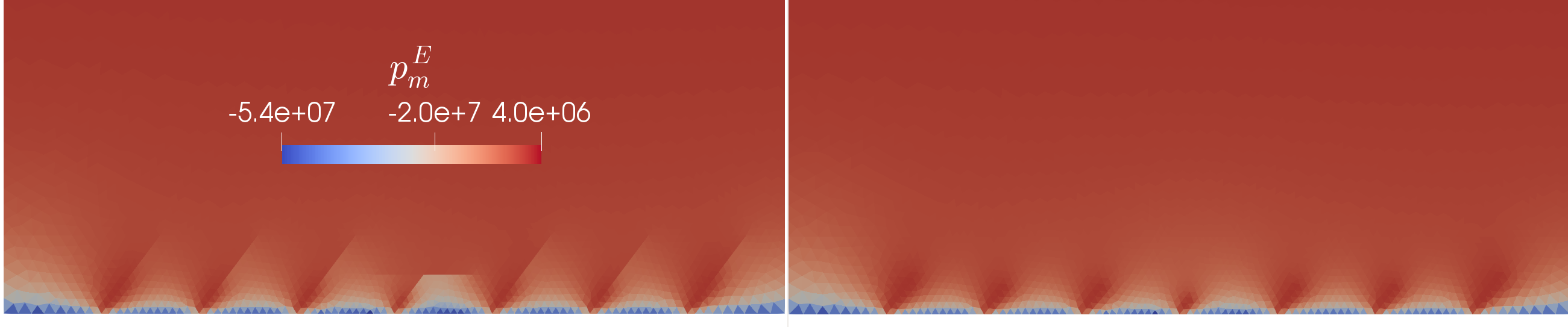}
\caption{Zoom on the matrix equivalent pressure (Pa) at final time for the discontinuous (left) and continuous (right) pressure models. }
\label{fig_PEm_df_tf_andra}
\end{center}
\end{figure}

\begin{figure}
\begin{center}
  \includegraphics[scale=.7]{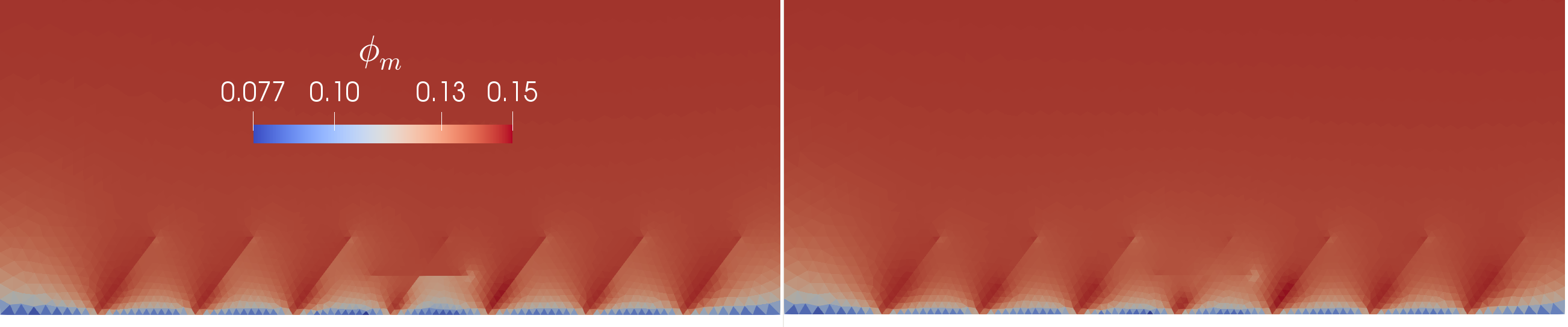}\\
\includegraphics[scale=.7]{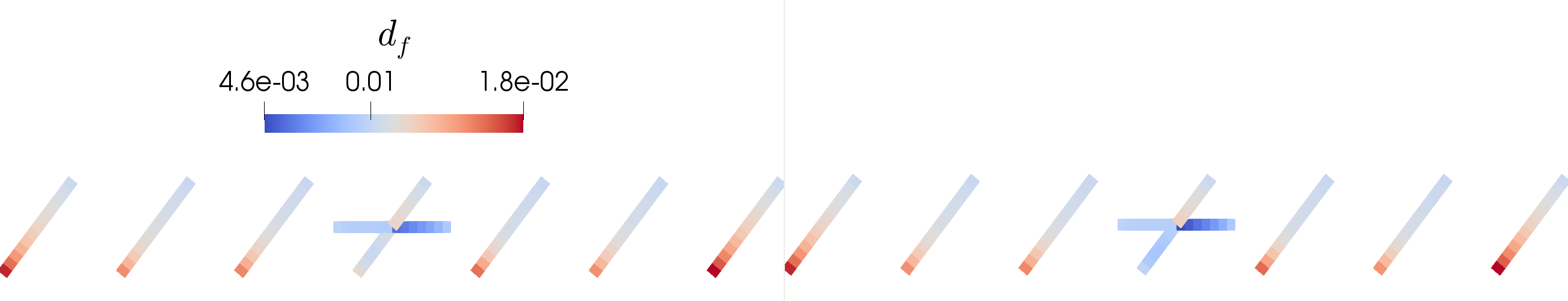}  
\caption{{Zoom on the matrix porosity $\phi_m$ and on the fracture aperture $d_f$ at final time for the discontinuous (left) and continuous (right) pressure models}. }
\label{fig_Phim_df_tf_andra}
\end{center}
\end{figure}

\begin{figure}
\begin{center}
\includegraphics[scale=.4]{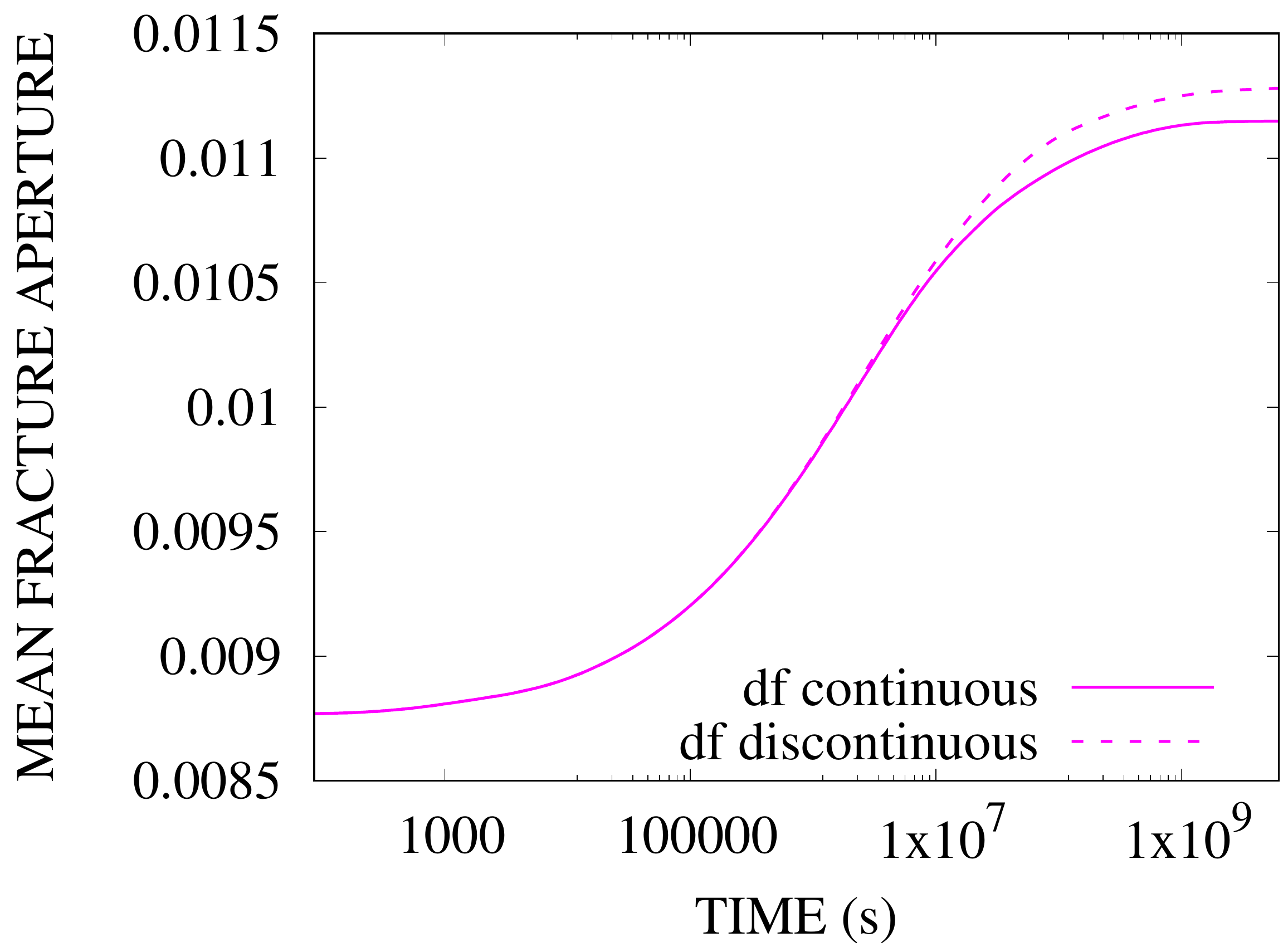}
\caption{Mean fracture aperture (\text{m}) as a function of time for both models.}
\label{fig_Snw_df_andra}
\end{center}
\end{figure}

\begin{figure}
\begin{center}
\includegraphics[scale=.1825]{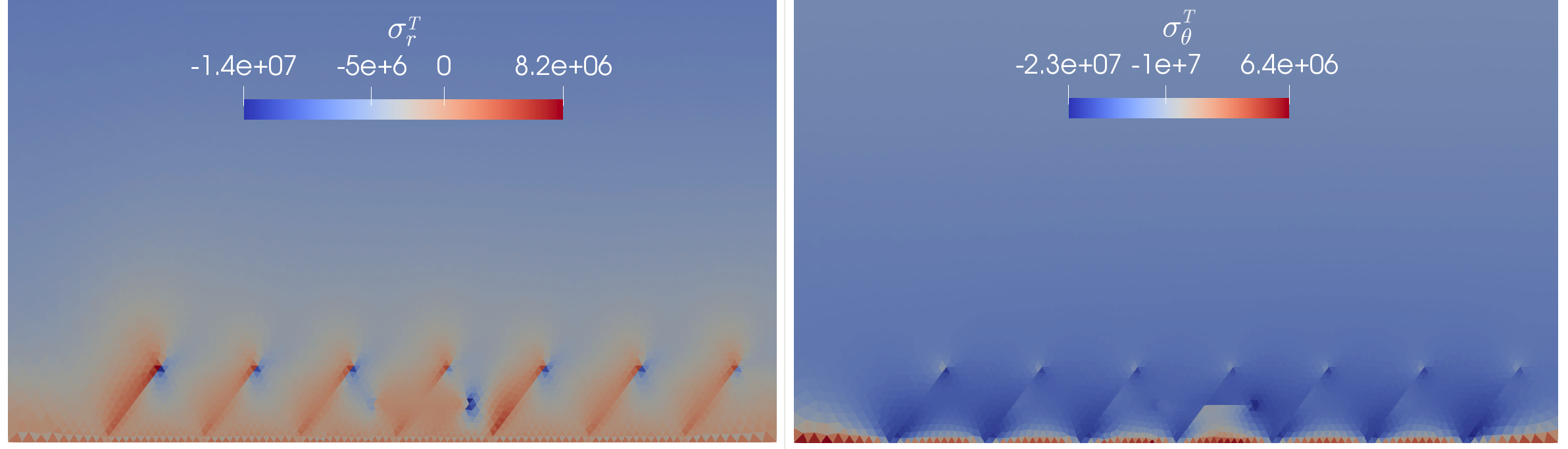}\\[5ex]
\includegraphics[scale=.1825]{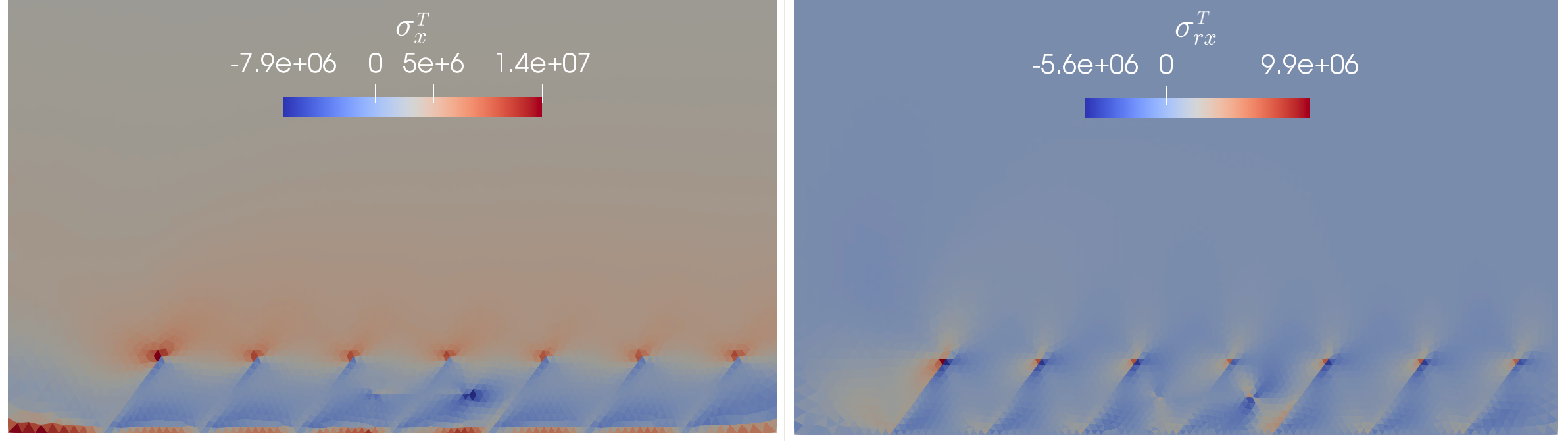}
\caption{Top: zoom on the radial and orthoradial total stresses $\sigma_r^T$ and $\sigma_\theta^T$ (Pa) at final time for the discontinuous pressure model. Bottom: zoom on the axial and shear total stresses $\sigma_x^T$ and $\sigma_{rx}^T$ (Pa) at final time for the discontinuous pressure model.}
\label{axi_stress}
\end{center}
\end{figure}

\begin{figure}
\begin{center}
\includegraphics[scale=1.25]{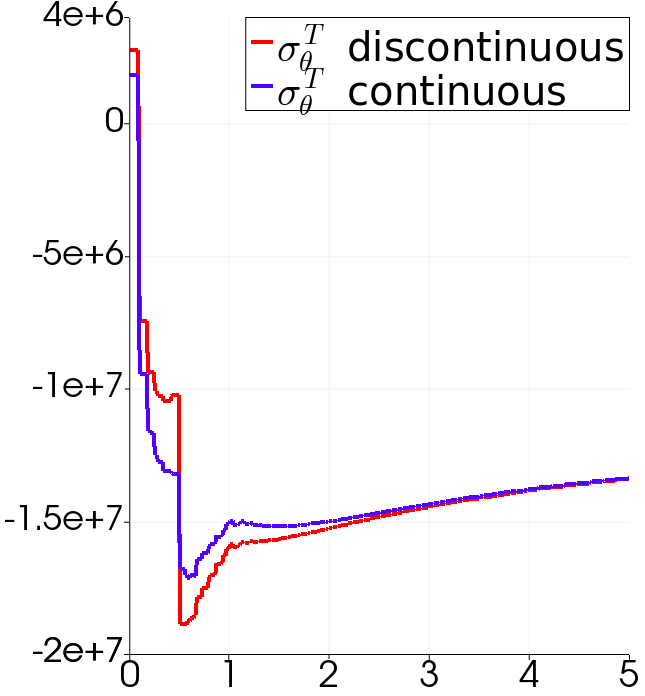}
\caption{Profile of the orthoradial total stress $\sigma_\theta^T$ (Pa) at final time along the line $x=5.5$ m as a function of the distance to the bottom boundary in m and for both models.}
\label{sigma_th_cont_disc}
\end{center}
\end{figure}

{As shown in Figures \ref{fig_Snw_tf_andra} and \ref{fig_PEm_df_tf_andra} strong capillary forces induce the desaturation of the matrix in the neighborhood of the inner surface combined with a high negative liquid pressure. As exhibited in Figures \ref{fig_Phim_df_tf_andra} and \ref{fig_Snw_df_andra} this negative liquid pressure triggers the contraction of the pores as well as the spreading of the fracture sides}. 
Figure~\ref{fig_Snw_tf_andra} also displays a comparison between the matrix non-wetting saturations obtained with the discontinuous and continuous pressure models at final time. It can be clearly seen that, unlike the continuous pressure model, the discontinuous pressure model is able to capture the barrier effect induced on the liquid phase by the fractures almost fully filled by the gas phase. This is particularly remarkable at the intersection of the horizontal and oblique fractures.
Figure~\ref{fig_PEm_df_tf_andra} shows a comparison of matrix equivalent pressures at final time obtained for the continuous and discontinuous pressure models; in the first case, discontinuities at the matrix-fracture interface can be clearly detected.
Figure~\ref{fig_Snw_df_andra} shows the time history of the average fracture aperture for the continuous and discontinuous pressure models, with significant differences induced by the equivalent pressures $p_m^E$ computed in the two models. Finally, in Figure~\ref{axi_stress} we display the radial, orthoradial, axial, and shear total stresses $\sigma_r^T$, $\sigma_\theta^T$, $\sigma_x^T$, and $\sigma_{rx}^T$ respectively in the matrix at the final time for the discontinuous pressure model. Again, stresses are concentrated in the neighborhood of fracture tips, as expected. {The arching effect is clearly visible by comparison of the radial and orthoradial stresses in the neighborhood of the inner surface. As expected, the radial stresses are transmitted across the horizontal fracture as opposed to the orthoradial stresses.} The comparison of the results given by the two models is shown in Figure~\ref{sigma_th_cont_disc}, where a different behavior in the orthoradial total stresses $\sigma_\theta^T$ given by the two models along the vertical line $x=5.5$ (intersecting the horizontal fracture) can be detected.

As in the previous subsection, we summarize also here the performance of our method in Table~\ref{perfs_axi}.
%%%
\begin{table}
\centering
{ 
  {
    \begin{tabular}{|c|c|c|}
\hline & Discontinuous pressure & Continuous pressure \\ \hline
NbCells & 28945 & 28945 \\ \hline
 {N}$_{\Delta t}$  & 169  & 176 \\ \hline
  N$_{\text{Chops}}$ & 0 & 1 \\ \hline
   N$_{\text{Newton}}$ & 2509  & 3758 \\ \hline
    N$_{\text{GMRes}}$  &  104329 & 122183 \\ \hline
    N$_{\text{NK}}$ & 693 & 721 \\ \hline
     CPU[s] & 1174.5 & 1413.6 \\ \hline
%     224 & 153  & 3298  & 13716 & 2351 \\ \hline
%     896 & 153  & 1857  & 9751 & 1308 \\ \hline
%     3584 & 153  & 1839  & 10792 & 1312 \\ \hline
%     14336 & 153  & 1976  & 14062 & 1324 \\ \hline
%     57344 & 153  & 2247  & 20330 & 1322 \\ \hline
%     229376 & 153  & 2772  & 34799 & 1329 \\ \hline
    \end{tabular}
    }
  }
    \caption{Performance of the method for the axisymmetric problem, in terms of the number of mesh elements, the number of successful time steps, the number of time step chops, the total number of Newton-Raphson iterations, the total number of GMRes iterations, the total number of Newton-Krylov iterations, and the total computational time.}
      \label{perfs_axi}
\end{table}

\section{Conclusions}

This work extends the gradient discretization and convergence analysis carried out in \cite{bonaldi:hal-02549111} to the case of hybrid-dimensional poro-mechanical models with discontinuous phase pressures at matrix fracture interfaces. 
The model considers a linear elastic mechanical model with open fractures coupled with a two-phase Darcy flow. The Poiseuille law is used for the tangential fracture conductivity and the dependence of the normal fracture transmissivity on the fracture aperture is frozen. The model accounts for a general network of planar fractures including immersed, non-immersed fractures and fracture intersections, and considers different rock types in the matrix and fracture network domains as well as at the matrix fracture interfaces. 

Two test cases were considered to compare the continuous pressure hybrid-dimensional poro-mechanical model investigated in \cite{bonaldi:hal-02549111} to the discontinuous pressure model studied in this work. The first test case simulates the gas injection in a cross-shaped fracture network immersed in a two-dimensional porous medium initially water saturated. The second test case is based on an axisymmetric DFM model and simulates the desaturation by suction at the interface between a ventilation tunnel and a Callovo-Oxfordian argilite fractured storage rock. In both cases, it is shown that the discontinuous pressure model provides a better accuracy at matrix fracture interfaces than the continuous pressure model and allows in particular to account for the barrier effect induced on the liquid phase by the gas filled fractures.

\begin{acknowledgements}
We are grateful to Andra and to the Australian Research Council's Discovery Projects (project DP170100605) funding scheme for partially supporting this work.
\end{acknowledgements}

%
%%%%%%%%%%%%%%%%%%
%
%%-----------------------------
%%      your bibliography
%%-----------------------------
%\bibliographystyle{abbrv}
%\clearpage
\bibliographystyle{plain}
\bibliography{Poromeca_disc}
\end{document}